\pdfoutput=1
\documentclass[11pt]{article}
\usepackage{preprint-layout}
\usepackage[T1]{fontenc}
\usepackage{lmodern}
\usepackage{mathrsfs}
\usepackage{microtype}
\usepackage[hidelinks]{hyperref}
\SetSymbolFont{stmry}{bold}{U}{stmry}{m}{n}
\graphicspath{{figures/}}

\theoremstyle{definitionbf}
\newtheorem{assumption}[thm]{Assumption}

\title{Energy and Localization of Nonlinear\\Multiscale Manifold Methods}
\author{Mats G. Larson\thanks{Department of Mathematics and Mathematical
Statistics, Ume{\aa}~University, SE-901\,87 Ume{\aa}, Sweden.
Email: \href{mailto:mats.larson@umu.se}{mats.larson@umu.se}.}
\and Anna Persson\thanks{Department of Information Technology,
Division of Scientific Computing, Uppsala University, Box 337,
75505 Uppsala, Sweden.
Email: \href{mailto:apersson@it.uu.se}{apersson@it.uu.se}.}}
\hypersetup{pdftitle={Energy and Localization of Nonlinear Multiscale Manifold Methods},
  pdfauthor={Mats G. Larson and Anna Persson}}
\date{26 November 2026}

\begin{document}
\maketitle

\begin{abstract}
We establish energy and localization estimates for a partition-of-unity
nonlinear multiscale manifold method for stationary elliptic problems.
The flux derives from a convex potential with rough spatial dependence:
measurable and uniformly bounded, with no periodicity or scale
separation. A local weak projection defines the fine-scale kernel, and
constrained nonlinear patch problems approximate global fine-scale
elimination. An exact secant equation and a kernel-preserving cutoff
yield exponential interior decay; weighted assembly transfers this decay
to the reconstruction. Direct energy comparison then bounds every global
minimizer on the localized manifold without a curvature or
chord-smallness condition. For uniformly positive and bounded Hessians,
the energy error decays exponentially in patch depth up to a polynomial
overlap factor, with the corresponding square-root rate in the gradient
norm. Under sufficient
oversampling, a
second estimate bounds the error at every coarse coordinate by its
reduced residual and a localization term. It controls all
tangent-stationary solutions without additional Hessian regularity.
Experiments for a saturated diffusion law examine projection choice,
contrast, load and mesh refinement.
\end{abstract}


\clearpage

\section{Introduction}
\label{sec:introduction}

Nonlinear elliptic problems with heterogeneous coefficients couple the
spatial structure of the medium to the state of the solution. A coarse
approximation must account for both the fine-scale variation of the
medium and its dependence on the solution. One approach is to determine
the unresolved scales from constrained nonlinear problems and to regard
their solutions as a reconstruction of the coarse state. Global
fine-scale elimination then defines an ideal nonlinear manifold, while
patch solves and partition-of-unity assembly define a localized manifold.
We ask how accurately the localized manifold represents the ideal
reconstruction and how localization affects energy minimization and
tangent-stationary coarse equations.

We prove such estimates for stationary elliptic equations whose flux is
the gradient of a convex potential. We treat potentials with
uniformly positive and bounded Hessian; the supplementary material
extends the analysis to potentials with \(p\)-structure, \(p\ge2\), whose
Hessian grows like \(|\nabla u|^{p-2}\). The coefficient may be rough: the potential is
only measurable in \(x\), with uniform bounds, and no periodicity or
scale separation is assumed. The reconstruction is the nonlinear
partition-of-unity construction of Larson and
Persson~\cite{larson-persson2026}, realized here with a fixed local weak
projection. Our analysis shows how the error introduced by restricting these
nonlinear problems to patches affects the final approximation. We
transfer the reconstruction error to solution accuracy in two ways:
energy comparison for global minimizers and a reduced-residual estimate
for computed coarse states.

\subsection{Related Work}
\label{subsec:related-work}

\paragraph{Variational Splitting and Localized Approximation.}
\label{par:literature-localization}
The approximation of unresolved scales by local variational problems
has an established role in numerical homogenization. The underlying
splitting into coarse and fine scales goes back to the variational
multiscale method~\cite{hughes1995,hughes-feijoo-mazzei-quincy1998}.
Hughes and Sangalli~\cite{hughes-sangalli2007} define the fine-scale
Green's function through a projection and show that its localization
depends on the choice of projection. Larson and
M\aa lqvist~\cite{larson-malqvist2007} approximate the fine component
by a sum of decoupled patch solutions and control localization through
a posteriori energy estimates. For linear elliptic equations,
M\aa lqvist and Peterseim~\cite{malqvist-peterseim2014} prove
exponential decay of localized correctors and show that logarithmic
oversampling preserves coarse convergence without solution regularity
or scale separation. Their construction provides the localization
mechanism underlying the localized orthogonal decomposition (LOD)
framework. The use of local integral degrees of freedom also has
classical precedents: Scott and Zhang~\cite{scott-zhang1990} construct
projections by local averaging that remain meaningful for nonsmooth
functions. Our local weak projection belongs to this family of
constructions; its role is to supply the kernel approximation and
support properties needed in the cutoff and assembly arguments.

\paragraph{Nonlinear Reconstructions and Enriched Coarse Spaces.}
\label{par:literature-nonlinear-reconstructions}
Nonlinear coarse-to-fine maps are not new to multiscale approximation.
Efendiev, Hou, and Ginting~\cite{efendiev-hou-ginting2004} formulate
nonlinear multiscale finite element methods using maps constructed from
local nonlinear problems, and analyze convergence and oversampling.
Efendiev, Galvis, Presho, and Zhou~\cite{efendiev-galvis-presho-zhou2014}
develop a different treatment of nonlinear monotone operators by
solving local nonlinear eigenvalue problems to enrich coarse spaces,
with applications to high-contrast approximation and domain decomposition.
Constrained nonlinear reconstruction also appears in the nonlinear
nonlocal multicontinua framework of Chung, Efendiev, Leung, and
Wheeler~\cite{chung-efendiev-leung-wheeler2018}. There, local solutions
on oversampled regions define nonlinear maps from prescribed coarse
averages to fine fields, which are used to formulate macroscopic
equations. Together, these works establish nonlinear local reconstruction and
oversampling as multiscale design principles that extend beyond any
one choice of coarse variables or assembly rule.

\paragraph{Linearized and Iterated Numerical Homogenization.}
\label{par:literature-linearization}
For semilinear elliptic equations, Henning, M\aa lqvist, and
Peterseim~\cite{henning-malqvist-peterseim2014} construct a linear
multiscale space from the heterogeneous diffusion operator and solve
the nonlinear coarse equation in that space. For strongly monotone
nonlinear diffusion, Verf\"urth~\cite{verfuerth2020} uses a linearization
to construct localized linear correctors and derives error estimates
that separate approximation and linearization errors. Khrais and
Verf\"urth~\cite{khrais-verfuerth2026} develop linearized LOD methods
for a class of nonmonotone quasilinear equations with solution-dependent
diffusion. The latter problem class differs from the uniformly convex
gradient-dependent energy problem considered here.

State dependence can also be addressed by updating the coarse spaces.
Liu, Chung, and Zhang~\cite{liu-chung-zhang2021} combine numerical
homogenization with iterative methods for monotone nonlinear equations,
including Leray--Lions problems, and develop residual regularization
and sparse updating of the coarse spaces. Thus an error associated
with a fixed linearization is not an unavoidable limitation of all
linearized approaches. Our analysis concerns a different object: the
exact nonlinear patch map and the nonlinear manifold obtained by its
partition-of-unity assembly. The patch equations retain the constitutive nonlinearity. A comparison
of computational efficiency with linearized or updated-space methods
requires a separate study.

\subsection{Positioning and Contributions}
\label{subsec:positioning-contributions}

The immediate starting point is the companion
paper~\cite{larson-persson2026}, which defines the ideal and localized
nonlinear manifolds, derives their tangent-space coarse equations, and
introduces optional learned approximations of restricted patch maps.
Its stationary perturbation analysis retains the localization defect
as an input and uses a local chord-smallness condition. The present
paper supplies a localization estimate for the patch-solved method
under additional energy structure. It covers uniformly bounded Hessians;
potentials with \(p\)-structure, \(p\ge2\), including the quartic
diffusion law, are treated in the supplementary material. It does not claim the scope of general monotone
formulations, and it does not cover degenerate laws such as the
\(p\)-Laplacian.

The distinction from the preceding nonlinear multiscale constructions
is the combination analyzed here: a fixed local projection kernel,
patch solutions approximating the complete eliminated nonlinear
correction with the same restricted load, and a conforming
partition-of-unity reconstruction. The analysis controls this specific
map and its reduced energy. The estimates connect patch-state decay to both energy minimization
and tangent-stationary coarse equations.

The main analytical contributions are as follows.
\begin{enumerate}
\item \textbf{Localization of complete nonlinear patch states.}
An exact secant equation and a kernel-preserving three-layer cutoff
provide interior decay. Partition-of-unity assembly transfers this
decay to a reconstruction bound with weighted patch overlap.

\item \textbf{Energy bounds for manifold minimizers.}
Direct comparison with the reconstruction at the exact coarse
coordinate gives energy and gradient error bounds without a
manifold-curvature or chord-smallness condition.

\item \textbf{Residual bounds for computed states.}
Ideal-kernel cancellation and stable coarse tests control the error
under sufficient oversampling. Extensions include inexact patch solves
and computable residual certificates. By a certificate we mean a bound
on the error of a computed state that can be evaluated from that state
alone, without the fine solution. The matrix formulas for the
certificates are given in the supplementary material.

\end{enumerate}
The supplementary material extends the energy and computed-state bounds
to potentials with \(p\)-structure, with a contraction factor that
depends only on scale-invariant local quantities of a weight determined
by the ideal state.

The numerical study examines projection choice, reference refinement,
contrast and load sensitivity. In the supplementary material,
certificates at saved outer-solver states distinguish patch evaluation
error from total discrete error, and experiments for the quartic law
examine the load dependence of the decay and of the weighted quantities,
coefficient channels, coarse and fine refinement, and frozen
linearizations.

\subsection{Main Result}
\label{subsec:main-result-overview}

Let \(V_H\subset V_h\subset H_0^1(\Omega)\) be nested piecewise-affine
finite element spaces in dimension \(d=2\) or \(3\), and let
\(\Pi_H:V_h\to V_H\) be the local weak projection defined in
Section~\ref{sec:local-projection}. Write \(V_f=\ker\Pi_H\).
The fine solution \(u_h\) minimizes the energy \(\mathcal E\) on
\(V_h\). For a coarse coordinate \(v_H\), the ideal reconstruction
\(M_h(v_H)\) minimizes this energy over \(v_H+V_f\).
The localized reconstruction \(M_h^\ell\) uses patches of \(\ell\)
coarse layers and is defined in Section~\ref{sec:reconstruction}.
Let \(u_{\rm ms}^\ell\) be any global minimizer of \(\mathcal E\)
on \(M_h^\ell(V_H)\), and let \(L_F\) denote the dual gradient norm
of the load on \(V_h\).

Under Assumption~\ref{ass:standing}, that is, the mesh assumptions and
uniform constitutive Hessian bounds \(0<m\le M<\infty\),
Theorem~\ref{thm:stationary-error} proves existence of such minimizers
and yields constants \(C,b>0\) such that
\begin{align}
 0\le\mathcal E(u_{\rm ms}^\ell)-\mathcal E(u_h)&\le C(1+\ell)^d e^{-b\ell}L_F^2 \quad \forall\ell\ge1 \label{eq:intro-energy-rate}\\
 \|\nabla(u_{\rm ms}^\ell-u_h)\|_\Omega&\le C(1+\ell)^{d/2}e^{-b\ell/2}L_F \quad \forall\ell\ge1 \label{eq:intro-gradient-rate}
\end{align}
where \(\|\cdot\|_\Omega\) is the \(L^2(\Omega)\) norm.
The constants depend on \(m,M\) and the stated mesh and local
projection constants. They are independent of state amplitude, load
size, and the fine-to-coarse mesh ratio. In particular, the estimates
are not asserted to be uniform as the ellipticity ratio \(M/m\)
increases.

The ideal manifold contains \(u_h\) exactly at its coordinate
\(\Pi_Hu_h\), that is, \(u_h=M_h(\Pi_Hu_h)\), because \(u_h\) minimizes
the energy over \(V_h\) and hence over \(\Pi_Hu_h+V_f\). Thus the estimates measure localization error relative
to the fine solution, with no separate coarse approximation term.
After scaling the coarse diameter \(H\) by a fixed domain length,
\(\ell\ge c\log(1/H)\), with \(c>2s/b\), gives gradient
localization error of order \(H^sL_F\) for any fixed \(s>0\).
This statement concerns accuracy; it is not a complexity estimate for
the nonlinear solves.

Theorem~\ref{thm:computed-state-error} gives a second result for
computed coarse coordinates. Let \(\rho_H(v_H)\) be the dual gradient
norm on \(V_H\) of the derivative of
\(\mathcal E\circ M_h^\ell\), defined in
\eqref{eq:reduced-residual-norm}. Under Assumption~\ref{ass:standing},
there is a fixed oversampling threshold \(\ell_*\) such that
\begin{equation}
 \|\nabla(M_h^\ell(v_H)-u_h)\|_\Omega\le3\frac{C_\Pi}{m}\rho_H(v_H)+C(1+\ell)^d e^{-b\ell/2}\frac{L_F}{m} \quad \forall\ell\ge\ell_* \label{eq:intro-stationary-rate}
\end{equation}
where \(C_\Pi\) bounds the projection in the gradient norm. Thus
every tangent-stationary solution is accurate under sufficient oversampling.
This result uses tangent stability and state localization; it assumes
neither global minimality nor a Lipschitz constitutive Hessian. Its
polynomial overlap factor is larger than in the energy-minimizer bound.
Proposition~\ref{prop:inexact-computed-state} extends this estimate to
inexact patch evaluations and stable coarse test maps.
Corollary~\ref{cor:continuous-error} separately combines fine
discretization and localization in a squared gradient error bound.

Section~\ref{sec:p-structure} of the supplementary material extends
\eqref{eq:intro-energy-rate} and the computed-state estimate
\eqref{eq:intro-stationary-rate} to potentials with \(p\)-structure, \(p\ge2\), whose Hessian grows like
\(|\nabla u|^{p-2}\). The ratio \(M/m\) is then replaced by
scale-invariant local quantities of a weight determined by the fine
solution, and Section~\ref{sec:p-constants} there studies their size.

\subsection{Organization}
\label{subsec:proof-structure}

Section~\ref{sec:energy-setting} specifies the problem and energy
distances. Section~\ref{sec:local-projection} verifies the local weak
projection and product interpolation estimates, and
Section~\ref{sec:reconstruction} defines the method and its tangent
maps. Section~\ref{sec:energy-comparison} proves the variational
comparison. Sections~\ref{sec:localization} and
\ref{sec:assembly-error} establish localization and the final error
bounds. Section~\ref{sec:stationary-residual} transfers the reduced residual
to a computed-state error bound under sufficient oversampling and extends
this bound to inexact patch solves and coarse test maps.
Section~\ref{sec:numerical-experiments} compares projections under fine
refinement, examines contrast and load sensitivity, and separates
comparator errors from computed-state and solver checks.
Section~\ref{sec:discussion} discusses the constitutive scope,
projection comparisons, and an extension to evolution problems.
Appendix~\ref{sec:buffered-cutoff} gives the complete cutoff construction.
The supplementary material gives matrix formulas and overlap-aware
assembly bounds for the residual certificates and evaluates them
numerically. It then replaces the bounded Hessian by \(p\)-structure,
proves the weighted energy and computed-state estimates, studies the size
of the weighted constants, and reports experiments for the quartic law.

\section{Energy Problem and Mesh Setting}
\label{sec:energy-setting}

The model is the weak form of a scalar nonlinear diffusion equation with
homogeneous Dirichlet data. For a volume load \(f\in L^2(\Omega)\),
its differential notation is
\begin{equation}
 -\operatorname{div}a(x,\nabla u)=f\text{ in }\Omega,\qquad u=0\text{ on }\partial\Omega \label{eq:model-diffusion}
\end{equation}
with the equation interpreted weakly when the spatial coefficients are
rough. The analysis is carried out relative to its fine Galerkin
solution and permits any linear load on the fine space, including
\(F(v)=\int_\Omega fv\,dx\). Exact integrals define all energies and
patch equations below; quadrature errors are outside the present scope.

\subsection{Spaces, Layers, and Boundary Patches}
\label{subsec:spaces-layers}

Let \(\Omega\subset\mathbb R^d\), with \(d=2\) or \(3\), be a bounded
connected polyhedral Lipschitz domain. Let \(\mathcal T_H\) be a conforming,
quasi-uniform, shape-regular simplicial mesh, and let \(\mathcal T_h\) be a
conforming shape-regular refinement. The continuous piecewise-affine
spaces with homogeneous Dirichlet conditions satisfy
\begin{equation}
 V_H\subset V_h\subset H_0^1(\Omega) \label{eq:nested-spaces}
\end{equation}
and the fine mesh need not be quasi-uniform. We use \(H\) for the largest
coarse diameter. All coarse diameters are comparable to \(H\). The
notation \(\|\cdot\|_D\) denotes the \(L^2(D)\) norm.

Two coarse elements are neighbors when their closures share a vertex.
For a set \(S\) of coarse elements, \(N^j(S)\) denotes its enlargement by
\(j\) steps in this graph, with \(N^0(S)=S\) and \(N(S)=N^1(S)\). We identify an element set
with the interior of its union when writing an integral or a support
condition. Differences of regions are taken elementwise. The number of
elements in a fixed-radius neighborhood is uniformly bounded. This layer
convention agrees with the closure-intersection neighborhoods in the
companion paper~\cite{larson-persson2026}.

We assume that the element patches used below are connected, have
diameter at most a constant times \(H\), and satisfy the following
Poincar\'e inequalities with a constant \(C_{\rm patch}\) that is uniform
over the mesh family. For an interior
patch \(\omega\), subtraction of its mean \(v_\omega\) gives
\begin{equation}
 \|v-v_\omega\|_\omega\le C_{\rm patch}H\|\nabla v\|_\omega \quad \forall v\in H^1(\omega) \label{eq:patch-poincare}
\end{equation}
with \(v_\omega=|\omega|^{-1}\int_\omega v\,dx\). If a coarse element
has a boundary vertex, its relevant enlargement can additionally be chosen
to contain a physical boundary facet and to satisfy
\begin{equation}
 \|v\|_\omega\le C_{\rm patch}H\|\nabla v\|_\omega \quad \forall v\in H^1(\omega),\quad v|_{\partial\Omega\cap\partial\omega}=0 \label{eq:boundary-poincare}
\end{equation}
where the facet used for the zero trace has size comparable to \(H\).
The boundary case is needed because the projection of
Section~\ref{sec:local-projection} does not reproduce constants near the
boundary; the proof of Lemma~\ref{lem:local-projection} uses the zero
trace there instead of a patch mean.
The enlargement radii and constants are uniform over the mesh family.
These are geometric assumptions, independent of the unknown solution.

\subsection{Potential and Constitutive Bounds}
\label{subsec:potential-bounds}

Let \(F\in V_h'\) be a linear functional. The potential
\(\psi:\Omega\times\mathbb R^d\to\mathbb R\) is measurable in its first
argument and twice continuously differentiable in its second argument
for almost every \(x\). We define the flux and fix its normalization by
\begin{equation}
 \psi(x,0)=0,\qquad a(x,0)=0,\qquad a(x,\xi)=D_\xi\psi(x,\xi) \label{eq:potential-normalization}
\end{equation}
and assume that the Hessian \(D_\xi a=D_\xi^2\psi\) of the potential is
uniformly bounded and positive definite: there are constants
\(0<m\le M<\infty\) such that
\begin{equation}
 m|\zeta|^2\le\zeta^T D_\xi a(x,\xi)\zeta\le M|\zeta|^2 \quad \forall\xi,\zeta\in\mathbb R^d \label{eq:hessian-bounds}
\end{equation}
for almost every \(x\). Being a Hessian, \(D_\xi a\) is symmetric, and
its eigenvalues lie in \([m,M]\). We write \(\chi=M/m\). No spatial derivatives of
\(a\), periodicity, or scale separation are assumed.

The normalization at zero simplifies the data bounds. If instead
\(a(x,0)=a_0(x)\in L^2(\Omega;\mathbb R^d)\) and
\(\psi(x,0)\) is integrable, subtract the affine part of the potential
and replace \(F(v)\) by \(F(v)-\int_\Omega a_0\cdot\nabla v\,dx\).
This leaves the Hessian and minimizers unchanged, up to an additive
energy constant, and puts the problem in the normalized form
\eqref{eq:potential-normalization}.

We measure the load in the dual of the gradient norm. The energy,
its nonlinear variational operator, and this load norm are
\begin{align}
 \mathcal E(v)&=\int_\Omega\psi(x,\nabla v)\,dx-F(v) \quad \forall v\in V_h \label{eq:energy-definition}\\
 \mathcal A(v)[w]&=\int_\Omega a(x,\nabla v)\cdot\nabla w\,dx \quad \forall v,w\in V_h \label{eq:operator-definition}\\
 L_F&=\|F\|_{-1,h}:=\sup_{0\ne w\in V_h}\frac{|F(w)|}{\|\nabla w\|_\Omega} \label{eq:load-norm}
\end{align}
with a zero supremum on a zero-dimensional space. Integration of
\eqref{eq:hessian-bounds} along a segment gives
\begin{align}
 [a(x,\xi)-a(x,\eta)]\cdot(\xi-\eta)&\ge m|\xi-\eta|^2 \label{eq:strong-monotonicity}\\
 |a(x,\xi)-a(x,\eta)|&\le M|\xi-\eta| \label{eq:flux-lipschitz}\\
 \frac m2|\xi|^2\le\psi(x,\xi)&\le\frac M2|\xi|^2 \label{eq:potential-quadratic-bounds}
\end{align}
for all \(\xi,\eta\). The symmetry of the Hessian is used to infer its
operator norm bound from \eqref{eq:hessian-bounds}.

\begin{prop}[Fine minimizer and basic bounds]
\label{prop:fine-minimizer}
The energy has a unique minimizer \(u_h\in V_h\), characterized by
\begin{equation}
 \mathcal A(u_h)[w]=F(w) \quad \forall w\in V_h \label{eq:fine-equation}
\end{equation}
and satisfying
\begin{equation}
 \|\nabla u_h\|_\Omega\le m^{-1}L_F,\qquad -\frac{L_F^2}{2m}\le\mathcal E(u_h)\le0 \label{eq:fine-data-bounds}
\end{equation}
under the normalization \eqref{eq:potential-normalization}.
\end{prop}

\begin{proof}
The lower bound \eqref{eq:potential-quadratic-bounds} gives
\begin{equation}
 \mathcal E(v)\ge\frac m2\|\nabla v\|_\Omega^2-L_F\|\nabla v\|_\Omega\ge-\frac{L_F^2}{2m} \label{eq:energy-coercivity}
\end{equation}
and coercivity follows. Finite dimensionality gives attainment, and strong
convexity gives uniqueness. Differentiation yields
\eqref{eq:fine-equation}. Testing with \(u_h\), using
\eqref{eq:strong-monotonicity} and \(a(x,0)=0\), proves the gradient
bound. The competitor zero gives the upper energy bound.
\end{proof}

\subsection{Energy Gap and Natural Distance}
\label{subsec:energy-distances}

We use two measures of error. The Bregman distance describes energy
comparison, while the natural distance is convenient for local decay
estimates. For \(v,w\in V_h\), define the potential Bregman distance by
\begin{equation}
 \mathcal B(v,w)=\int_\Omega\big[\psi(x,\nabla v)-\psi(x,\nabla w)-a(x,\nabla w)\cdot\nabla(v-w)\big]\,dx \label{eq:bregman-definition}
\end{equation}
and the squared natural distance on an element region \(D\) by
\begin{equation}
 \mathscr D_D(v,w)=\int_D[a(x,\nabla v)-a(x,\nabla w)]\cdot\nabla(v-w)\,dx \label{eq:natural-definition}
\end{equation}
where ``squared'' describes its quadratic scaling. No triangle inequality
for \(\mathscr D_D\) is assumed. Taylor's formula and
\eqref{eq:fine-equation} give
\begin{align}
 \mathcal E(v)-\mathcal E(u_h)&=\mathcal B(v,u_h) \label{eq:energy-bregman}\\
 \frac m2\|\nabla(v-u_h)\|_\Omega^2\le\mathcal B(v,u_h)&\le\frac M2\|\nabla(v-u_h)\|_\Omega^2 \label{eq:energy-gap-equivalence}\\
 m\|\nabla(v-w)\|_D^2\le\mathscr D_D(v,w)&\le M\|\nabla(v-w)\|_D^2 \label{eq:natural-equivalence}
\end{align}
and direct addition shows that
\begin{equation}
 \mathcal B(v,w)+\mathcal B(w,v)=\mathscr D_\Omega(v,w) \label{eq:bregman-symmetrization}
\end{equation}
for every pair. The energy gap and natural distance are different
quantities; in particular, neither Bregman term is generally half their sum.

\begin{example}[A saturated constitutive law]
\label{example:saturated-law}
Let \(0<a_-\le a_\varepsilon(x)\le a_+\) and \(\alpha,\beta>0\). The law
\begin{equation}
 a(x,\xi)=a_\varepsilon(x)\left(1+\frac{\alpha|\xi|^2}{1+\beta|\xi|^2}\right)\xi \label{eq:saturated-law}
\end{equation}
is generated by the explicit potential
\begin{equation}
 \psi(x,\xi)=\frac{a_\varepsilon(x)}2\left[\left(1+\frac\alpha\beta\right)|\xi|^2-\frac\alpha{\beta^2}\log(1+\beta|\xi|^2)\right] \label{eq:saturated-potential}
\end{equation}
which vanishes at zero. Here \(a_\varepsilon\) is any measurable
coefficient with the displayed bounds; the subscript indicates possible
fine spatial variation and imposes no periodic structure. Its tangential
Jacobian multiplier lies between \(1\) and \(1+\alpha/\beta\). Writing
\(t=\beta|\xi|^2\), the radial multiplier is
\begin{equation}
 1+\frac\alpha\beta\frac{t(t+3)}{(1+t)^2}\le1+\frac{9\alpha}{8\beta} \quad \forall t\ge0 \label{eq:saturated-radial-bound}
\end{equation}
because the rational factor has its maximum \(9/8\) at \(t=3\).
Thus \eqref{eq:hessian-bounds} holds with \(m=a_-\) and
\(M=a_+(1+9\alpha/(8\beta))\). The localization proof below uses these
bounds, not the particular formula \eqref{eq:saturated-law}.
\end{example}

\subsection{Assumptions}
\label{subsec:standing-hypotheses}

The results below use the following conditions. The detailed definitions
are those given above and in the indicated sections.

\begin{assumption}
\label{ass:standing}
\leavevmode
\begin{enumerate}[label=(A\arabic*),ref=(A\arabic*)]
\item \label{ass:potential}
The potential is measurable in space and \(C^2\) in its gradient
argument, satisfies \eqref{eq:potential-normalization}, and has the
uniform symmetric Hessian bounds \eqref{eq:hessian-bounds}. The load
is any \(F\in V_h'\).
\item \label{ass:mesh}
The meshes and domain satisfy Section~\ref{subsec:spaces-layers},
including uniform fixed-patch Poincar\'e inequalities at the physical
boundary. Coarse elements are quasi-uniform; the fine refinement is
shape-regular and may be graded.
\item \label{ass:projection}
The coarse projection is the fixed local weak operator of
Section~\ref{sec:local-projection}. Its reproduction, locality, gradient
stability, and kernel approximation are proved there. Both corrected
cutoffs and assembly use \(I-\Pi_H\).
\item \label{ass:patches}
Patches are fixed as the coarse input varies. Their zero-extended
corrections solve the exact same-load equations in the global fine
kernel. Assembly uses the full partition of unity, including
boundary hats, as defined in Section~\ref{sec:reconstruction}.
\end{enumerate}
\end{assumption}

Section~\ref{sec:p-structure} of the supplementary material replaces
\ref{ass:potential} by the \(p\)-structure condition
\eqref{eq:p-structure} and keeps \ref{ass:mesh}--\ref{ass:patches}.

\section{Local Weak Projection}
\label{sec:local-projection}

The projection must preserve coarse functions and provide local
control of its kernel. We construct it from elementwise integral
degrees of freedom, then verify the estimates needed for the cutoff
and partition-of-unity assembly.

\subsection{Local Dual Degrees of Freedom}
\label{subsec:local-duals}

Let \(\mathcal N_H^{\rm int}\) be the interior coarse vertices and
\(\lambda_z\) their standard continuous coarse nodal hats. For every
\(z\in\mathcal N_H^{\rm int}\), select one coarse simplex \(K_z\)
containing \(z\). The selection is fixed independently of the state and
load. For a vertex \(z\) of a simplex \(K\), let \(\lambda_z^K\)
denote its barycentric coordinate and define
\begin{equation}
 \theta_z^K=\frac{d+1}{|K|}\big((d+2)\lambda_z^K-1\big) \label{eq:local-dual-formula}
\end{equation}
and write \(\mathcal N(K)\) for the vertices of \(K\). The identities
for integrals of barycentric coordinates give
\begin{equation}
 \int_K\theta_z^K\lambda_y^K\,dx=\delta_{zy} \quad \forall y\in\mathcal N(K) \label{eq:local-duality}
\end{equation}
where \(\delta_{zy}\) is the Kronecker delta. The weak local coefficients
and the coarse projection are
\begin{equation}
 \ell_z(v)=\int_{K_z}\theta_z^{K_z}v\,dx,\qquad \Pi_Hv=\sum_{z\in\mathcal N_H^{\rm int}}\ell_z(v)\lambda_z \label{eq:projection-definition}
\end{equation}
for \(v\in L^2(\Omega)\). The dual functions \(\theta_z^{K_z}\) may change sign. This operator
uses local integrals rather than point values, and only interior coarse
hats enter its output. We call \(\Pi_H\) the local weak projection. It
is a projection of Scott--Zhang type~\cite{scott-zhang1990}: one simplex
is selected for each vertex, and the coefficients are integrals against
the \(L^2\)-dual basis of the nodal basis on that simplex. It differs
from the original operator in two respects. The averaging domain of a
vertex is a full coarse element and not a face, so that \(\Pi_H\) is
defined on \(L^2(\Omega)\), and boundary vertices are omitted. For
functions with zero trace, the omitted coefficients vanish also for the
Scott--Zhang operator with boundary faces.

The fine space and the correction operator are
\begin{equation}
 V_f=\ker(\Pi_H|_{V_h}),\qquad \mathcal K_H=I-\Pi_H:V_h\longrightarrow V_f \label{eq:fine-kernel-definition}
\end{equation}
where \(I\) is the identity operator. The following lemma verifies the
geometric patch estimates used below. These estimates lead to the
projection and kernel bounds in Lemma~\ref{lem:local-projection}.

\begin{lem}[One-layer patch inequalities]
\label{lem:one-layer-patches}
For the vertex-sharing element graph of
Section~\ref{subsec:spaces-layers}, the patches in
\eqref{eq:patch-poincare} and \eqref{eq:boundary-poincare} may be
chosen as \(N(T)\). Their constants are uniform over the stated mesh
family.
\end{lem}

\begin{proof}
The simplices incident to a vertex are connected through full
\((d-1)\)-faces. Indeed, intersect a sufficiently small neighborhood
of that vertex with the domain. At an interior vertex this is a ball;
at a boundary vertex the Lipschitz polyhedral domain has a connected
local link. A path in its interior, perturbed to avoid faces of
codimension at least two, passes between incident simplices through
full faces. Thus each vertex star is face-connected. The set \(N(T)\)
is the union of the vertex stars of \(T\), each containing \(T\),
and is therefore face-connected. Shape regularity and quasi-uniformity
bound its number of elements and its diameter by constants and \(CH\),
respectively.

Let \(v_K=|K|^{-1}\int_Kv\,dx\). Scaling the simplex Poincar\'e and
trace inequalities, for two elements \(K,L\) sharing a full face,
gives
\begin{equation}
 |v_K-v_L|\le C H^{1-d/2}\|\nabla v\|_{K\cup L} \label{eq:neighbor-mean-control}
\end{equation}
because the two traces agree and can be compared with their common
face mean. Connect every element of \(N(T)\) to \(T\) by a face
chain of uniformly bounded length, sum
\eqref{eq:neighbor-mean-control}, and use the simplex Poincar\'e
inequality. This bounds \(\|v-v_T\|_{N(T)}\) by
\(CH\|\nabla v\|_{N(T)}\). Subtraction of the patch mean can only
decrease that norm.

If \(T\) has a boundary vertex, its vertex star contains a physical
boundary facet. Let \(K_0\subset N(T)\) contain such a facet.
Its diameter and facet measure are comparable to \(H\) and
\(H^{d-1}\). Zero trace on that facet and the simplex trace inequality
give \(|v_{K_0}|\le CH^{1-d/2}\|\nabla v\|_{K_0}\).
The same face-chain argument, now starting at \(K_0\), proves the
zero-trace version. These arguments involve only coarse simplices and
apply to every \(H^1\) function on the patches.
\end{proof}

\begin{lem}[Local projection and approximation]
\label{lem:local-projection}
The operator \(\Pi_H\) reproduces \(V_H\). One may take the fixed radii
\(r=r_P=1\), with uniform constants \(C_{\rm loc},C_P,C_\Pi,C_J\)
such that
\begin{align}
 \Pi_Hv_H&=v_H \quad \forall v_H\in V_H \label{eq:coarse-reproduction}\\
 \|\nabla\Pi_Hv\|_T+H^{-1}\|v-\Pi_Hv\|_T&\le C_{\rm loc}\|\nabla v\|_{N^r(T)} \quad \forall v\in V_h \label{eq:local-approximation}\\
 H^{-1}\|w_f\|_T&\le C_P\|\nabla w_f\|_{N^{r_P}(T)} \quad \forall w_f\in V_f \label{eq:kernel-poincare}\\
 \|\nabla\Pi_Hv\|_\Omega&\le C_\Pi\|\nabla v\|_\Omega \quad \forall v\in V_h \label{eq:projection-global-stability}\\
 \|\nabla\mathcal K_Hv\|_\Omega&\le C_J\|\nabla v\|_\Omega \quad \forall v\in V_h \label{eq:kernel-global-stability}
\end{align}
for every coarse element \(T\). One may take \(C_J=1+C_\Pi\).
The value of \(\mathcal K_Hv\) on \(T\) depends only on \(v\) in
\(N^r(T)\), and an input supported in an element set \(S\) has output
supported in \(N^r(S)\). The constants and radii are independent of the
fine-to-coarse mesh ratio and of the choices of \(K_z\) within the
stated mesh family.
\end{lem}

\begin{proof}
For a coarse affine function on \(K_z\), the duality
\eqref{eq:local-duality} extracts its coefficient at \(z\). This proves
\eqref{eq:coarse-reproduction}. The output coefficient at \(z\) uses
only \(K_z\), while its output hat is supported on the star of \(z\).
These elements are neighbors in the chosen graph, proving fixed locality
and support propagation.

Scaling \eqref{eq:local-dual-formula} gives
\begin{equation}
 \|\theta_z^{K_z}\|_{K_z}\le C H^{-d/2} \label{eq:dual-scaling}
\end{equation}
and the coarse hat bounds therefore yield
\begin{equation}
 \|\Pi_Hv\|_T+H\|\nabla\Pi_Hv\|_T\le C\|v\|_{N(T)} \label{eq:projection-local-ltwo}
\end{equation}
by Cauchy--Schwarz and the bounded number of vertices. If every vertex
of \(T\) is interior, take \(\omega_T=N(T)\), which contains \(T\) and all its
relevant selected elements. Subtract the mean
\(c=v_{\omega_T}\). The functionals reproduce constants, and the hats
on \(T\) sum to one, so
\begin{equation}
 (\Pi_H(v-c))|_T=(\Pi_Hv-c)|_T \label{eq:local-constant-reproduction}
\end{equation}
although constants need not be reproduced on boundary elements. The
operator in \eqref{eq:projection-definition} is defined on \(L^2\),
which makes this subtraction legitimate. Applying
\eqref{eq:projection-local-ltwo} to \(v-c\) and then
\eqref{eq:patch-poincare} proves
\begin{align}
 \|\nabla\Pi_Hv\|_T&\le CH^{-1}\|v-c\|_{\omega_T}\le C\|\nabla v\|_{\omega_T} \label{eq:projection-gradient-proof}\\
 \|v-\Pi_Hv\|_T&\le\|v-c\|_T+\|\Pi_H(v-c)\|_T\le CH\|\nabla v\|_{\omega_T} \label{eq:projection-approximation-proof}
\end{align}
using Lemma~\ref{lem:one-layer-patches}. If \(T\) has a boundary
vertex, use \(c=0\) and the boundary inequality
\eqref{eq:boundary-poincare} on the same one-layer patch. The same estimates
follow without claiming reproduction of a nonzero constant there.

This proves \eqref{eq:local-approximation} with \(r=1\).
Setting \(\Pi_Hw_f=0\) gives \eqref{eq:kernel-poincare} with
\(r_P=1\). Summation with bounded
overlap proves \eqref{eq:projection-global-stability}, and the triangle
inequality gives \eqref{eq:kernel-global-stability}.
\end{proof}

Consequently \(\mathcal K_Hw_f=w_f\) for \(w_f\in V_f\). With the same radius
\(r=1\), we also have the weaker local form
\begin{equation}
 \|\nabla\mathcal K_Hv\|_T\le C_K\big(\|\nabla v\|_{N^r(T)}+H^{-1}\|v\|_{N^r(T)}\big) \quad \forall v\in V_h \label{eq:kernel-local-stability}
\end{equation}
which is convenient for the cutoff proof. Its local \(L^2\) term is
controlled using \eqref{eq:kernel-poincare} when the input is built from
a fine-kernel function.

\begin{rem}[One operator for the local method]
\label{rem:one-operator}
The same \(\mathcal K_H=I-\Pi_H\) will correct cutoffs and assemble
the partition-of-unity blend. Lemma~\ref{lem:local-projection} supplies
both locality and gradient stability. No interior fine-element bubbles
or additional refinement threshold are needed for this projection.
If \(V_h=V_H\), the kernel is zero and the correction results are
trivial. Other local weak operators can be used after verifying the
same reproduction, locality, and approximation properties.
\end{rem}

\subsection{Interpolation of Weighted Fine Functions}
\label{subsec:product-interpolation}

Let \(I_h\) denote nodal interpolation onto continuous piecewise-affine
fine functions. We apply it only to products of coarse affine and fine
affine functions, not to arbitrary \(H^1\) inputs.

\begin{lem}[Product interpolation]
\label{lem:product-interpolation}
For a fine element \(t\), a coarse piecewise-affine function \(\eta\),
and \(w\in V_h\), one has
\begin{align}
 \|I_h(\eta w)\|_t&\le C\|\eta\|_{L^\infty(t)}\|w\|_t \label{eq:product-ltwo}\\
 \|\nabla I_h(\eta w)\|_t&\le C\big(\|\eta\|_{L^\infty(t)}\|\nabla w\|_t+\|\nabla\eta\|_{L^\infty(t)}\|w\|_t\big) \label{eq:product-gradient}
\end{align}
with a constant depending only on fine-element shape regularity.
\end{lem}

\begin{proof}
Nestedness makes \(\eta\) affine on \(t\). At edge endpoints
\(x_j,x_k\), write \(\eta_j=\eta(x_j)\) and \(w_j=w(x_j)\). Then
\begin{equation}
 \eta_jw_j-\eta_kw_k=\eta_j(w_j-w_k)+w_k(\eta_j-\eta_k) \label{eq:product-edge-identity}
\end{equation}
and fine-element scaling bounds the interpolated gradient by these edge
differences divided by the element diameter. Differences of the affine
factors contribute that same diameter. Nodal norm equivalence for
affine \(w\) then proves \eqref{eq:product-gradient}. Nodal norm
equivalence also proves \eqref{eq:product-ltwo}. In particular, a cutoff
with gradient of size \(H^{-1}\) introduces \(H^{-1}\|w\|\), with no
additional factor \(H/h\).
\end{proof}

\section{Nonlinear Multiscale Reconstruction}
\label{sec:reconstruction}

With the coarse projection fixed, we first eliminate the fine kernel
globally. We then replace this global correction by independent patch
corrections and assemble them into a conforming state with the same
coarse coordinate.

\subsection{Global Fine-Scale Elimination}
\label{subsec:global-elimination}

For a coarse coordinate \(v_H\in V_H\), define the global correction
and reconstruction by
\begin{equation}
 q(v_H)=\operatorname*{argmin}_{q\in V_f}\mathcal E(v_H+q),\qquad M_h(v_H)=v_H+q(v_H) \label{eq:global-reconstruction}
\end{equation}
which are well defined by coercivity and strict convexity on the affine
space \(v_H+V_f\). Equivalently,
\begin{equation}
 \mathcal A(M_h(v_H))[z_f]=F(z_f) \quad \forall z_f\in V_f \label{eq:global-correction-equation}
\end{equation}
and \(\Pi_HM_h(v_H)=v_H\). The ideal reconstruction is exact at the
coarse coordinate of the fine solution:
\begin{equation}
 v_H^\dagger=\Pi_Hu_h,\qquad M_h(v_H^\dagger)=u_h \label{eq:exact-comparator-coordinate}
\end{equation}
because \(u_h-\Pi_Hu_h\in V_f\) solves
\eqref{eq:global-correction-equation} and that correction is unique.

\subsection{Patch Problems and Partition of Unity}
\label{subsec:patch-assembly-definition}

Let \(\mathcal I\) index all coarse vertices, including boundary vertices,
and let \(\phi_i\) be their standard continuous coarse nodal hats. For
an interior vertex, \(\phi_i\) is the hat function \(\lambda_z\) of
Section~\ref{sec:local-projection}; the index set \(\mathcal I\) also
contains the boundary vertices, whose hats do not belong to \(V_H\).
Define
\begin{equation}
 \sum_{i\in\mathcal I}\phi_i=1,\qquad \omega_i=\operatorname{int}(\operatorname{supp}\phi_i) \label{eq:full-partition}
\end{equation}
and call \(\omega_i\) the active support of patch \(i\). Choose a coarse
element patch \(D_i\supset\omega_i\). For an element
region \(D\), the supported fine space is
\begin{equation}
 V_f^0(D)=\{z_f\in V_f:\operatorname{supp}z_f\subset\overline D\} \label{eq:supported-fine-space}
\end{equation}
so each member is already a global conforming fine-kernel function. For
each \(v_H\), let \(q_i(v_H)\in V_f^0(D_i)\) solve
\begin{equation}
 \int_{D_i}a(x,\nabla(v_H+q_i(v_H)))\cdot\nabla z_f\,dx=F(z_f) \quad \forall z_f\in V_f^0(D_i) \label{eq:patch-equation}
\end{equation}
and regard it as zero outside \(D_i\). This is the constrained energy
minimizer over \(v_H+V_f^0(D_i)\), hence unique. Each patch uses the restriction of the same global load. The
partition-of-unity function \(\phi_i\) enters only during assembly,
not in the patch equation.

For the collection \(\mathcal D=(D_i)_{i\in\mathcal I}\), define
\begin{equation}
 M_h^{\mathcal D}(v_H)=v_H+(I-\Pi_H)I_h\left(\sum_{i\in\mathcal I}\phi_iq_i(v_H)\right) \label{eq:pou-reconstruction}
\end{equation}
as in~\cite{larson-persson2026}. Only the restriction of \(q_i\) to
\(\omega_i\) enters the product. Fine interpolation restores membership
in \(V_h\), and \(I-\Pi_H\) restores membership of the correction in
\(V_f\). Thus
\begin{equation}
 \Pi_HM_h^{\mathcal D}(v_H)=v_H \quad \forall v_H\in V_H \label{eq:localized-section}
\end{equation}
so the reconstruction preserves its prescribed coarse coordinate. A map
\(\widetilde M:V_H\to V_h\) with \(\Pi_H\widetilde M(v_H)=v_H\) for all
\(v_H\in V_H\) is called a section of \(\Pi_H\); we refer to
\eqref{eq:localized-section} as the section property. The global \(q\) also satisfies
\(I_hq=q\) and \((I-\Pi_H)q=q\). It follows that
\begin{equation}
 M_h^{\mathcal D}(v_H)-M_h(v_H)=(I-\Pi_H)I_h\left(\sum_{i\in\mathcal I}\phi_i(q_i(v_H)-q(v_H))\right) \label{eq:exact-assembly-error}
\end{equation}
by \eqref{eq:full-partition}. The inclusion of boundary hats in the
partition is needed for this identity, although there are no additional
boundary coarse unknowns.

Uniform-depth patches and their reconstruction are denoted by
\begin{equation}
 D_i=N^\ell(\omega_i),\qquad M_h^\ell=M_h^{\mathcal D} \label{eq:uniform-patches}
\end{equation}
for an integer \(\ell\ge0\). Variable depths will be allowed in the
weighted reconstruction estimate
(Theorem~\ref{thm:weighted-reconstruction}), but every patch remains fixed as the coarse
input varies.

\subsection{Continuity and Tangent Maps}
\label{subsec:map-regularity}

The dependence of the patch solutions on the coarse coordinate
provides the tangent directions used in the reduced equation.
Keeping the patches fixed is essential to the following argument.

\begin{lem}[Regularity of fixed-patch reconstructions]
\label{lem:reconstruction-regularity}
Every fixed patch map \(q_i:V_H\to V_f^0(D_i)\) is continuously
differentiable. Consequently \(M_h^{\mathcal D}\) is continuously
differentiable and satisfies \eqref{eq:localized-section}. For inputs
\(v_H,w_H\), its individual patch states satisfy
\begin{equation}
 \|\nabla[(v_H+q_i(v_H))-(w_H+q_i(w_H))]\|_{D_i}\le\chi\|\nabla(v_H-w_H)\|_{D_i} \label{eq:patch-state-lipschitz}
\end{equation}
and the corresponding correction difference is bounded by \(1+\chi\)
times the input gradient difference on the patch.
\end{lem}

\begin{proof}
Write \(t=(v_H+q_i(v_H))-(w_H+q_i(w_H))\). Subtracting the patch
equations and testing with \(q_i(v_H)-q_i(w_H)\) gives
\begin{align}
 m\|\nabla t\|_{D_i}^2&\le\int_{D_i}[a(x,\nabla(v_H+q_i(v_H)))-a(x,\nabla(w_H+q_i(w_H)))]\cdot\nabla(v_H-w_H)\,dx \label{eq:patch-lipschitz-test}\\
 &\le M\|\nabla t\|_{D_i}\|\nabla(v_H-w_H)\|_{D_i} \label{eq:patch-lipschitz-bound}
\end{align}
which proves \eqref{eq:patch-state-lipschitz}; the correction bound
follows by subtraction of \(v_H-w_H\). In fixed bases the assembled
patch residual is \(C^1\): continuity of \(D_\xi a\) and its uniform
bound \(M\) justify differentiation under the integral. Its derivative
with respect to the fine coefficients is positive definite by
\eqref{eq:hessian-bounds}. The implicit-function theorem therefore gives
a local \(C^1\) solution map. Uniqueness joins these local maps over
all coarse inputs. The reconstruction uses a finite sum of fixed linear
operations and therefore has the same regularity.
\end{proof}

For \(\tau_H\in V_H\), let \(\dot q_i=Dq_i(v_H)[\tau_H]\). The
derivative is characterized by
\begin{equation}
 \int_{D_i}\nabla z_f^T D_\xi a(x,\nabla(v_H+q_i(v_H)))\nabla(\tau_H+\dot q_i)\,dx=0 \quad \forall z_f\in V_f^0(D_i) \label{eq:patch-tangent}
\end{equation}
with \(\dot q_i\in V_f^0(D_i)\). Therefore
\begin{equation}
 DM_h^{\mathcal D}(v_H)[\tau_H]=\tau_H+(I-\Pi_H)I_h\left(\sum_i\phi_iDq_i(v_H)[\tau_H]\right) \label{eq:assembled-tangent}
\end{equation}
is the derivative of the same reconstruction used for the state.

\subsection{Reduced Energy and Evaluation of the Method}
\label{subsec:reduced-method}

For fixed patches the coarse objective is
\begin{equation}
 J_{\mathcal D}(v_H)=\mathcal E(M_h^{\mathcal D}(v_H)) \quad \forall v_H\in V_H \label{eq:reduced-energy-definition}
\end{equation}
and the localized method seeks a global minimizer of this finite-dimensional
function. Evaluating it at one coarse input consists of solving
\eqref{eq:patch-equation} on each patch, forming
\eqref{eq:pou-reconstruction}, and evaluating the original energy of the
assembled state. The patch problems are independent for a fixed input.
Their energy values are not summed to define the coarse objective,
because blending and projection change the reconstructed state.

The chain rule gives its derivative:
\begin{equation}
 DJ_{\mathcal D}(v_H)[\tau_H]=(\mathcal A(M_h^{\mathcal D}(v_H))-F)[DM_h^{\mathcal D}(v_H)[\tau_H]] \label{eq:reduced-energy-derivative}
\end{equation}
where each tangent correction is obtained from \eqref{eq:patch-tangent}.
Thus the tangent-space coarse equation is exactly the stationarity
condition for this objective.
We call \(DJ_{\mathcal D}(v_H)\) the reduced residual at \(v_H\). Its
dual norm \(\rho_H(v_H)\) is defined in \eqref{eq:reduced-residual-norm},
and the term reduced residual is also used for this number.
Global minimality gives the energy comparison in
Section~\ref{sec:energy-comparison}. For a computed stationary point,
Section~\ref{sec:stationary-residual} instead proves an error estimate
from tangent stability and sufficient oversampling. A nonzero reduced
residual enters that estimate as an explicit additional term.

The local support constraint in \eqref{eq:supported-fine-space} must
also be retained in an implementation. In terms of the weak degrees
of freedom, a zero-extended patch function satisfies
\begin{equation}
 \ell_z(q_i)=0 \quad \forall z\in\mathcal N_H^{\rm int} \label{eq:patch-global-moments}
\end{equation}
including any selected element \(K_z\) that intersects the patch,
even if its associated vertex is outside the patch interior. Constraints
whose selected elements have no overlap with the support are automatic.
This gives a local constrained system with the same global kernel used
in the proof. Replacing these constraints by an independently defined
patch projection would define a different reconstruction.

\section{Direct Energy Comparison}
\label{sec:energy-comparison}

The localized reconstruction of the fine solution, taken at the coarse
coordinate \(\Pi_Hu_h\), lies on the manifold, so the energy of a
manifold minimizer is at most its energy. This transfers the accuracy of
the reconstruction to the minimizer.
The next theorem isolates this variational step; the localization
estimate needed to bound the reconstruction is proved in Sections~\ref{sec:localization}
and~\ref{sec:assembly-error}. Let \(\widetilde M:V_H\to V_h\) be any continuous map satisfying
\begin{equation}
 \Pi_H\widetilde M(v_H)=v_H \quad \forall v_H\in V_H \label{eq:general-section}
\end{equation}
that is, a continuous section of \(\Pi_H\). In the application,
\(\widetilde M=M_h^{\mathcal D}\) for fixed patches. Define the trial
manifold and its analytical minimizers by
\begin{equation}
 \widetilde{\mathcal M}=\widetilde M(V_H),\qquad u_{\rm ms}\in\operatorname*{argmin}_{w\in\widetilde{\mathcal M}}\mathcal E(w) \label{eq:manifold-minimizer}
\end{equation}
where global minimality, rather than only stationarity, is part of the
definition.

\begin{thm}[Energy comparison on a section]
\label{thm:energy-comparison}
Let \(\widetilde M\) be a continuous section of \(\Pi_H\) as in
\eqref{eq:general-section}, with the trial manifold and its minimizers
defined in \eqref{eq:manifold-minimizer}, and let the constitutive bounds
\eqref{eq:potential-normalization}--\eqref{eq:hessian-bounds} hold. Then
the minimum in \eqref{eq:manifold-minimizer} is attained. Put
\(w_* =\widetilde M(\Pi_Hu_h)\). Every global minimizer is a best approximation in the potential Bregman
distance:
\begin{equation}
 \mathcal E(u_{\rm ms})-\mathcal E(u_h)=\min_{w\in\widetilde{\mathcal M}}\mathcal B(w,u_h) \label{eq:bregman-best-approximation}
\end{equation}
and satisfies
\begin{align}
 0\le\mathcal E(u_{\rm ms})-\mathcal E(u_h)&\le\mathcal E(w_*)-\mathcal E(u_h) \label{eq:energy-comparator}\\
 \|\nabla(u_{\rm ms}-u_h)\|_\Omega&\le\sqrt{M/m}\|\nabla(w_*-u_h)\|_\Omega \label{eq:gradient-quasioptimality}
\end{align}
In particular, if
\(\|\nabla(w_*-u_h)\|_\Omega^2\le G\), then
\begin{equation}
 \mathcal E(u_{\rm ms})-\mathcal E(u_h)\le\frac M2G,\qquad \|\nabla(u_{\rm ms}-u_h)\|_\Omega^2\le\frac MmG \label{eq:energy-gradient-insertion}
\end{equation}
with the same reconstruction bound \(G\).
\end{thm}

\begin{proof}
Take a minimizing sequence \(\widetilde M(v_H^j)\). Its energies are
bounded above and \eqref{eq:energy-coercivity} bounds the reconstructed
states in \(V_h\). The section property \eqref{eq:general-section} bounds
\(v_H^j=\Pi_H\widetilde M(v_H^j)\) in \(V_H\). Finite dimensionality
gives a convergent subsequence of coarse coordinates; continuity of
\(\widetilde M\) and \(\mathcal E\) gives attainment.

Identity \eqref{eq:energy-bregman} and global minimality prove
\eqref{eq:bregman-best-approximation}. The state \(w_*\) lies on
\(\widetilde{\mathcal M}\), proving
\eqref{eq:energy-comparator}. Applying
\eqref{eq:energy-gap-equivalence} on both sides yields
\begin{equation}
 \frac m2\|\nabla(u_{\rm ms}-u_h)\|_\Omega^2\le\mathcal E(u_{\rm ms})-\mathcal E(u_h)\le\frac M2\|\nabla(w_*-u_h)\|_\Omega^2 \label{eq:energy-comparison-chain}
\end{equation}
and proves the remaining statements.
\end{proof}

The theorem uses no smallness assumption, and its proof does not
differentiate the reconstruction. If
\(\widetilde M\) is \(C^1\) and
\(u_{\rm ms}=\widetilde M(v_H^*)\), differentiation at a global
minimizer gives the tangent equation
\begin{equation}
 \mathcal A(u_{\rm ms})[D\widetilde M(v_H^*)[\tau_H]]=F(D\widetilde M(v_H^*)[\tau_H]) \quad \forall\tau_H\in V_H \label{eq:tangent-stationarity}
\end{equation}
as a necessary optimality condition. Strict convexity of \(\mathcal E\)
on \(V_h\) does not imply convexity of its composition with the nonlinear
map \(\widetilde M\), and no converse to this stationarity statement is asserted.

\begin{cor}[Outer energy error]
\label{cor:outer-energy-error}
Suppose \(\widetilde u\in\widetilde{\mathcal M}\) satisfies
\(\mathcal E(\widetilde u)\le\mathcal E(u_{\rm ms})+\eta_E\) with
\(\eta_E\ge0\). Under the reconstruction bound in
Theorem~\ref{thm:energy-comparison},
\begin{equation}
 \mathcal E(\widetilde u)-\mathcal E(u_h)\le\frac M2G+\eta_E,\qquad \|\nabla(\widetilde u-u_h)\|_\Omega^2\le\frac MmG+\frac{2\eta_E}{m} \label{eq:outer-gap-bound}
\end{equation}
\end{cor}

\begin{proof}
Insert the assumed energy suboptimality into
\eqref{eq:energy-comparator}, then use
\eqref{eq:energy-gap-equivalence}.
\end{proof}

Corollary~\ref{cor:outer-energy-error} requires a bound \(\eta_E\) on
the energy gap, not merely a small reduced residual.

\begin{rem}[Signed energy]
\label{rem:energy-and-residual}
For nonzero load, zero need not belong to the localized manifold; it
cannot automatically be used as a competitor there. The bounds above
control the energy gap. In particular, a localization estimate for the energy gap does not
by itself assert that the signed multiscale energy is nonpositive.
\end{rem}

\begin{prop}[Full residual certificate]
\label{prop:full-residual-certificate}
For any \(\widetilde u\in V_h\), define
\(\rho=\|F-\mathcal A(\widetilde u)\|_{-1,h}\). Then
\begin{align}
 \|\nabla(u_h-\widetilde u)\|_\Omega&\le m^{-1}\rho \label{eq:full-residual-error}\\
 0\le\mathcal E(\widetilde u)-\mathcal E(u_h)&\le\frac{\rho^2}{2m} \label{eq:full-residual-energy}
\end{align}
where the dual norm tests all of \(V_h\).
\end{prop}

\begin{proof}
Strong monotonicity, the fine equation, and duality give
\begin{equation}
 m\|\nabla(\widetilde u-u_h)\|_\Omega^2\le(\mathcal A(\widetilde u)-F)[\widetilde u-u_h]\le\rho\|\nabla(\widetilde u-u_h)\|_\Omega \label{eq:residual-monotonicity-proof}
\end{equation}
which proves the gradient bound, including the zero-error case.
Strong convexity expanded at \(\widetilde u\) gives
\begin{equation}
 \mathcal E(\widetilde u)-\mathcal E(u_h)\le\rho\|\nabla(\widetilde u-u_h)\|_\Omega-\frac m2\|\nabla(\widetilde u-u_h)\|_\Omega^2\le\frac{\rho^2}{2m} \label{eq:residual-energy-proof}
\end{equation}
by completing the square. Minimality of \(u_h\) gives nonnegativity.
\end{proof}

We call this certificate full because its residual is tested against all
of \(V_h\), in contrast to the reduced residual
\eqref{eq:reduced-energy-derivative}, which is tested only in tangent
directions. It applies whether
or not the state is a minimizer on a multiscale manifold. A residual restricted
to the tangent space can vanish at a stationary root without controlling
\(\rho\). Computing or bounding the full dual norm is a separate
certification task, rather than a consequence of solving the coarse
equation. Its matrix evaluation is given in
Section~\ref{subsec:computable-certificates} of the supplementary
material.

For a computed coarse state, the estimates of this section apply through
global energy comparison, a bound on the
energy gap, or the full fine-space residual. Section~\ref{sec:stationary-residual}
provides a further route: a reduced-residual bound under sufficient
oversampling, including states that are not known to be global minimizers.

\section{Interior Localization}
\label{sec:localization}

We now estimate how a patch correction differs from the global
correction at the same coarse coordinate. Subtracting their equations
gives an exact linear equation for this difference, with a coefficient
that depends on both nonlinear states. A kernel-preserving cutoff
then turns that equation into interior decay. We assume that
Assumption~\ref{ass:standing} holds throughout this section and
Section~\ref{sec:assembly-error}.

\subsection{The Secant Equation}
\label{subsec:full-chord-equation}

Fix one coarse input \(v_H\), one patch \(D\), and the corresponding
global and patch states. Write
\begin{equation}
 U=v_H+q(v_H),\qquad U_D=v_H+q_D(v_H),\qquad e=U_D-U\in V_f \label{eq:localization-pair}
\end{equation}
where \(q_D\) is zero-extended. Define the secant tensor by
\begin{equation}
 B(x)=\int_0^1D_\xi a(x,\nabla U+t\nabla e)\,dt \label{eq:full-chord-tensor}
\end{equation}
which is fixed when the endpoints are fixed. The fundamental theorem of
calculus gives
\begin{equation}
 a(x,\nabla U_D)-a(x,\nabla U)=B(x)\nabla e \label{eq:full-chord-flux}
\end{equation}
and \eqref{eq:hessian-bounds} gives the same bounds for \(B\).

\begin{lem}[Same-load secant harmonicity]
\label{lem:secant-harmonicity}
For the pair \eqref{eq:localization-pair},
\begin{align}
 \int_D\nabla e^TB\nabla z_f\,dx&=0 \quad \forall z_f\in V_f^0(D) \label{eq:secant-harmonicity}\\
 \mathscr D_S(U_D,U)&=\int_S\nabla e^TB\nabla e\,dx \label{eq:secant-natural-energy}
\end{align}
for every element region \(S\). These identities also hold on nested
subregions with their supported test spaces.
\end{lem}

\begin{proof}
Every \(z_f\in V_f^0(D)\) is an admissible global test. Its support
restricts the global integral to \(D\). The global and patch equations
therefore have exactly the same right-hand side on this test space.
Subtract them and use \eqref{eq:full-chord-flux} to prove
\eqref{eq:secant-harmonicity}. Multiplication of
\eqref{eq:full-chord-flux} by \(\nabla e\) proves
\eqref{eq:secant-natural-energy}. Supported test spaces of subregions
are contained in \(V_f^0(D)\), giving the final assertion.
\end{proof}

We say that the pair \eqref{eq:localization-pair} is secant-harmonic on
an element region \(S\subset D\) if \eqref{eq:secant-harmonicity} holds
with \(S\) in place of \(D\). By the last assertion of
Lemma~\ref{lem:secant-harmonicity}, this is the case for every such
\(S\). The secant equation is exact: its coefficient averages the Hessian
between the two nonlinear states. Its right-hand side vanishes because
the loads agree on the common test space, not because the patch residual
is assumed to be supported near the boundary. Different local loads
would leave their difference on the right-hand side.

\subsection{A Fine-Kernel Cutoff}
\label{subsec:kernel-cutoff}

Let \(S\) be a set of coarse elements, and let \(\eta_H\) be a coarse
cut-off function: continuous and piecewise affine on the coarse mesh,
with \(0\le\eta_H\le1\), equal to \(1\) on \(N(S)\) and to \(0\) outside
\(N^2(S)\). For \(w_f\in V_f\), the product \(\eta_Hw_f\) need not belong
to \(V_f\). We use the corrected cutoff of \(w_f\),
\(z_f=(I-\Pi_H)I_h(\eta_Hw_f)\), instead; the locality of
\(\mathcal K_H=I-\Pi_H\) keeps its support within \(N^3(S)\). The
construction is given in \eqref{eq:geometric-cutoff} and
\eqref{eq:corrected-cutoff-definition}, and the following lemma is
proved with explicit buffers in Appendix~\ref{sec:buffered-cutoff}.

\begin{lem}[Buffered corrected cutoff]
\label{lem:corrected-cutoff}
For the local dual projection \eqref{eq:projection-definition},
one may take \(Q=3\). There is a uniform constant \(C_{\rm cut}\) such
that, for every coarse element set \(S\) and \(w_f\in V_f\), the
corrected cutoff \(z_f\) of \(w_f\) satisfies
\begin{align}
 D^-&=S,\qquad D^+=N^Q(S),\qquad A=D^+\setminus\overline{D^-} \label{eq:cutoff-regions}\\
 z_f&\in V_f^0(D^+),\qquad z_f=w_f\text{ on }D^- \label{eq:cutoff-matching}\\
 \|\nabla z_f\|_A^2&\le C_{\rm cut}\|\nabla w_f\|_A^2 \label{eq:annulus-cutoff-estimate}
\end{align}
under the element-union convention.
\end{lem}

The constant \(C_{\rm cut}\) is independent of the mesh ratio, the set
\(S\), and \(w_f\), and may be chosen as in
\eqref{eq:three-layer-cutoff-constant}.
The annulus on the right of \eqref{eq:annulus-cutoff-estimate} is
essential. A bound involving the energy in all of \(D^+\) does not
directly give the contraction in Proposition~\ref{prop:block-contraction}. The three-layer construction keeps every neighborhood
used to estimate the transition part inside this annulus. The key observation is that projection coefficients vanish on selected
elements where the cut-off function is constant. Only the transition set therefore
contributes to the projection correction.

\subsection{Contraction and Iteration}
\label{subsec:contraction}

\begin{prop}[One-block contraction]
\label{prop:block-contraction}
Suppose the pair \eqref{eq:localization-pair} is secant-harmonic on
\(D^+=N^Q(D^-)\subset D\), that is, \eqref{eq:secant-harmonicity} holds
with \(D^+\) in place of \(D\). Define
\begin{equation}
 \kappa=\sqrt{\chi C_{\rm cut}},\qquad \theta=\frac{\kappa}{1+\kappa}\in(0,1) \label{eq:contraction-factor}
\end{equation}
with \(\chi=M/m\). Then
\begin{equation}
 \mathscr D_{D^-}(U_D,U)\le\theta\mathscr D_{D^+}(U_D,U) \label{eq:block-contraction}
\end{equation}
\end{prop}

\begin{proof}
Let \(z_f\) be the corrected cutoff of \(e\) from
Lemma~\ref{lem:corrected-cutoff}, with \(w_f=e\). Let
\(E_-=\int_{D^-}\nabla e^TB\nabla e\,dx\) and
\(E_A=\int_A\nabla e^TB\nabla e\,dx\). By
\eqref{eq:secant-natural-energy}, \(E_-=\mathscr D_{D^-}(U_D,U)\) and
\(E_-+E_A=\mathscr D_{D^+}(U_D,U)\). Secant harmonicity on \(D^+\) and
the matching \(z_f=e\) on \(D^-\) from \eqref{eq:cutoff-matching} give
\begin{align}
 E_-&=-\int_A\nabla e^TB\nabla z_f\,dx \label{eq:inner-annular-identity}\\
 &\le E_A^{1/2}\left(\int_A\nabla z_f^TB\nabla z_f\,dx\right)^{1/2} \label{eq:annular-cauchy}\\
 &\le E_A^{1/2}\big(MC_{\rm cut}\|\nabla e\|_A^2\big)^{1/2}\le\kappa E_A \label{eq:annular-energy-comparison}
\end{align}
where the last step uses the lower bound \(m\) for \(B\).
Since the outer energy is \(E_-+E_A\), rearrangement gives
\begin{equation}
 E_-\le\frac{\kappa}{1+\kappa}(E_-+E_A) \label{eq:energy-contraction-rearrangement}
\end{equation}
which is \eqref{eq:block-contraction}. No positive-energy division was
used, so zero energies are included.
\end{proof}

The factor \(\theta\) depends only on \(\chi\) and \(C_{\rm cut}\), so
the contraction is uniform in the state amplitudes.

\begin{cor}[Interior decay]
\label{cor:interior-decay}
For a fixed pair and an inner element set \(S\), if
\(N^{nQ}(S)\subset D\), then
\begin{equation}
 \mathscr D_S(U_D,U)\le\theta^n\mathscr D_D(U_D,U) \quad \forall n\in\mathbb N_0,\quad N^{nQ}(S)\subset D \label{eq:interior-decay}
\end{equation}
where \(\mathbb N_0\) denotes the nonnegative integers.
\end{cor}

\begin{proof}
Apply Proposition~\ref{prop:block-contraction} to the
successive pairs \(N^{jQ}(S)\subset N^{(j+1)Q}(S)\), for
\(j=0,\ldots,n-1\), and then use
nonnegativity to bound the final regional energy by the energy on \(D\).
For \(n=0\), the estimate follows directly from \(S\subset D\).
The endpoints and \(B\) remain fixed at every step. No separately
computed patch states are introduced at the intermediate radii.
\end{proof}

If \(D=\Omega\), uniqueness gives \(q_D=q\) and the error vanishes
exactly. This exact terminal case will be assigned patch weight zero in the
assembly estimate, independently of the positive factors in
\eqref{eq:contraction-factor}.

\section{Assembly and Stationary Error Bounds}
\label{sec:assembly-error}

Interior decay (Corollary~\ref{cor:interior-decay}) controls each patch
error \(q_i-q\) on the active support
\(\omega_i\) of its partition-of-unity function \(\phi_i\), which is the
only part of the patch solution that enters the reconstruction
\eqref{eq:pou-reconstruction}.
To obtain a global reconstruction estimate, we first bound the effect
of the partition of unity and the projection, then sum the patch contributions using
their actual overlap. At the fine solution's coarse coordinate, this
estimate supplies the reconstruction bound \(G\) required by
Theorem~\ref{thm:energy-comparison}.

\subsection{Stability on Enlarged Active Supports}
\label{subsec:active-support-stability}

Let \(r_P=1\), as verified in Lemma~\ref{lem:local-projection}. The
estimates below use the enlarged active supports
\begin{equation}
 N(\omega_i)=N^{r_P}(\omega_i) \label{eq:enlarged-active-support}
\end{equation}
of the patches. Summing the kernel Poincar\'e
estimate \eqref{eq:kernel-poincare} over the elements of \(\omega_i\) gives a uniform \(C_{FP}\) such that
\begin{equation}
 H^{-2}\|w_f\|_{\omega_i}^2\le C_{FP}\|\nabla w_f\|_{N(\omega_i)}^2 \quad \forall w_f\in V_f \label{eq:active-kernel-poincare}
\end{equation}
and Lemma~\ref{lem:product-interpolation} gives a uniform \(C_I\) with
\begin{equation}
 \|\nabla I_h(\phi_iw_f)\|_{\omega_i}^2\le C_I\big(\|\nabla w_f\|_{\omega_i}^2+H^{-2}\|w_f\|_{\omega_i}^2\big) \label{eq:active-product-estimate}
\end{equation}
because the coarse hats have gradients bounded by \(CH^{-1}\).
Let \(N_{\rm pu}\) be the largest number of active supports meeting
one fine element, and set
\begin{equation}
 C_{\rm asm}=C_J^2N_{\rm pu}C_I(1+C_{FP}) \label{eq:assembly-constant}
\end{equation}
using \(C_J\) from \eqref{eq:kernel-global-stability}.

\begin{lem}[Partition-of-unity assembly]
\label{lem:pou-assembly}
For any family \((w_i)_{i\in\mathcal I}\) of global members of \(V_f\),
\begin{equation}
 \left\|\nabla(I-\Pi_H)I_h\left(\sum_i\phi_iw_i\right)\right\|_\Omega^2\le C_{\rm asm}\sum_i\|\nabla w_i\|_{N(\omega_i)}^2 \label{eq:pou-assembly-estimate}
\end{equation}
with the fixed enlarged supports \eqref{eq:enlarged-active-support}.
\end{lem}

\begin{proof}
Linearity of \(I_h\), the stability \eqref{eq:kernel-global-stability} of \(I-\Pi_H\), and the pointwise
overlap bound yield
\begin{equation}
 \left\|\nabla(I-\Pi_H)\sum_iI_h(\phi_iw_i)\right\|_\Omega^2\le C_J^2N_{\rm pu}\sum_i\|\nabla I_h(\phi_iw_i)\|_{\omega_i}^2 \label{eq:assembly-overlap-step}
\end{equation}
because each function \(I_h(\phi_iw_i)\) is supported in the active
support \(\omega_i\). Apply \eqref{eq:active-product-estimate} and
\eqref{eq:active-kernel-poincare} to conclude.
\end{proof}

\subsection{Outer Energy and Weighted Patch Overlap}
\label{subsec:outer-energy}

In this subsection assume \(N(\omega_i)\subset D_i\). Fix \(v_H\), and
abbreviate \(q=q(v_H)\), \(q_i=q_i(v_H)\), \(U=v_H+q\), and
\(U_i=v_H+q_i\). For a proper patch choose any integer \(n_i\ge0\)
with \(N^{n_iQ+r_P}(\omega_i)\subset D_i\), and set \(t_i=\theta^{n_i}\).
For a full-domain patch set \(t_i=0\), since \(q_i=q\). Define
\begin{equation}
 \mathcal N_t(\mathcal D)=\left\|\sum_{i\in\mathcal I}t_i\mathbf1_{D_i}\right\|_{L^\infty(\Omega)} \label{eq:weighted-overlap}
\end{equation}
where \(\mathbf1_{D_i}\) is the indicator of the patch. Thus \(\mathcal N_t\) measures overlap after each patch contribution
has been reduced by its decay factor. We also define the residual of
the coarse state, tested against the fine kernel, and its dual norm:
\begin{equation}
 r_{v_H}(w_f)=F(w_f)-\mathcal A(v_H)[w_f],\qquad R(v_H)=\sup_{0\ne w_f\in V_f}\frac{|r_{v_H}(w_f)|}{\|\nabla w_f\|_\Omega} \label{eq:fine-residual-norm}
\end{equation}
with zero norm on a zero-dimensional kernel.

\begin{lem}[Orthogonal outer-amplitude estimate]
\label{lem:outer-amplitude}
For each pair \((U_i,U)\), let \(B_i\) be its secant tensor.
The exact corrections satisfy
\begin{align}
 \mathscr D_{D_i}(U_i,U)+\int_{D_i}\nabla q_i^TB_i\nabla q_i\,dx&=\int_{D_i}\nabla q^TB_i\nabla q\,dx \label{eq:outer-pythagoras}\\
 \sum_i t_i\mathscr D_{D_i}(U_i,U)&\le M\mathcal N_t(\mathcal D)\|\nabla q\|_\Omega^2\le\frac{M}{m^2}\mathcal N_t(\mathcal D)R(v_H)^2 \label{eq:weighted-outer-energy}
\end{align}
\end{lem}

\begin{proof}
Since \(q_i\in V_f^0(D_i)\), secant harmonicity
\eqref{eq:secant-harmonicity} on \(D_i\) gives
\begin{equation}
 \int_{D_i}\nabla(q_i-q)^TB_i\nabla q_i\,dx=0 \label{eq:outer-orthogonality}
\end{equation}
by testing with the patch correction itself. Expanding the squared
\(B_i\)-seminorm of \(q=q_i-(q_i-q)\) proves
\eqref{eq:outer-pythagoras}. Positivity of its second term and the
upper bound in \eqref{eq:hessian-bounds} yield
\begin{equation}
 \mathscr D_{D_i}(U_i,U)\le M\|\nabla q\|_{D_i}^2 \label{eq:outer-local-amplitude}
\end{equation}
which sums to the first inequality in \eqref{eq:weighted-outer-energy}
by the definition \eqref{eq:weighted-overlap} of the weighted overlap. Finally, test the global
correction equation after subtracting \(\mathcal A(v_H)\) with \(q\).
Strong monotonicity gives
\begin{equation}
 \|\nabla q\|_\Omega\le m^{-1}R(v_H) \label{eq:correction-amplitudes}
\end{equation}
including the zero-correction case. This proves the second inequality.
No separate estimates of local residual norms are needed.
\end{proof}

In each identity of Lemma~\ref{lem:outer-amplitude}, the tensor \(B_i\)
of the patch is held fixed.

\subsection{The Reconstruction Estimate}
\label{subsec:reconstruction-estimate}

The assembly lemma (Lemma~\ref{lem:pou-assembly}) controls the error on
the active supports \(\omega_i\). Interior decay
(Corollary~\ref{cor:interior-decay}) relates those errors to the patch
energies, and the orthogonality estimate
(Lemma~\ref{lem:outer-amplitude}) bounds their weighted sum. Combining these three steps gives the reconstruction
estimate.

\begin{thm}[Weighted reconstruction error]
\label{thm:weighted-reconstruction}
Under Assumption~\ref{ass:standing},
for the fixed patches and patch weights \(t_i\) of
Section~\ref{subsec:outer-energy} and every \(v_H\in V_H\),
\begin{equation}
 \|\nabla(M_h^{\mathcal D}(v_H)-M_h(v_H))\|_\Omega^2\le C_{\rm asm}\chi\mathcal N_t(\mathcal D)\|\nabla q(v_H)\|_\Omega^2 \label{eq:weighted-correction-reconstruction}
\end{equation}
and consequently
\begin{align}
 \|\nabla(M_h^{\mathcal D}(v_H)-M_h(v_H))\|_\Omega^2&\le\frac{C_{\rm asm}M}{m^3}\mathcal N_t(\mathcal D)R(v_H)^2 \label{eq:weighted-gradient-reconstruction}\\
 \mathscr D_\Omega(M_h^{\mathcal D}(v_H),M_h(v_H))&\le\frac{C_{\rm asm}M^2}{m^3}\mathcal N_t(\mathcal D)R(v_H)^2 \label{eq:weighted-natural-reconstruction}
\end{align}
\end{thm}

\begin{proof}
For each proper patch, Corollary~\ref{cor:interior-decay} gives
\begin{equation}
 \mathscr D_{N(\omega_i)}(U_i,U)\le t_i\mathscr D_{D_i}(U_i,U) \label{eq:patch-to-active-decay}
\end{equation}
and the same inequality is exact for a full-domain patch with patch
weight zero. Use \eqref{eq:exact-assembly-error},
Lemma~\ref{lem:pou-assembly}, and the lower bound in
\eqref{eq:natural-equivalence} to obtain
\begin{align}
 \|\nabla(M_h^{\mathcal D}(v_H)-U)\|_\Omega^2&\le C_{\rm asm}\sum_i\|\nabla(q_i-q)\|_{N(\omega_i)}^2 \label{eq:assembly-gradient-chain}\\
 &\le\frac{C_{\rm asm}}m\sum_i t_i\mathscr D_{D_i}(U_i,U) \label{eq:assembly-natural-chain}
\end{align}
Apply the two bounds of Lemma~\ref{lem:outer-amplitude} to prove
\eqref{eq:weighted-correction-reconstruction} and
\eqref{eq:weighted-gradient-reconstruction}, respectively. The global upper bound in
\eqref{eq:natural-equivalence} gives
\eqref{eq:weighted-natural-reconstruction}. Each patch has its own
secant tensor; the proof does not assign one patch tensor to the
assembled error.
\end{proof}

For uniform patches \(D_i=N^\ell(\omega_i)\), require \(\ell\ge r_P\)
so that \(N(\omega_i)\subset D_i\), and define
\begin{equation}
 n_\ell=\left\lfloor\frac{\ell-r_P}{Q}\right\rfloor,\qquad C_{\rm rec}=\frac{C_{\rm asm}C_{\rm ol}M}{m^3} \label{eq:uniform-rate-constants}
\end{equation}
where \(C_{\rm ol}\) is a quasi-uniform packing constant.
For the present projection, \(r_P=1\) by Lemma~\ref{lem:local-projection}
and \(Q=3\) by Lemma~\ref{lem:corrected-cutoff}, which
give \(n_\ell=\lfloor(\ell-1)/3\rfloor\). Each patch
lies within distance \(CH(1+\ell)\) of its seed, and counting coarse
vertices in such a ball gives
\begin{equation}
 \mathcal N_t(\mathcal D)\le C_{\rm ol}(1+\ell)^d\theta^{n_\ell} \label{eq:uniform-weighted-overlap}
\end{equation}
by choosing \(n_i=n_\ell\) on proper patches. Hence
\begin{equation}
 \|\nabla(M_h^\ell(v_H)-M_h(v_H))\|_\Omega^2\le C_{\rm rec}(1+\ell)^d\theta^{n_\ell}R(v_H)^2 \label{eq:uniform-gradient-reconstruction}
\end{equation}
with every overlap factor displayed.

\begin{rem}[Shallow patches]
\label{rem:shallow-patches}
The assumption \(\ell\ge r_P\) is needed when applying interior
decay (Corollary~\ref{cor:interior-decay}) down to \(N(\omega_i)\). For
\(\ell<r_P\), the assembly lemma (Lemma~\ref{lem:pou-assembly}) remains
valid, but this inner
region need not lie inside \(D_i\). A separate nondecaying bound follows
by estimating \(q_i-q\) directly on \(N(\omega_i)\): \(q_i\) is supported
in \(D_i\subset N(\omega_i)\), the \(N(\omega_i)\) have bounded overlap, and the
identity \eqref{eq:outer-pythagoras} bounds
\(\|\nabla q_i\|_{D_i}^2\le\chi\|\nabla q\|_{D_i}^2\).
This gives a bound \(CR(v_H)^2/m^2\) for the squared reconstruction
gradient error, with a possibly different constant. We do not extend
\eqref{eq:uniform-gradient-reconstruction} to shallow patches merely by
truncating its exponent.
\end{rem}

\subsection{Complete Stationary Estimate}
\label{subsec:stationary-estimate}

At the exact coarse coordinate the global correction is known
algebraically in terms of the fine solution:
\begin{equation}
 q(\Pi_Hu_h)=(I-\Pi_H)u_h,\qquad \|\nabla q(\Pi_Hu_h)\|_\Omega\le C_Jm^{-1}L_F \label{eq:comparator-data-bound}
\end{equation}
by \eqref{eq:kernel-global-stability} and \eqref{eq:fine-data-bounds}.
At this coordinate, the amplitude is therefore bounded directly by
the load norm; the residual \(R(v_H)\) is not needed.
For fixed, possibly nonuniform patches define
\begin{equation}
 G_{\mathcal D}=\frac{C_{\rm asm}C_J^2\chi}{m^2}\mathcal N_t(\mathcal D)L_F^2 \label{eq:weighted-data-bound}
\end{equation}
with the patch weights \(t_i\) of Section~\ref{subsec:outer-energy}.

\begin{thm}[Stationary energy and gradient error]
\label{thm:stationary-error}
Let Assumption~\ref{ass:standing} hold. Choose fixed
coarse-element patches with \(N(\omega_i)\subset D_i\).
For each proper patch choose \(n_i\in\mathbb N_0\) such that
\(N^{n_iQ+r_P}(\omega_i)\subset D_i\), and set \(t_i=\theta^{n_i}\);
set \(t_i=0\) for full-domain patches. A global energy minimizer
\(u_{\rm ms}^{\mathcal D}\) on \(M_h^{\mathcal D}(V_H)\) exists,
and every such minimizer satisfies
\begin{align}
 0\le\mathcal E(u_{\rm ms}^{\mathcal D})-\mathcal E(u_h)&\le\frac M2G_{\mathcal D} \label{eq:weighted-stationary-energy}\\
 \|\nabla(u_{\rm ms}^{\mathcal D}-u_h)\|_\Omega^2&\le\chi G_{\mathcal D} \label{eq:weighted-stationary-gradient}
\end{align}
where \(G_{\mathcal D}\) is given by \eqref{eq:weighted-data-bound}.
For uniform patches \(D_i=N^\ell(\omega_i)\), \(\ell\ge r_P\), one
may replace \(G_{\mathcal D}\) by
\begin{equation}
 G_\ell=\frac{C_{\rm asm}C_{\rm ol}C_J^2\chi}{m^2}(1+\ell)^d\theta^{n_\ell}L_F^2 \label{eq:final-data-bound}
\end{equation}
with \(n_\ell=\lfloor(\ell-r_P)/Q\rfloor\). In particular,
\begin{align}
 0\le\mathcal E(u_{\rm ms}^\ell)-\mathcal E(u_h)&\le\frac M2G_\ell \label{eq:final-energy-error}\\
 \|\nabla(u_{\rm ms}^\ell-u_h)\|_\Omega^2&\le\chi G_\ell \label{eq:final-gradient-error}\\
 \|\nabla u_{\rm ms}^\ell\|_\Omega&\le m^{-1}L_F+\sqrt{\chi G_\ell} \label{eq:multiscale-state-bound}
\end{align}
with constants independent of state amplitude and the fine-to-coarse mesh
ratio. The constants \(C_{\rm asm},C_J,C_{\rm ol},Q,C_{\rm cut}\)
depend only on the stated geometry and local operators;
\(\theta\) depends additionally on \(\chi=M/m\).
\end{thm}

\begin{proof}
The exact comparator identity \eqref{eq:exact-comparator-coordinate},
\eqref{eq:weighted-correction-reconstruction}, and
\eqref{eq:comparator-data-bound} bound the squared gradient error of
\(M_h^{\mathcal D}(\Pi_Hu_h)\) by \(G_{\mathcal D}\).
Continuity (Lemma~\ref{lem:reconstruction-regularity}) and the section
property \eqref{eq:localized-section} permit application of
Theorem~\ref{thm:energy-comparison}, proving existence and
\eqref{eq:weighted-stationary-energy}--\eqref{eq:weighted-stationary-gradient}.
For uniform patches, insert \eqref{eq:uniform-weighted-overlap}.
The state bound follows by the triangle inequality and
\eqref{eq:fine-data-bounds}.
\end{proof}

For a fixed domain scale define \(b=|\log\theta|/Q>0\). Absorbing
the fixed offset and the floor into a constant gives
\begin{align}
 \mathcal E(u_{\rm ms}^\ell)-\mathcal E(u_h)&\le C(1+\ell)^d e^{-b\ell}L_F^2 \label{eq:exponential-energy-form}\\
 \|\nabla(u_{\rm ms}^\ell-u_h)\|_\Omega&\le C(1+\ell)^{d/2}e^{-b\ell/2}L_F \label{eq:exponential-gradient-form}
\end{align}
Thus the exponential decay exponent for the gradient error is half
that for the energy gap.
After nondimensionalizing \(H\) by a fixed domain length, choosing
\(\ell\ge c\log(1/H)\) with \(c>2s/b\) gives an
\(O(H^s)L_F\) localization error in the gradient norm; the strict
inequality absorbs the polynomial logarithmic factor. The ideal map
represents \(u_h\) exactly, so no independent coarse approximation
term is present relative to \(u_h\). Comparing to a continuous weak
solution adds a separate fine-discretization error and requires its own
approximation hypotheses.

\begin{cor}[Separation of discretization and localization]
\label{cor:continuous-error}
Assume in addition that \(F\) is the restriction of a bounded linear
functional on \(H_0^1(\Omega)\) and let \(u\) be the continuous
energy minimizer. Under the hypotheses of Theorem~\ref{thm:stationary-error},
\begin{equation}
 \|\nabla(u-u_{\rm ms}^\ell)\|_\Omega^2\le\chi\inf_{v_h\in V_h}\|\nabla(u-v_h)\|_\Omega^2+\chi G_\ell \label{eq:continuous-total-error}
\end{equation}
with \(G_\ell\) defined in \eqref{eq:final-data-bound} using the
restricted load norm. No rate for the first term is assumed.
\end{cor}

\begin{proof}
The potential bounds and strong convexity give existence and uniqueness
of the continuous minimizer by coercivity and weak lower semicontinuity
on \(H_0^1(\Omega)\). Its variational equation gives the continuous
counterpart of \eqref{eq:energy-gap-equivalence}. The energy gaps split
exactly at the fine minimizer:
\begin{align}
 \mathcal E(u_{\rm ms}^\ell)-\mathcal E(u)&=\mathcal E(u_{\rm ms}^\ell)-\mathcal E(u_h)+\mathcal E(u_h)-\mathcal E(u) \label{eq:continuous-energy-splitting}\\
 &\le\frac M2G_\ell+\frac M2\inf_{v_h\in V_h}\|\nabla(u-v_h)\|_\Omega^2 \label{eq:continuous-energy-upper-bound}
\end{align}
because \(u_h\) minimizes the same energy on \(V_h\). The lower
bound \(\mathcal E(u_{\rm ms}^\ell)-\mathcal E(u)\ge
(m/2)\|\nabla(u-u_{\rm ms}^\ell)\|_\Omega^2\) proves the result.
Thus discretization and localization enter additively in energy and
in quadrature in the gradient estimate; no gradient orthogonality
between the two errors is assumed.
\end{proof}

\subsection{Meaning of the Constants}
\label{subsec:constant-dependence}

The weighted estimate is the primary result. The factor \((1+\ell)^d\)
enters only when replacing the actual weighted overlap by the uniform
packing bound \eqref{eq:uniform-weighted-overlap}. Full-domain patches contribute zero to the weighted
estimate, so it records exact recovery when every patch fills the domain.
The uniform polynomial bound can remain positive even in that exact
recovery case. The
proof does not establish that its polynomial factor is unavoidable.

The explicit contraction exponent is
\begin{equation}
 b=\frac1Q\log\left(1+\frac1{\sqrt{\chi C_{\rm cut}}}\right) \label{eq:contrast-explicit-exponent}
\end{equation}
by \eqref{eq:contraction-factor}. For fixed geometry it deteriorates
like \(\chi^{-1/2}\) as \(\chi\) increases. This is the dependence
proved by the unweighted cutoff argument, not an optimality claim.
The uniform gradient estimate \eqref{eq:final-gradient-error} can be written with its amplitude explicit:
\begin{equation}
 \|\nabla(u_{\rm ms}^\ell-u_h)\|_\Omega\le\sqrt{C_{\rm asm}C_{\rm ol}}C_J\frac\chi m(1+\ell)^{d/2}\theta^{n_\ell/2}L_F \label{eq:explicit-gradient-amplitude}
\end{equation}
so contrast affects both the prefactor and the decay rate.
Appendix~\ref{sec:buffered-cutoff} proves \(Q=3\) for the local
dual projection, and Lemma~\ref{lem:local-projection} gives \(r_P=1\).
Thus the explicit rate is
\begin{equation}
 b=\frac13\log\left(1+\frac1{\sqrt{\chi C_{\rm cut}}}\right),
 \qquad n_\ell=\left\lfloor\frac{\ell-1}{3}\right\rfloor \label{eq:three-layer-exponent}
\end{equation}
The block width is independent of the fine refinement and the selected
coarse simplices. Its effect on the exponent must be considered together
with the cutoff constant in \eqref{eq:three-layer-cutoff-constant}.
Three is a sufficient block width, not a claimed minimal patch depth.

\section{Reduced Residual and Stationary Solutions}
\label{sec:stationary-residual}

The energy-minimizer estimate (Theorem~\ref{thm:stationary-error}) does
not require a lower oversampling
threshold beyond the geometric support condition. We now obtain a
separate estimate for every coarse coordinate from its reduced residual,
provided localization is sufficiently small. The proof uses stability
of the localized tangent and the fine equation at the ideal reconstruction
with the same coordinate. It does not require convergence of tangent maps.
Throughout this section the patches are fixed and satisfy \(N(\omega_i)\subset D_i\),
and Assumption~\ref{ass:standing} holds. Enlarge the valid constants \(C_\Pi,C_J\) to at least
one if necessary.

\subsection{Stability of the Localized Tangent}
\label{subsec:localized-tangent-stability}

For any \(v_H\in V_H\), define
\begin{equation}
 U=M_h(v_H),\qquad \widetilde u=M_h^{\mathcal D}(v_H),\qquad T_v=DM_h^{\mathcal D}(v_H) \label{eq:residual-transfer-states}
\end{equation}
and denote the unweighted patch overlap by
\begin{equation}
 \mathcal N_0(\mathcal D)=\left\|\sum_i\mathbf1_{D_i}\right\|_{L^\infty(\Omega)},\qquad C_T=1+\sqrt{C_{\rm asm}\chi\mathcal N_0(\mathcal D)} \label{eq:tangent-stability-constant}
\end{equation}
where \(C_{\rm asm}\) is the assembly constant in
\eqref{eq:assembly-constant}.

\begin{lem}[Uniform tangent stability]
\label{lem:localized-tangent-stability}
For every coarse coordinate and direction,
\begin{equation}
 \Pi_HT_v\tau_H=\tau_H,\qquad \|\nabla T_v\tau_H\|_\Omega\le C_T\|\nabla\tau_H\|_\Omega \quad \forall\tau_H\in V_H \label{eq:localized-tangent-stability}
\end{equation}
with \(C_T\) independent of the state and load amplitudes and the
fine-to-coarse mesh ratio.
\end{lem}

\begin{proof}
The first identity follows by differentiating \eqref{eq:localized-section}.
For a direction \(\tau_H\), put \(s_i=Dq_i(v_H)[\tau_H]\) and
\(A_i(x)=D_\xi a(x,\nabla(v_H+q_i(v_H)))\).
The tangent equation \eqref{eq:patch-tangent} tested with
\(s_i\in V_f^0(D_i)\) gives
\begin{equation}
 \int_{D_i}\nabla(\tau_H+s_i)^TA_i\nabla s_i\,dx=0 \label{eq:patch-tangent-orthogonality}
\end{equation}
and expansion of \(\tau_H=(\tau_H+s_i)-s_i\) yields
\begin{align}
 \int_{D_i}\nabla\tau_H^TA_i\nabla\tau_H\,dx&=\int_{D_i}\nabla(\tau_H+s_i)^TA_i\nabla(\tau_H+s_i)\,dx+\int_{D_i}\nabla s_i^TA_i\nabla s_i\,dx \label{eq:patch-tangent-pythagoras}\\
 \|\nabla s_i\|_{D_i}^2&\le\chi\|\nabla\tau_H\|_{D_i}^2 \label{eq:patch-tangent-correction-bound}
\end{align}
by symmetry and the Hessian bounds. Applying Lemma~\ref{lem:pou-assembly}
to the global fine-kernel functions \(s_i\) gives
\begin{equation}
 \|\nabla(T_v\tau_H-\tau_H)\|_\Omega^2\le C_{\rm asm}\sum_i\|\nabla s_i\|_{N(\omega_i)}^2\le C_{\rm asm}\chi\mathcal N_0(\mathcal D)\|\nabla\tau_H\|_\Omega^2 \label{eq:tangent-assembly-stability}
\end{equation}
since \(N(\omega_i)\subset D_i\). The triangle inequality proves
\eqref{eq:localized-tangent-stability}. Each \(A_i\) is the fixed
endpoint Hessian for its patch; no comparison between Hessians at different
states is used.
\end{proof}

\subsection{Transfer at the Computed Coordinate}
\label{subsec:residual-transfer}

We first relate the error at a coarse coordinate to two quantities:
its reduced residual and the difference between its ideal and localized
reconstructions. Define the reduced residual in the coarse dual gradient norm by
\begin{equation}
 \rho_H(v_H)=\sup_{0\ne\tau_H\in V_H}\frac{|(\mathcal A(\widetilde u)-F)[T_v\tau_H]|}{\|\nabla\tau_H\|_\Omega} \label{eq:reduced-residual-norm}
\end{equation}
with zero supremum for a zero-dimensional coarse space. This is the norm
of \(DJ_{\mathcal D}(v_H)\) induced by the gradient norm on \(V_H\),
not an unscaled Euclidean norm of its coefficient vector. Set
\begin{equation}
 X=\|\nabla(U-u_h)\|_\Omega,\qquad \delta_v=\|\nabla(\widetilde u-U)\|_\Omega,\qquad \beta_{\mathcal D}=\chi C_\Pi C_T \label{eq:transfer-quantities}
\end{equation}
so that \(\beta_{\mathcal D}\ge1\).

\begin{prop}[Reduced-residual transfer]
\label{prop:reduced-residual-transfer}
At every coarse coordinate,
\begin{align}
 X&\le\frac{C_\Pi}{m}\rho_H(v_H)+\beta_{\mathcal D}\delta_v \label{eq:ideal-state-residual-transfer}\\
 \|\nabla(\widetilde u-u_h)\|_\Omega&\le\frac{C_\Pi}{m}\rho_H(v_H)+(1+\beta_{\mathcal D})\delta_v \label{eq:localized-state-residual-transfer}
\end{align}
\end{prop}

\begin{proof}
Let \(e=U-u_h\), \(\tau_H=\Pi_He\), and \(r_U=\mathcal A(U)-F\).
The ideal fine equation gives \(r_U[z_f]=0\) for all \(z_f\in V_f\).
The identity \(\Pi_HT_v\tau_H=\tau_H\) in
\eqref{eq:localized-tangent-stability} implies
\begin{equation}
 e-T_v\tau_H\in V_f,\qquad r_U[e]=r_U[T_v\tau_H] \label{eq:ideal-residual-cancellation}
\end{equation}
which is the cancellation needed for transfer. Strong monotonicity and
the fine equation then give
\begin{align}
 mX^2&\le r_U[e]=(\mathcal A(\widetilde u)-F)[T_v\tau_H]+(\mathcal A(U)-\mathcal A(\widetilde u))[T_v\tau_H] \label{eq:residual-transfer-identity}\\
 &\le\big(\rho_H(v_H)+MC_T\delta_v\big)\|\nabla\tau_H\|_\Omega \label{eq:residual-transfer-duality}\\
 &\le C_\Pi\big(\rho_H(v_H)+MC_T\delta_v\big)X \label{eq:residual-transfer-monotonicity}
\end{align}
using the flux Lipschitz bound \eqref{eq:flux-lipschitz}, the tangent bound
\eqref{eq:localized-tangent-stability} and the projection stability
\eqref{eq:projection-global-stability}. For \(X>0\),
divide by \(mX\); the case \(X=0\) satisfies the conclusion directly.
The triangle inequality proves the second estimate.
\end{proof}

\subsection{Sufficient Oversampling and a Data Bound}
\label{subsec:residual-oversampling}

The transfer estimate (Proposition~\ref{prop:reduced-residual-transfer})
assumes neither minimality nor a smallness condition, but it
still contains the unknown localization error
at the computed coordinate. We now bound it in terms of the ideal-state
error and the load, then absorb the error term when the patches are
sufficiently large. Retain the weighted overlap \(\mathcal N_t(\mathcal D)\) of
\eqref{eq:weighted-overlap}, including zero weights for full-domain
patches, and define
\begin{equation}
 \varepsilon_{\mathcal D}=C_J\sqrt{C_{\rm asm}\chi\mathcal N_t(\mathcal D)} \label{eq:residual-localization-factor}
\end{equation}
The constants \(C_T\) and \(\beta_{\mathcal D}\) are those of
\eqref{eq:tangent-stability-constant} and \eqref{eq:transfer-quantities}.

\begin{thm}[Computed-state error under sufficient oversampling]
\label{thm:computed-state-error}
Suppose \(\beta_{\mathcal D}\varepsilon_{\mathcal D}<1\). Then for
every \(v_H\in V_H\) and \(\widetilde u=M_h^{\mathcal D}(v_H)\),
\begin{equation}
 \|\nabla(\widetilde u-u_h)\|_\Omega\le\frac{(1+\varepsilon_{\mathcal D})C_\Pi}{m(1-\beta_{\mathcal D}\varepsilon_{\mathcal D})}\rho_H(v_H)+\frac{\varepsilon_{\mathcal D}(1+\beta_{\mathcal D})}{1-\beta_{\mathcal D}\varepsilon_{\mathcal D}}\frac{L_F}{m} \label{eq:computed-state-closed-bound}
\end{equation}
under Assumption~\ref{ass:standing}. If
\(\beta_{\mathcal D}\varepsilon_{\mathcal D}\le1/2\), define
\begin{equation}
 S_{\mathcal D}(v_H)=3\frac{C_\Pi}{m}\rho_H(v_H)+4\beta_{\mathcal D}\varepsilon_{\mathcal D}\frac{L_F}{m} \label{eq:computed-state-simplified-bound}
\end{equation}
and obtain
\begin{equation}
 \|\nabla(\widetilde u-u_h)\|_\Omega\le S_{\mathcal D}(v_H),\qquad 0\le\mathcal E(\widetilde u)-\mathcal E(u_h)\le\frac M2 S_{\mathcal D}(v_H)^2 \label{eq:computed-state-energy-bound}
\end{equation}
\end{thm}

\begin{proof}
The section property \(\Pi_HM_h(v_H)=v_H\) of the ideal reconstruction gives
\(q(v_H)=(I-\Pi_H)U\). By
\eqref{eq:weighted-correction-reconstruction}, the stability \eqref{eq:kernel-global-stability} of \(I-\Pi_H\), and
\eqref{eq:fine-data-bounds},
\begin{equation}
 \delta_v\le\sqrt{C_{\rm asm}\chi\mathcal N_t(\mathcal D)}\|\nabla q(v_H)\|_\Omega\le\varepsilon_{\mathcal D}\left(X+\frac{L_F}{m}\right) \label{eq:computed-coordinate-localization}
\end{equation}
which controls localization at the computed coordinate. Inserting this
bound into \eqref{eq:ideal-state-residual-transfer} and absorbing gives
\begin{equation}
 X\le\frac{(C_\Pi/m)\rho_H(v_H)+\beta_{\mathcal D}\varepsilon_{\mathcal D}L_F/m}{1-\beta_{\mathcal D}\varepsilon_{\mathcal D}} \label{eq:computed-state-absorption}
\end{equation}
and \(\|\nabla(\widetilde u-u_h)\|_\Omega\le X+\delta_v\)
proves \eqref{eq:computed-state-closed-bound}. If the product is at
most \(1/2\), then \(\varepsilon_{\mathcal D}\le1/2\) because
\(\beta_{\mathcal D}\ge1\). The two coefficients are bounded by
\(3C_\Pi/m\) and \(4\beta_{\mathcal D}\varepsilon_{\mathcal D}\),
respectively. Finally apply \eqref{eq:energy-gap-equivalence}.
\end{proof}

Theorem~\ref{thm:computed-state-error} does not assume that the coarse
coordinate minimizes the reduced energy.

\begin{cor}[Uniform-depth stationary solutions]
\label{cor:stationary-root-rate}
There is a fixed depth \(\ell_*\ge r_P\), depending on contrast and
the stated geometry and operator constants, such that uniform patches
with \(\ell\ge\ell_*\) satisfy
\begin{equation}
 \|\nabla(M_h^\ell(v_H)-u_h)\|_\Omega\le3\frac{C_\Pi}{m}\rho_H(v_H)+C(1+\ell)^d e^{-b\ell/2}\frac{L_F}{m} \label{eq:uniform-reduced-residual-rate}
\end{equation}
for every coarse coordinate, with \(b=|\log\theta|/Q\).
In particular, every tangent-stationary solution satisfies this bound
with \(\rho_H(v_H)=0\). The depth and constants are independent of
state and load amplitudes and the fine-to-coarse mesh ratio.
\end{cor}

\begin{proof}
The counting argument that gives \eqref{eq:uniform-weighted-overlap},
applied with unit weights, yields the first bound below; together with
\eqref{eq:uniform-weighted-overlap} this gives the second:
\begin{align}
 \mathcal N_0(\mathcal D)&\le C_{\rm ol}(1+\ell)^d \label{eq:uniform-unweighted-overlap}\\
 \beta_{\mathcal D}\varepsilon_{\mathcal D}&\le C(1+\ell)^d\theta^{n_\ell/2}\le C'(1+\ell)^d e^{-b\ell/2} \label{eq:uniform-absorption-factor}
\end{align}
by the definitions \eqref{eq:tangent-stability-constant},
\eqref{eq:transfer-quantities} and \eqref{eq:residual-localization-factor}
of \(C_T\), \(\beta_{\mathcal D}\) and \(\varepsilon_{\mathcal D}\). The right-hand side tends to zero,
so it is at most \(1/2\) beyond a fixed depth. Theorem~\ref{thm:computed-state-error}
then proves the assertion. The constants \(C,C'\) may depend on contrast
but not on \(\ell\), the load, or mesh resolution.
\end{proof}

The additional factor \((1+\ell)^{d/2}\), compared with the
energy-minimizer gradient bound \eqref{eq:exponential-gradient-form},
comes from the tangent stability estimate
(Lemma~\ref{lem:localized-tangent-stability}).
It is not claimed sharp. Logarithmic oversampling still gives order
\(H^s\) localization with \(\ell\ge\max\{\ell_*,c\log(1/H)\}\)
and \(c>2s/b\), after scaling \(H\) by a fixed domain length.

The theorem asserts accuracy of every stationary root under sufficient
oversampling, without identifying the roots as global minimizers or
proving uniqueness or convergence of an outer solver. At least one root
exists by the existence of an energy minimizer
(Theorem~\ref{thm:stationary-error}) and the differentiability of the
reconstruction (Lemma~\ref{lem:reconstruction-regularity}).
If every patch fills the domain, \(\varepsilon_{\mathcal D}=0\), and
\eqref{eq:computed-state-closed-bound} recovers the ideal residual bound;
every ideal stationary state is therefore the fine solution. The next
subsection includes defects from approximate patch solves and test maps.

\subsection{Inexact Patch Solves and Coarse Tests}
\label{subsec:inexact-solves}

The preceding estimates (Theorem~\ref{thm:computed-state-error} and
Corollary~\ref{cor:stationary-root-rate}) use exact patch maps. We now separate their
localization error from the error of evaluating them. The exact patch
maps remain comparison objects; the computed corrections below need
not satisfy their equations exactly. Fix a coarse
coordinate \(v_H\) and let \(\widehat q_i\in V_f^0(D_i)\) be
approximate patch corrections that retain the global kernel constraints.
Define their residual norms and the assembled state by
\begin{align}
 \sigma_i&=\sup_{0\ne z_f\in V_f^0(D_i)}\frac{|F[z_f]-\mathcal A(v_H+\widehat q_i)[z_f]|}{\|\nabla z_f\|_{D_i}} \label{eq:inexact-patch-residual}\\
 \widehat u&=v_H+(I-\Pi_H)I_h\left(\sum_i\phi_i\widehat q_i\right),\qquad
 \zeta=\frac{\sqrt{C_{\rm asm}}}{m}\left(\sum_i\sigma_i^2\right)^{1/2} \label{eq:inexact-state-definition}
\end{align}
with zero residual norm when the test space is zero dimensional.
These are dual gradient norms on the constrained patch spaces, not
unscaled norms of algebraic residual vectors.

\begin{lem}[Patch residuals control assembly error]
\label{lem:inexact-patch-assembly}
Under Assumption~\ref{ass:standing},
\begin{equation}
 \|\nabla(\widehat u-M_h^{\mathcal D}(v_H))\|_\Omega\le\zeta \label{eq:inexact-assembly-bound}
\end{equation}
\end{lem}

\begin{proof}
Put \(e_i=\widehat q_i-q_i(v_H)\in V_f^0(D_i)\). Subtract the
exact patch equation and test with \(e_i\). Strong monotonicity gives
\begin{equation}
 m\|\nabla e_i\|_{D_i}^2\le(\mathcal A(v_H+\widehat q_i)-F)[e_i]\le\sigma_i\|\nabla e_i\|_{D_i} \label{eq:inexact-patch-error}
\end{equation}
and hence \(\|\nabla e_i\|_{D_i}\le\sigma_i/m\).
Apply Lemma~\ref{lem:pou-assembly} to the family \((e_i)_i\).
Since \(e_i\) vanishes outside \(D_i\), we have
\(\|\nabla e_i\|_{N(\omega_i)}\le\|\nabla e_i\|_{D_i}\) for every patch.
\end{proof}

Lemma~\ref{lem:inexact-patch-assembly} requires no smallness of the patch
residuals. The inclusion \(N(\omega_i)\subset D_i\) is not used in this
lemma, so it also applies to shallow patches.

Let \(\widehat T:V_H\to V_h\) be any linear coarse test map satisfying
\begin{equation}
 \Pi_H\widehat T\tau_H=\tau_H,\qquad
 \|\nabla\widehat T\tau_H\|_\Omega\le\widehat C_T\|\nabla\tau_H\|_\Omega \quad\forall\tau_H\in V_H \label{eq:inexact-test-section}
\end{equation}
where a valid \(\widehat C_T\ge1\) is used. Define
\begin{equation}
 \widehat\rho_H=\sup_{0\ne\tau_H\in V_H}\frac{|(\mathcal A(\widehat u)-F)[\widehat T\tau_H]|}{\|\nabla\tau_H\|_\Omega},\qquad
 \widehat\beta=\chi C_\Pi\widehat C_T \label{eq:inexact-coarse-residual}
\end{equation}
with the same convention for a zero-dimensional space. The test map
need not be the derivative of the approximate reconstruction: its
section property and stability are the properties used below.

\begin{prop}[Error with patch and coarse residuals]
\label{prop:inexact-computed-state}
Let \(\varepsilon_{\mathcal D}\) be defined in
\eqref{eq:residual-localization-factor}. If
\(\widehat\beta\varepsilon_{\mathcal D}<1\), then
\begin{align}
 \|\nabla(\widehat u-u_h)\|_\Omega
 &\le\frac{(1+\varepsilon_{\mathcal D})C_\Pi}{m(1-\widehat\beta\varepsilon_{\mathcal D})}\widehat\rho_H
 +\frac{1+\widehat\beta}{1-\widehat\beta\varepsilon_{\mathcal D}}
 \left(\varepsilon_{\mathcal D}\frac{L_F}{m}+\zeta\right) \label{eq:inexact-computed-error}
\end{align}
under Assumption~\ref{ass:standing}. If the right-hand side
is denoted by \(\widehat S\), the energy error also satisfies
\begin{equation}
 0\le\mathcal E(\widehat u)-\mathcal E(u_h)\le\frac M2\widehat S^2 \label{eq:inexact-computed-energy}
\end{equation}
\end{prop}

\begin{proof}
Write \(U=M_h(v_H)\) and \(X=\|\nabla(U-u_h)\|_\Omega\).
The ideal residual annihilates \(V_f\), and
\(U-u_h-\widehat T\Pi_H(U-u_h)\in V_f\) by
\eqref{eq:inexact-test-section}. Repeating the test in
\eqref{eq:residual-transfer-identity} with \(\widehat u,\widehat T\)
gives
\begin{equation}
 X\le\frac{C_\Pi}{m}\widehat\rho_H+\widehat\beta\|\nabla(\widehat u-U)\|_\Omega \label{eq:inexact-transfer}
\end{equation}
whereas Lemma~\ref{lem:inexact-patch-assembly} and
\eqref{eq:computed-coordinate-localization} give
\begin{equation}
 \|\nabla(\widehat u-U)\|_\Omega\le\varepsilon_{\mathcal D}(X+L_F/m)+\zeta \label{eq:inexact-localization}
\end{equation}
Absorb \(\widehat\beta\varepsilon_{\mathcal D}X\) in
\eqref{eq:inexact-transfer} and use
\(\|\nabla(\widehat u-u_h)\|_\Omega\le
(1+\varepsilon_{\mathcal D})X+\varepsilon_{\mathcal D}L_F/m+\zeta\).
Collecting coefficients proves \eqref{eq:inexact-computed-error};
\eqref{eq:energy-gap-equivalence} gives the energy estimate.
\end{proof}

This estimate separates localization, patch solution error and the coarse
residual, and it does not assume global minimality. Approximate tangent solves are covered if their assembled test
map retains \eqref{eq:inexact-test-section} with a justified stability
bound. No estimate of the difference between exact and approximate
constitutive Hessians is needed. Violations of the kernel constraints
or the test-map section property are not included. Stopping tolerances alone do not certify these dual norms;
Section~\ref{sec:supp-certificates} of the supplementary material
gives their matrix evaluation, assembly bounds that use the overlap of
the patch supports, and a numerical evaluation.


\section{Numerical Experiments}
\label{sec:numerical-experiments}

The experiments address three questions: how localization depends on
the projection and patch depth, how it responds to contrast and load,
and whether patch and outer-solver errors affect the observed accuracy.
We use the saturated law of Example~\ref{example:saturated-law} and
the same partition-of-unity construction as the analysis.
Section~\ref{subsec:p-numerics} of the supplementary material adds
the quartic law. The experiments address
accuracy relative to a discrete fine solution; fine-reference
refinement is reported separately. The numerical results illustrate
the estimates without determining their constants or the sufficient
oversampling threshold in Theorem~\ref{thm:computed-state-error}.

\subsection{Problem, Discretization, and Error Measures}
\label{subsec:numerical-setup}

Let \(\Omega=(0,1)^2\), impose homogeneous Dirichlet conditions, and
set \(\alpha=4\), \(\beta=1\) in
\eqref{eq:saturated-law}--\eqref{eq:saturated-potential}. With
coefficient contrast \(\kappa_a=a_+/a_-\ge1\) and constant load \(f_0>0\),
the potential and load are
\begin{equation}
 \psi(x,\xi)=\frac{a_\varepsilon(x)}2\big[5|\xi|^2-4\log(1+|\xi|^2)\big],\qquad F(w)=f_0\int_\Omega w\,dx \label{eq:numerical-model}
\end{equation}
where \(a_\varepsilon\) takes values \(1\) and \(\kappa_a\) on a
fixed \(12\times12\) array of square pixels. The baseline has
\(\kappa_a=10\) and \(f_0=20\). The binary array is generated by
NumPy's default random generator with seed 20260909, using independent
equiprobable binary entries assigned row by row from bottom to top. The Hessian bounds are \(m=1\) and
\(M=5.5\kappa_a\), independently of the load.

A grid with \(n\) subdivisions in each coordinate direction is
triangulated by dividing each square along its lower-left to
upper-right diagonal. Write \(n_H\) and \(n_h\) for the coarse
and fine subdivision counts. The quantities \(H=1/n_H\) and
\(h=1/n_h\) used in the numerical data are square edge lengths; the
corresponding triangle diameters are \(\sqrt2H\) and \(\sqrt2h\).
We use nested conforming piecewise-affine spaces. The baseline coarse
grid has \(n_H=6\), with 25 interior coarse unknowns. Fine grids have
\(n_h=24,48,96\), and a further reference solve uses \(n_h=192\).
Every coefficient interface is resolved on every fine grid. Since the
coefficient and discrete gradients are constant on each fine triangle,
element integration of the potential, flux, and Hessian is exact.

For the local weak projector, each interior coarse vertex selects its
first incident triangle in the mesh ordering: squares are ordered row
by row, with the lower triangle before the upper triangle. The integral
degrees of freedom are those of \eqref{eq:projection-definition} and
are evaluated exactly. The comparison projectors are nodal interpolation
and the global \(L^2\) projection onto \(V_H\). Each choice defines
its own fine kernel and its own coordinate \(\Pi_Hu_h\). On a patch,
all global kernel constraints are restricted to the supported fine
unknowns; dependent constraints are removed before solving. This is the
supported global kernel \(V_f^0(D_i)\), rather than a separately
defined patch kernel.

Depth zero is the support of a coarse nodal hat. Each layer
adds coarse triangles sharing a vertex with the current patch. All
coarse nodal hats, including boundary hats, enter the partition of unity.
The corrections solve \eqref{eq:patch-equation} with the same restricted
load and are assembled using \eqref{eq:pou-reconstruction}.

For a state \(w_h\in V_h\), define the relative gradient and energy
errors by
\begin{align}
 e_\nabla(w_h)&=\frac{\|\nabla(w_h-u_h)\|_\Omega}{\|\nabla u_h\|_\Omega} \label{eq:numerical-gradient-error}\\
 e_{\mathcal E}(w_h)&=\frac{\mathcal B(w_h,u_h)}{\int_\Omega\psi(x,\nabla u_h)\,dx} \label{eq:numerical-energy-error}
\end{align}
where \(\mathcal B\) is the Bregman distance from
Section~\ref{sec:energy-setting}. Its numerator is evaluated directly,
which avoids subtracting nearby signed total energies. For an exact
fine solution it equals the total-energy gap by
\eqref{eq:energy-bregman}. The computed fine references have relative
residual in the dual gradient norm at most \(1.39\times10^{-11}\) in all
experiments reported here. The positive internal energy in the denominator fixes the
normalization used throughout the numerical figures.

We call \(M_h^\ell(\Pi_Hu_h)\) the \emph{comparator reconstruction}.
It uses the fine reference to isolate the error caused by localization.
A \emph{computed state} instead uses a coarse coordinate obtained from
an outer solve of the reduced energy. These states are distinguished in
the results below; no computed state is assumed to be a global minimizer.

\subsection{Localization and Fine-Reference Refinement}
\label{subsec:numerical-refinement}

\begin{figure}[!htb]
\centering
\includegraphics[width=\textwidth]{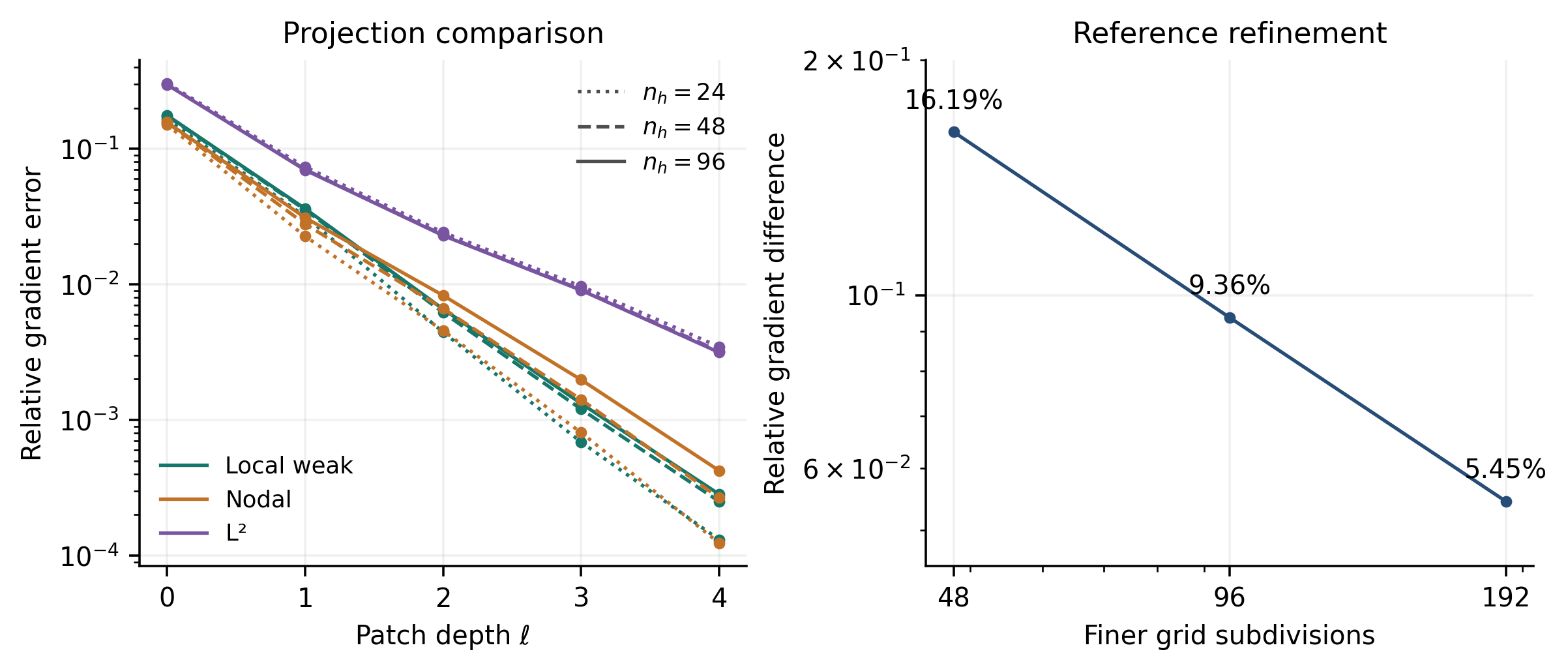}
\caption{Baseline problem with contrast 10 and load 20. Left: relative
 gradient error \(e_\nabla\) of comparator reconstructions
 \(M_h^\ell(\Pi_Hu_h)\) for local weak, nodal, and global \(L^2\)
 projections, on fine grids with 24, 48, and 96 subdivisions per edge;
 the coarse grid has six subdivisions per edge. Right: successive
 fine-reference differences after exact affine prolongation to the
 finer mesh, divided by the finer reference's gradient norm. These
 differences measure reference refinement, not localization.}
\label{fig:numerical-refinement}
\end{figure}

Figure~\ref{fig:numerical-refinement} compares the three projectors at
fixed \(n_H=6\) for three fine meshes, each with its own fine reference
solution \(u_h\). The comparator errors decrease rapidly with patch
depth on all three fine meshes. Depths zero through four contain no
full-domain patches. At depth six, which is omitted from the figure, 19 of the 49 patches are the full
domain. The plotted range therefore avoids exact global recovery as
an explanation of the observed decay, although the larger patches
already occupy a substantial part of this small domain.

At depth four on the fine grid with \(n_h=96\), the relative gradient
errors are \(2.84\times10^{-4}\) for the local weak projection,
\(4.23\times10^{-4}\) for nodal interpolation, and
\(3.15\times10^{-3}\) for the global \(L^2\) projection. The local
weak projector has smaller errors than nodal interpolation at depths
two through four on this mesh, whereas nodal interpolation remains
better at depths zero and one. On the fine grid with \(n_h=24\), nodal
interpolation has a slightly smaller depth-four error. Thus the
projector ordering depends on both depth and fine resolution; the
results support the local weak choice without establishing a universal
ranking.

As a separate diagnostic, the discrete operator norm of \(\Pi_H\)
from the fine gradient norm to the coarse gradient norm was computed
at fixed \(n_H=6\). As \(n_h\) increases through \(12,24,48,96\),
the norms are \(1.446,1.701,1.831,1.868\) for the local weak
projection, \(1.390,1.689,1.924,2.128\) for nodal interpolation,
and \(1.204,1.299,1.326,1.334\) for the global \(L^2\) projection.
These finite-mesh measurements do not prove uniform stability or its
failure. They also show that a smaller global operator norm need not
give a smaller localization error.

The successive reference differences for
\(24\to48\), \(48\to96\), and \(96\to192\) subdivisions are
\(16.19\%\), \(9.36\%\), and \(5.45\%\), respectively.
They decrease under refinement, but the finest reference is not a
certified continuum solution. In particular, the small comparator
errors in Figure~\ref{fig:numerical-refinement} quantify accuracy
relative to \(u_h\), as in the theory, and must not be interpreted
as total error relative to the continuous solution.

\subsection{Coarse Resolution, Contrast, and Nonlinear Load}
\label{subsec:numerical-parameters}

To examine coarse resolution, we keep the fine grid with \(n_h=48\),
the pixel field, and the load fixed and use \(n_H=4,6,12\). The fine
reference is identical in these three runs. In addition to
\(e_\nabla\), we report the error normalized by the ideal correction:
\begin{equation}
 e_q^\ell=\frac{\|\nabla(M_h^\ell(\Pi_Hu_h)-u_h)\|_\Omega}{\|\nabla(I-\Pi_H)u_h\|_\Omega} \label{eq:numerical-correction-normalization}
\end{equation}
which corresponds to the amplitude in the reconstruction estimate
\eqref{eq:weighted-correction-reconstruction}.
Table~\ref{tab:numerical-coarse} shows rapid decay for all three coarse
grids. A layer has a different physical size on each grid, and the
\(n_H=4\), depth-four case has 13 full-domain patches out of 25.
Consequently the table is a fixed-reference localization comparison,
not an estimate of a coarse-mesh convergence rate or computational
complexity.

\begin{table}[!htb]
\centering
\begin{tabular}{ccccc}
\toprule
\(n_H\) & \(e_\nabla\), \(\ell=2\) & \(e_\nabla\), \(\ell=4\) & \(e_q^2\) & \(e_q^4\) \\
\midrule
4  & \(1.08\times10^{-2}\) & \(1.75\times10^{-4}\) & \(1.89\times10^{-2}\) & \(3.07\times10^{-4}\) \\
6  & \(6.21\times10^{-3}\) & \(2.51\times10^{-4}\) & \(1.14\times10^{-2}\) & \(4.59\times10^{-4}\) \\
12 & \(6.17\times10^{-3}\) & \(2.60\times10^{-4}\) & \(1.74\times10^{-2}\) & \(7.33\times10^{-4}\) \\
\bottomrule
\end{tabular}
\caption{Comparator errors for local weak projection at fixed fine
 resolution \(n_h=48\), contrast 10, and load 20. The coarse subdivision
 count is \(n_H\); \(e_\nabla\) and \(e_q^\ell\) are defined in
 \eqref{eq:numerical-gradient-error} and
 \eqref{eq:numerical-correction-normalization}. Only the \(n_H=4\),
 \(\ell=4\) entry includes full-domain patches.}
\label{tab:numerical-coarse}
\end{table}

\begin{figure}[!htb]
\centering
\includegraphics[width=\textwidth]{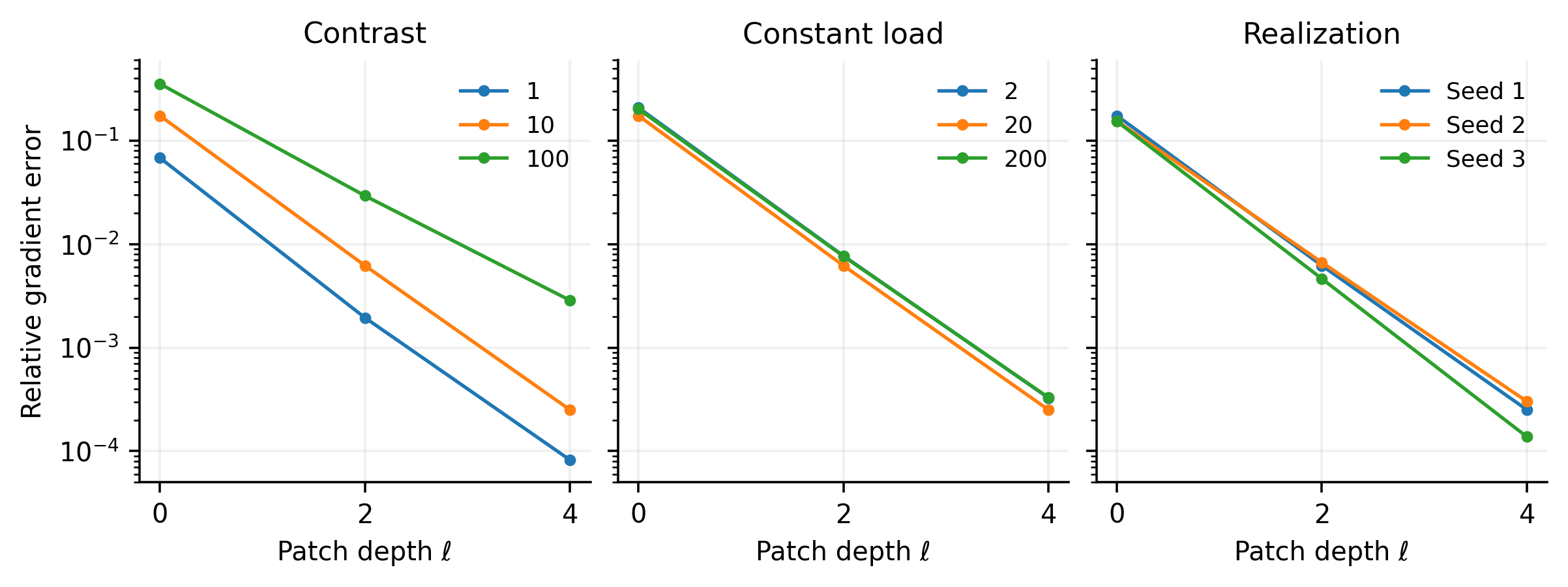}
\caption{Relative gradient errors of local weak comparator
 reconstructions at depths \(\ell=0,2,4\), with \(n_H=6\) and
 \(n_h=48\). Each panel varies one factor from the baseline contrast
 10, load 20, and coefficient seed 20260909. Seed 1, Seed 2, and Seed 3
 denote 20260909, 20260910, and 20260911. Each error uses the fine reference for its
 own coefficient and load.}
\label{fig:numerical-parameters}
\end{figure}

Figure~\ref{fig:numerical-parameters} varies one factor at a time on
the \(n_H=6\), \(n_h=48\) mesh pair, using local weak projection.
The contrast values are \(1,10,100\); the constant loads are
\(2,20,200\); and the three coefficient realizations use seeds
20260909, 20260910, and 20260911. Each changed problem is compared with
its own fine solution. Increasing the contrast from 10 to 100 increases
the depth-four relative gradient error from
\(2.51\times10^{-4}\) to \(2.87\times10^{-3}\) and slows the observed
decay. This is consistent with the contrast dependence of the theoretical
constants; no contrast-independent rate is claimed.

Load changes and the two additional realizations have smaller effects
in this sample. For loads \(2,20,200\), the median magnitudes of the
fine-solution gradient are \(0.094,0.511,2.261\), respectively. Thus
the load study samples gradients below, around, and above the saturation
scale \(|\xi|=1\) of the chosen law. The depth-four errors for loads
2 and 200 are \(3.27\times10^{-4}\) and \(3.30\times10^{-4}\).
The two additional realizations give \(3.02\times10^{-4}\) and
\(1.37\times10^{-4}\). These tests examine a range of states for one
bounded-Hessian potential; three realizations do not constitute a
statistical robustness study.

\subsection{Energy Errors and Computed States}
\label{subsec:numerical-computed-states}

\begin{figure}[!htb]
\centering
\includegraphics[width=\textwidth]{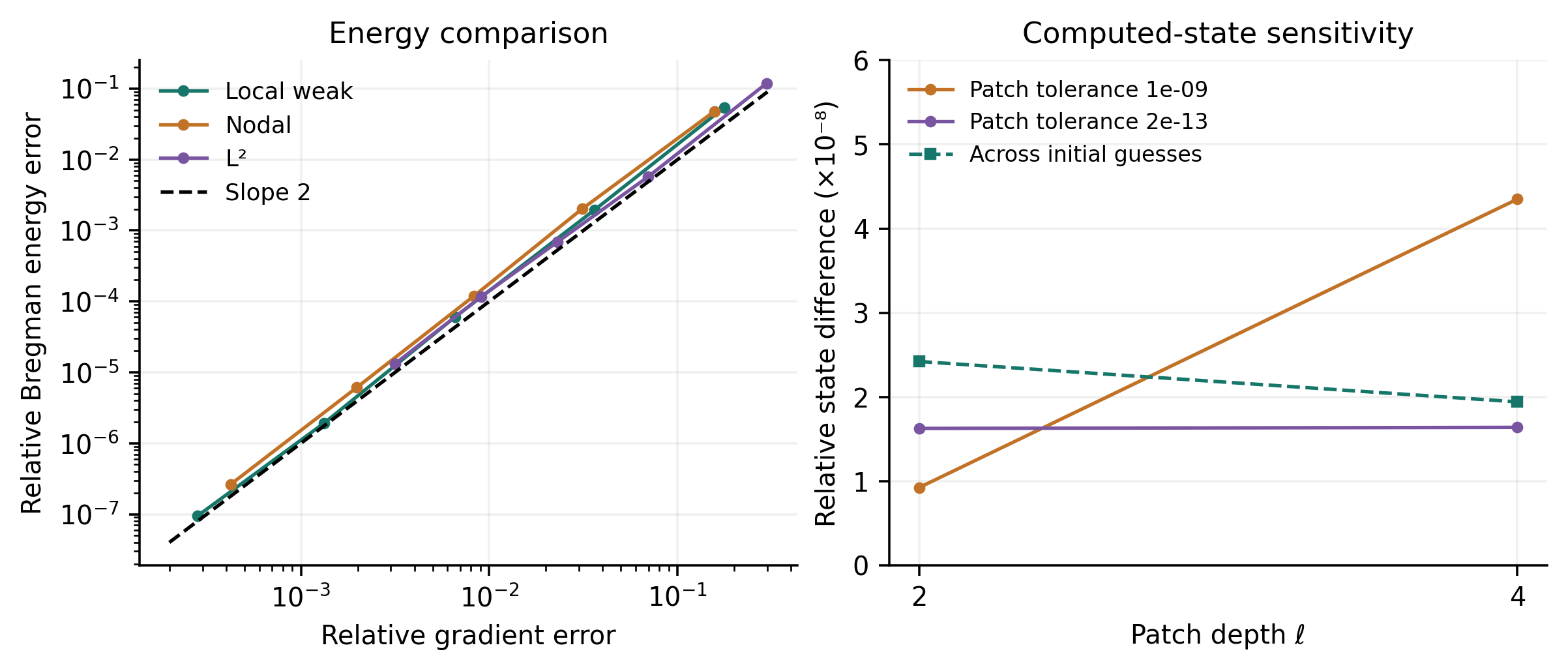}
\caption{Left: relative Bregman energy error versus relative gradient
 error for baseline comparator reconstructions with \(n_h=96\),
 \(n_H=6\), and depths zero through four. The dashed line
 has slope two and is a reference slope, not a theoretical bound.
 Right: direct gradient-norm differences between computed local weak
 states with \(n_h=48\), divided by \(\|\nabla u_h\|_\Omega\).
 Solid curves compare zero-start states at the indicated patch tolerance
 with the default \(2\times10^{-11}\); the dashed curve gives the
 largest difference among three initial guesses at the default tolerance.}
\label{fig:numerical-energy-solver}
\end{figure}

The left panel of Figure~\ref{fig:numerical-energy-solver} plots the
relative Bregman energy error against the relative gradient error for
the comparator states on the fine grid with \(n_h=96\). The curves
are consistent with the quadratic relation
\eqref{eq:energy-gap-equivalence}. Across the
60 comparator reconstructions and 10 outer solves of this section, all 70
states satisfy
the lower and upper Bregman bounds with a numerical allowance of
\(0.01\%\). The observed unnormalized ratio
\(\mathcal B(w_h,u_h)/\|\nabla(w_h-u_h)\|_\Omega^2\) ranges from
1.128 to 7.111 across the different contrasts, loads, and meshes.
This check concerns the energy-distance equivalence, not certification
of global minimality on the localized manifold.

Patch problems are solved by constrained Newton iteration with a
line search. If \(J_i\) is the current patch Hessian matrix,
\(q_i\) its correction vector, and \(s_i\) the constrained Newton
step, the stopping rule is
\begin{equation}
 (s_i^TJ_is_i)^{1/2}<\tau_{\rm patch}\big[1+(q_i^TJ_iq_i)^{1/2}\big],\qquad \tau_{\rm patch}=2\times10^{-11} \label{eq:numerical-patch-tolerance}
\end{equation}
with at most 40 Newton iterations. Tangents are obtained by differentiating
the same constrained equations, as in \eqref{eq:patch-tangent}.
The outer iteration uses the Broyden--Fletcher--Goldfarb--Shanno (BFGS)
method for the reduced energy with its analytical gradient. Coarse
coordinates are scaled by a Cholesky factor of the coarse Laplace
stiffness matrix, so the Euclidean norm of the scaled gradient equals
\(\rho_H\) in \eqref{eq:reduced-residual-norm}. For the outer solves, BFGS uses its default infinity-norm
stopping test on the scaled gradient with tolerance \(2\times10^{-9}\)
and at most 150 iterations. This stopping norm differs from the Euclidean
norm defining \(\rho_H\). Reported residuals are computed separately
as \(\rho_H/L_F\).

Table~\ref{tab:numerical-computed} gives the zero-start computed states
for local weak projection on the fine grid with \(n_h=48\). Their
errors closely follow the comparator errors, while the reduced residuals
are much smaller. The finest-grid projector comparison in
Figure~\ref{fig:numerical-refinement} uses comparator reconstructions.

\begin{table}[!htb]
\centering
\begin{tabular}{ccccc}
\toprule
\(\ell\) & Comparator \(e_\nabla\) & Computed \(e_\nabla\) & Computed \(e_{\mathcal E}\) & \(\rho_H/L_F\) \\
\midrule
2 & \(6.206\times10^{-3}\) & \(6.186\times10^{-3}\) & \(5.363\times10^{-5}\) & \(4.852\times10^{-9}\) \\
4 & \(2.508\times10^{-4}\) & \(2.508\times10^{-4}\) & \(7.543\times10^{-8}\) & \(3.610\times10^{-9}\) \\
\bottomrule
\end{tabular}
\caption{Comparator and zero-start computed-state errors for local
 weak projection, \(n_H=6\), \(n_h=48\), contrast 10, and load 20.
 Error normalizations are \eqref{eq:numerical-gradient-error} and
 \eqref{eq:numerical-energy-error}; \(\rho_H/L_F\) is the reduced
 residual divided by the fine-space dual gradient norm of the load.
 Both outer solves carry a precision-loss flag; their measured
 residuals are reported explicitly.}
\label{tab:numerical-computed}
\end{table}

At depths two and four, the outer solve is repeated from zero, a random
coarse vector, and the reference-informed coordinate \(\Pi_Hu_h\).
The random vector has independent standard normal entries multiplied
by \(L_F\) in the gradient-scaled coordinates, with seed 987. The last
start is included only as a diagnostic, since it requires the fine solution. Direct comparisons of the saved
states give a maximum relative gradient difference below
\(2.43\times10^{-8}\) across these starts. Repeating the zero-start
runs with patch tolerances \(10^{-9}\) and \(2\times10^{-13}\)
changes the computed states by at most \(4.36\times10^{-8}\).
These differences are well below the measured localization errors.
Across the ten outer solves, the largest relative reduced residual is
\(2.64\times10^{-8}\), and the largest residual of the scaled kernel
constraints is \(1.52\times10^{-16}\).

Nine of these ten outer runs return BFGS's precision-loss flag at the
strict requested tolerance. We retain these flags and assess the states
using the measured residuals and direct comparisons. The agreement
across initial guesses does not establish uniqueness, global minimality,
or convergence of the outer method. Nor do the experiments certify
that the sufficient theoretical oversampling threshold has been reached.


\section{Discussion and Outlook}
\label{sec:discussion}

Localization supports two accuracy arguments: energy comparison for
manifold minimizers and residual transfer for computed states. The
matrix certificates of the supplementary material connect the latter to
numerical evaluation, while
the full fine-space residual bounds total discrete error independently
of the analytical localization constants. The remaining questions
concern the projection, constitutive growth and evolution problems.

\subsection{Choice of Projection and Numerical Comparisons}
\label{subsec:projection-comparisons}

The main theorem uses the local weak projection
\eqref{eq:projection-definition}. Its local degrees of freedom preserve
coarse functions exactly and provide uniform approximation estimates.
The selection of \(K_z\) may affect actual errors even when the
theorem constants are bounded uniformly over the allowed choices. No
comparison theorem ranking those choices is claimed.

The comparisons in Section~\ref{sec:numerical-experiments} show that
a smaller global gradient operator norm need not give better
localization. This supports examining locality and kernel approximation
together with stability. The observed projection ranking and decay
rates remain problem-dependent.

The present cutoff proof depends on a projection with a fixed local
dependence radius. It therefore does not apply directly to the global
\(L^2\) projection. For that choice, the missing ingredient is a
local operation that preserves its moment constraints with the required
annular stability. Nodal interpolation likewise requires kernel estimates
uniform in the fine resolution. Good accuracy for a fixed pair of meshes
does not by itself establish those estimates. The local weak projection
used here supplies the properties needed by both localization and assembly
without an additional fine-resolution threshold.

\subsection{Growth and Remaining Constitutive Regimes}
\label{subsec:unbounded-growth}

The bounds of this paper use the state-independent constants in
\eqref{eq:hessian-bounds}. For laws with growth, such as the quartic
diffusion law, the upper constant is unavailable.
Section~\ref{sec:p-structure} of the supplementary material replaces
the global ratio \(M/m\) by scale-invariant local quantities of a weight
determined by the fine solution, for potentials with \(p\)-structure,
\(p\ge2\). The price is that the constants are quantities of the fine
solution. Their uniformity in the load and the mesh is only partly
settled; see Section~\ref{sec:p-discussion} there. Degenerate laws such
as the \(p\)-Laplacian and laws with \(p<2\) are not covered.

The energy comparison in Theorem~\ref{thm:energy-comparison} separates
these issues from approximation on a nonlinear trial set: comparison of
energy values follows from minimality, while conversion to gradient
norms and the explicit rate use the constitutive bounds.

\subsection{Time-Dependent Energy Problems}
\label{subsec:evolution-outlook}

The energy setting provides a natural route to an evolution problem
with an \(L^2\) time derivative:
\begin{equation}
 (\partial_tu_h,w)_\Omega+\mathcal A(u_h)[w]=F(t)[w] \quad \forall w\in V_h \label{eq:evolution-equation}
\end{equation}
where \(F(t)\) is a time-dependent load. For a time-independent load
and a differentiable trajectory, testing with \(\partial_tu_h\) gives
\begin{equation}
 \frac{d}{dt}\mathcal E(u_h(t))=-\|\partial_tu_h(t)\|_\Omega^2 \label{eq:energy-dissipation}
\end{equation}
as the energy dissipation identity. A backward-Euler step of size
\(\tau>0\) is associated with the incremental energy
\begin{equation}
 \mathcal E_n(v)=\frac{1}{2\tau}\|v-u_h^{n-1}\|_\Omega^2+\int_\Omega\psi(x,\nabla v)\,dx-F(t_n)[v] \label{eq:incremental-energy}
\end{equation}
where \(u_h^{n-1}\) is the preceding fine state and \(t_n\) the new
time. Reusing a stationary reconstruction, as explored in the companion
paper~\cite{larson-persson2026}, and eliminating fine scales for this incremental energy are
different choices. The latter introduces a mass term into the secant
form and the history into the load. A time-dependent analysis must control the mass contribution to
localization, the propagation of errors through the history,
and the time-discretization error. The stationary energy framework
provides a starting point for these questions; the estimates proved in
this paper concern the stationary problem.
\appendix

\section{Buffered Cutoff Construction}
\label{sec:buffered-cutoff}

For the local dual projection \eqref{eq:projection-definition}, the
cutoff in Lemma~\ref{lem:corrected-cutoff} needs only three coarse
layers. The proof uses the vanishing local moments of a fine-kernel
function before estimating the projection correction.

\subsection{Cutoff and Exact Support}
\label{subsec:three-layer-cutoff}

Write \(S_j=N^j(S)\) and choose
\begin{equation}
 Q=3,\qquad A=S_3\setminus S \label{eq:buffer-widths}
\end{equation}
with all differences interpreted elementwise. Give every coarse vertex
belonging to \(S_1\) cutoff value one and every remaining vertex
value zero. The resulting continuous coarse affine function \(\eta_H\)
includes boundary vertices and satisfies
\begin{equation}
 0\le\eta_H\le1,\qquad \eta_H|_{S_1}=1,\qquad
 \operatorname{supp}\eta_H\subset\overline{S_2},\qquad
 \|\nabla\eta_H\|_{L^\infty(T)}\le CH^{-1} \label{eq:geometric-cutoff}
\end{equation}
on every coarse element. Its nonconstant elements form the set
\(\mathcal T_\eta\), for which graph distances give
\begin{equation}
 \mathcal T_\eta\subset S_2\setminus S_1,\qquad
 N(\mathcal T_\eta)\subset S_3\setminus S=A \label{eq:transition-buffer-inclusion}
\end{equation}
since every transition element has distance exactly two from \(S\).

For \(w_f\in V_f\), set
\begin{equation}
 v=I_h(\eta_Hw_f),\qquad z_f=(I-\Pi_H)v \label{eq:corrected-cutoff-definition}
\end{equation}
so \(z_f\in V_f\). Both functions have zero physical trace.
Let \(c_z=\ell_z(v)\). Linear independence of the interior coarse hats
implies \(\ell_z(w_f)=0\) for every interior vertex \(z\).
Writing \(\eta_z=\eta_H(z)\), linearity and \(I_hw_f=w_f\) give
\begin{equation}
 c_z=\ell_z\big(I_h((\eta_H-\eta_z)w_f)\big),\qquad
 c_z=0\quad\hbox{if }K_z\notin\mathcal T_\eta \label{eq:cutoff-moment-cancellation}
\end{equation}
because \(\eta_H\) is constant on such a selected simplex.
Every hat with a remaining coefficient is supported in \(N(K_z)\).
Consequently
\begin{equation}
 \operatorname{supp}\Pi_Hv\subset
 \overline{N(\mathcal T_\eta)}\subset\overline A \label{eq:transition-influence}
\end{equation}
and this support contains no element of \(S\). Since \(v=w_f\) on
\(S_1\), we obtain the exact statements
\begin{align}
 z_f|_S&=w_f|_S \label{eq:buffer-inner-matching}\\
 \operatorname{supp}z_f&\subset\overline{S_3},\qquad z_f\in V_f^0(S_3) \label{eq:buffer-outer-support}
\end{align}
including all artificial-boundary constraints by global conformity.
Support in the closure is precisely the definition of \(V_f^0(S_3)\);
no additional empty layer is required.

\subsection{Annular Stability}
\label{subsec:three-layer-stability}

First, each selected simplex has a Poincar\'e estimate of its own
for fine-kernel functions. The dual formula \eqref{eq:local-dual-formula} gives
\begin{equation}
 \int_K\theta_z^K\,dx=1,\qquad
 \|\theta_z^K\|_K^2=\frac{(d+1)^2}{|K|},\qquad
 \|\theta_z^K-|K|^{-1}\|_K^2=\frac{d(d+2)}{|K|} \label{eq:selected-dual-norms}
\end{equation}
by barycentric integration. Let \(C_{\rm simp}\) be a uniform
constant in \(\|w-w_K\|_K\le C_{\rm simp}H\|\nabla w\|_K\).
When \(K=K_z\) and \(\ell_z(w_f)=0\), subtracting the mean inside
that moment and using orthogonality to constants yields
\begin{align}
 |K|| (w_f)_K |^2&\le d(d+2)\|w_f-(w_f)_K\|_K^2 \label{eq:selected-mean-control}\\
 \|w_f\|_K^2&\le(d+1)^2 C_{\rm simp}^2H^2\|\nabla w_f\|_K^2 \label{eq:selected-simplex-poincare}
\end{align}
without enlarging \(K\).

The local dual norm, product interpolation, and
\eqref{eq:selected-simplex-poincare} therefore imply
\begin{equation}
 |c_z|\le C H^{-d/2}\|w_f\|_{K_z}
       \le C H^{1-d/2}\|\nabla w_f\|_{K_z}
       \quad\hbox{if }K_z\in\mathcal T_\eta \label{eq:transition-coefficient-bound}
\end{equation}
with a constant independent of the fine refinement. Coarse hat scaling,
bounded overlap of hats, and the fact that a simplex can be selected
by at most \(d+1\) of its vertices give a uniform constant \(C_{\rm corr}\)
such that
\begin{equation}
 \|\nabla\Pi_Hv\|_\Omega^2
 \le C_{\rm corr}\|\nabla w_f\|_{\mathcal T_\eta}^2 \label{eq:transition-retraction-bound}
\end{equation}
for every seed set. All nonzero coefficients were confined to
\(\mathcal T_\eta\) before this global norm was estimated.

For the uncorrected product, choose \(C_{\rm prod}\ge1\) so that
the squared product estimate of Lemma~\ref{lem:product-interpolation}
and the cutoff gradient bound give
\begin{equation}
 \|\nabla v\|_{\mathcal T_\eta}^2
 \le C_{\rm prod}\big(\|\nabla w_f\|_{\mathcal T_\eta}^2+
                  H^{-2}\|w_f\|_{\mathcal T_\eta}^2\big) \label{eq:transition-product-bound}
\end{equation}
Summation of the kernel Poincar\'e estimate \eqref{eq:kernel-poincare},
with \(r_P=1\), gives a uniform
constant \(C_{\rm star}\) with
\begin{equation}
 H^{-2}\|w_f\|_{\mathcal T_\eta}^2
 \le C_{\rm star}\|\nabla w_f\|_{N(\mathcal T_\eta)}^2
 \le C_{\rm star}\|\nabla w_f\|_A^2 \label{eq:transition-poincare-bound}
\end{equation}
Here \(C_{\rm star}\) includes the bounded overlap of the one-layer
patches. Outside the transition set, \(v\) is exactly \(w_f\) or zero,
element by element. Thus
\begin{align}
 \|\nabla v\|_{A\setminus\mathcal T_\eta}^2
 &\le\|\nabla w_f\|_{A\setminus\mathcal T_\eta}^2 \label{eq:exact-region-bound}\\
 \|\nabla v\|_A^2
 &\le C_{\rm prod}(1+C_{\rm star})\|\nabla w_f\|_A^2 \label{eq:uncorrected-annular-bound}
\end{align}
Combining this estimate with \eqref{eq:transition-retraction-bound}
and the triangle inequality proves
\begin{equation}
 \|\nabla z_f\|_A^2\le C_{\rm cut}\|\nabla w_f\|_A^2,\qquad
 C_{\rm cut}=\big(\sqrt{C_{\rm prod}(1+C_{\rm star})}
                    +\sqrt{C_{\rm corr}}\big)^2 \label{eq:three-layer-cutoff-constant}
\end{equation}
which is \eqref{eq:annulus-cutoff-estimate}. Every gradient on the
right belongs to \(A\); no inner-region energy has been absorbed into
the constant. If \(\mathcal T_\eta\) is empty, all projection
coefficients vanish and the exact-region argument suffices.
Physical boundary vertices cause no extra term, since \(w_f\) has
zero trace and the cut-off function was defined there as well.

\subsection{Scope of the Width Improvement}
\label{subsec:cutoff-width-scope}

The three layers have distinct roles: \(S_1\) protects inner matching,
\(S_2\setminus S_1\) contains the transition, and the third layer
contains both the correction hats and the one-layer Poincar\'e
neighborhoods. The exact coefficient cancellation
\eqref{eq:cutoff-moment-cancellation} confines the projection correction
to the transition neighborhood. Together with the kernel
estimate \eqref{eq:kernel-poincare}, it gives \(Q=3\) and \(r_P=1\), uniformly over the
permitted selections of \(K_z\), without a lower bound on \(H/h\).

If one retains only a kernel Poincar\'e estimate with radius \(r_P\ge1\),
the same coefficient proof uses a cutoff equal to one on \(S_{r_P}\).
The transition elements have graph distance \(r_P+1\) from \(S\),
the correction hats lie in
\(S_{r_P+2}\setminus S\), and its Poincar\'e neighborhoods lie in
\(S_{2r_P+1}\setminus S\). It therefore gives the intermediate bound
\(Q=2r_P+1\) without bundling the projection and Poincar\'e radii.
For a different retraction, the coefficient cancellation must be
established before using this width.

No optimality of \(Q=3\) is asserted. Simply moving the present
cutoff one layer inward would allow the Poincar\'e neighborhoods to
enter \(S\), invalidating the annular estimate \eqref{eq:annulus-cutoff-estimate}. The numerical value
of the decay exponent also depends on \(C_{\rm cut}\), whose
dependence is displayed in \eqref{eq:three-layer-cutoff-constant}.


\paragraph{Acknowledgements.}
Mats G. Larson was supported in part by the Swedish Research Council
under grants 2021-04925 and 2025-05562 and the Swedish Research
Programme Essence. Anna Persson acknowledges support from the Swedish
Research Council under grant 2022-03543.

\paragraph{Use of Generative AI.}
OpenAI Codex and Claude (Anthropic) were used to assist with the
development and checking of mathematical arguments, with literature
searches, with drafting, editing and reorganizing the manuscript, and
with the numerical software and figures. This applies to the main paper
and to the supplementary material, including its \(p\)-structure
analysis. The authors take responsibility for the mathematical content,
references, numerical results, and conclusions of the manuscript.

\paragraph{Software and Reproducibility.}
The numerical experiments use a custom finite element implementation in
Python 3.12 with NumPy 1.26.2, SciPy 1.11.4, and Matplotlib 3.8.2.
The overlap-aware bounds of Section~\ref{subsec:overlap-certificates}
of the supplementary material were evaluated with Python 3.13,
NumPy 2.5.3 and SciPy 1.18.1.
The quartic-law studies of Section~\ref{subsec:p-numerics} and the
nondegeneracy check after Theorem~\ref{thm:p-fixed-mesh} there were run
with the same software versions.
Section~\ref{sec:numerical-experiments} describes the discretization,
solvers, verification checks, and numerical limitations.

\bibliographystyle{plain}
\bibliography{references}

\clearpage
\setcounter{section}{0}
\setcounter{table}{0}
\setcounter{figure}{0}
\renewcommand{\thesection}{S\arabic{section}}
\renewcommand{\thetable}{S\arabic{table}}
\renewcommand{\thefigure}{S\arabic{figure}}
\renewcommand{\theHsection}{S\arabic{section}}
\renewcommand{\theHtable}{S\arabic{table}}
\renewcommand{\theHfigure}{S\arabic{figure}}
\phantomsection
\pdfbookmark[1]{Supplementary Material}{supplementary-material}
\begin{center}
{\LARGE\bfseries Supplementary Material\par}
\end{center}
\vspace{1.5ex}

Sections, results and equations of this supplement carry the prefix S.
Numbers without this prefix refer to the main paper, whose notation and
assumptions (Assumption~\ref{ass:standing}) are used throughout.


\section{Computable Residual Certificates}
\label{sec:supp-certificates}

This section supplements Section~\ref{subsec:inexact-solves} of the main
paper. It gives matrix formulas for the residual norms used there,
assembly bounds that use the overlap of the patch supports, and a
numerical evaluation.

\subsection{Computable Finite-Dimensional Certificates}
\label{subsec:computable-certificates}

To evaluate the residual bounds, we solve linear Riesz problems in
the same finite element spaces. The patch solves must retain the kernel
constraints; coarse and full residuals use their respective gradient
stiffness matrices.
Let \(\mathsf K\) be the fine stiffness matrix for the gradient inner
product, and let \(\mathsf K_i\) be its principal submatrix on the
interior patch degrees of freedom. Let \(\mathsf B_i\) contain
independent rows of the patch kernel constraints. If \(\mathbf r_i\)
is the coefficient vector of the patch residual in
\eqref{eq:inexact-patch-residual}, solve
\begin{equation}
 \begin{pmatrix}\mathsf K_i&\mathsf B_i^T\\\mathsf B_i&0\end{pmatrix}
 \begin{pmatrix}\mathbf w_i\\\lambda_i\end{pmatrix}
 =\begin{pmatrix}\mathbf r_i\\0\end{pmatrix},\qquad
 \sigma_i^2=\mathbf w_i^T\mathsf K_i\mathbf w_i \label{eq:patch-riesz-certificate}
\end{equation}
The multiplier removes the component that annihilates the constrained
space. Thus a nonzero unconstrained residual need not indicate patch
solution error. The norm is independent of the choice or scaling of
independent constraint rows. A zero-dimensional kernel gives
\(\sigma_i=0\).

One can also replace the analytical assembly constant by finite-mesh
operator norms. Let \(\mathsf S_i\) map an arbitrary patch coefficient
vector to the fine coefficient vector of its zero extension multiplied
by \(\phi_i\), interpolated by \(I_h\), and projected by \(I-\Pi_H\).
Define
\begin{equation}
 a_i^2=\lambda_{\max}(\mathsf S_i^T\mathsf K\mathsf S_i,\mathsf K_i),
 \qquad \zeta_{\rm mat}=\frac1m\sum_i a_i\sigma_i \label{eq:matrix-assembly-certificate}
\end{equation}
Valid upper bounds for the \(a_i\) may be substituted in this sum.
Here \(\lambda_{\max}(A,B)\) denotes the largest generalized
eigenvalue of \(Ax=\lambda Bx\), with \(a_i=0\) for an empty
patch space. The eigenvalue is taken over the full patch space, so
it also bounds the operator on the constrained subspace.

\begin{prop}[Matrix residual certificates]
\label{prop:matrix-certificates}
For the constraint-preserving approximate corrections in
Section~\ref{subsec:inexact-solves},
\begin{equation}
 \|\nabla(\widehat u-M_h^{\mathcal D}(v_H))\|_\Omega\le\zeta_{\rm mat} \label{eq:matrix-assembly-error}
\end{equation}
Consequently \(\zeta_{\rm mat}\) may replace \(\zeta\) in
Proposition~\ref{prop:inexact-computed-state}.
If \(\mathbf r\) represents \(\mathcal A(\widehat u)-F\),
\(\mathsf T\) represents the coarse test map, and \(\mathsf K_H\)
is the coarse gradient stiffness matrix, then
\begin{align}
 \widehat\rho_H^2&=(\mathsf T^T\mathbf r)^T\mathsf K_H^{-1}(\mathsf T^T\mathbf r) \label{eq:matrix-coarse-dual-norm}\\
 \|\widehat T\|_{\nabla}^2&=\lambda_{\max}(\mathsf T^T\mathsf K\mathsf T,\mathsf K_H) \label{eq:matrix-test-norm}\\
 \rho^2&=\mathbf r^T\mathsf K^{-1}\mathbf r \label{eq:matrix-full-dual-norm}
\end{align}
where \(\|\widehat T\|_{\nabla}\) is its induced gradient norm and
\(\rho\) is the full residual norm of
Proposition~\ref{prop:full-residual-certificate}.
\end{prop}

\begin{proof}
The constrained Riesz representative in
\eqref{eq:patch-riesz-certificate} satisfies
\(\mathbf z^T\mathbf r_i=\mathbf z^T\mathsf K_i\mathbf w_i\)
for every \(\mathbf z\in\ker\mathsf B_i\). Cauchy--Schwarz and
\(\mathbf z=\mathbf w_i\) prove the norm formula. For the error
\(e_i=\widehat q_i-q_i(v_H)\), the generalized eigenvalue bound and
\eqref{eq:inexact-patch-error} give
\begin{equation}
 \left\|\nabla\sum_i S_i e_i\right\|_\Omega
 \le\sum_i a_i\|\nabla e_i\|_{D_i}
 \le\frac1m\sum_i a_i\sigma_i \label{eq:matrix-assembly-proof}
\end{equation}
where \(S_i\) is the finite element operator represented by
\(\mathsf S_i\). This proves \eqref{eq:matrix-assembly-error}.
The remaining identities are the Riesz norm and Rayleigh quotient
formulas for the fine and coarse gradient inner products.
\end{proof}

In particular, if \(\mathsf K_H=\mathsf L\mathsf L^T\), then
\(\widehat\rho_H=\|\mathsf L^{-1}\mathsf T^T\mathbf r\|_2\).
This is the Euclidean norm of the scaled gradient used in the code;
the optimizer's infinity-norm stopping test is a different quantity.
A computed test norm, enlarged to at least one, supplies
\(\widehat C_T\), provided the section property is verified.
The matrix assembly estimate \eqref{eq:matrix-assembly-error} is a finite-mesh certificate and is not
claimed to have uniform constants or optimal computational cost.

The full-residual bound \(\rho/m\) of
Proposition~\ref{prop:full-residual-certificate} certifies the total gradient error
relative to \(u_h\) without evaluating localization constants or an
oversampling threshold; it costs a fine-grid linear Riesz solve.
The separated bound in \eqref{eq:inexact-computed-error} additionally
requires explicit valid localization and projection constants. Computing
its residual terms alone does not certify its absorption condition.
The formulas above are exact-arithmetic statements. Floating-point
Riesz solves, eigenvalues and constraint checks provide numerical
certificates; rigorous enclosures would additionally require verified
linear algebra and rounding-error bounds.

\subsection{Overlap-Aware Assembly Certificates}
\label{subsec:overlap-certificates}

The summed certificate \eqref{eq:matrix-assembly-certificate} combines
all patch contributions by the triangle inequality. It therefore grows
with the number of patches and does not use the bounded overlap of
their supports. We now replace the triangle inequality by a weighted
Cauchy--Schwarz inequality on each coarse element and use an identity
for the norm of the complementary projection.

For coarse element sets \(W_i\) and positive weights
\(\nu=(\nu_i)_{i\in\mathcal I}\), define the overlap load
\begin{equation}
 \Theta_i^W(\nu)=\max_{K\in W_i}\sum_{j:\,K\in W_j}\nu_j \label{eq:overlap-load}
\end{equation}
where \(K\) ranges over coarse elements. For unit weights, this is the
largest number of sets \(W_j\) that contain a single coarse element of
\(W_i\).

\begin{lem}[Weighted overlap inequality]
\label{lem:weighted-overlap}
Let \(g_i\in V_h\) vanish outside \(W_i\) for every \(i\). Then
\begin{equation}
 \Big\|\nabla\sum_ig_i\Big\|_\Omega^2\le\sum_i\frac{\Theta_i^W(\nu)}{\nu_i}\|\nabla g_i\|_\Omega^2 \label{eq:weighted-overlap-inequality}
\end{equation}
for all positive weights. If \(\|\nabla g_i\|_\Omega\le y_i\) with
\(y_i>0\), the choice \(\nu_i=y_i\) bounds the right-hand side by
\((\sum_iy_i)^2\).
\end{lem}

\begin{proof}
Fix a coarse element \(K\) and let \(J_K=\{j:K\in W_j\}\). Only the
functions \(g_j\) with \(j\in J_K\) can be nonzero on \(K\).
The weighted Cauchy--Schwarz inequality gives
\begin{equation}
 \Big|\sum_{j\in J_K}\nabla g_j\Big|^2\le\Big(\sum_{j\in J_K}\nu_j\Big)\sum_{j\in J_K}\frac{|\nabla g_j|^2}{\nu_j} \label{eq:weighted-cauchy-schwarz}
\end{equation}
almost everywhere on \(K\). The first factor is at most
\(\Theta_j^W(\nu)\) for every \(j\in J_K\). Integration over \(K\)
and summation over all coarse elements prove
\eqref{eq:weighted-overlap-inequality}. For \(\nu_i=y_i\), every overlap
load is at most \(\sum_jy_j\), and the right-hand side is at most
\(\sum_i\Theta_i^W(y)y_i\le(\sum_iy_i)^2\).
\end{proof}

Indices with \(y_i=0\) correspond to \(g_i=0\) and may be omitted
before the weights are chosen. The weights must be fixed from computed
quantities, not from the unknown norms \(\|\nabla g_i\|_\Omega\);
the right-hand side of \eqref{eq:weighted-overlap-inequality} is
nondecreasing in each of these norms.

\begin{lem}[Norm of the complementary projection]
\label{lem:complementary-projection}
Equip \(V_h\) with the gradient inner product, regard \(\Pi_H\) as an
operator on \(V_h\), and let \(\|\cdot\|_\nabla\) denote the induced
operator norm. If \(V_H\ne\{0\}\) and \(V_f\ne\{0\}\), then
\begin{equation}
 \|I-\Pi_H\|_\nabla=\|\Pi_H\|_\nabla \label{eq:complementary-projection-norm}
\end{equation}
and consequently \(C_J=\max\{1,C_\Pi\}\) is admissible in
\eqref{eq:kernel-global-stability}.
\end{lem}

\begin{proof}
By \eqref{eq:coarse-reproduction}, \(\Pi_H\) is idempotent on \(V_h\),
with range \(V_H\) and kernel \(V_f\), and \(V_h=V_H\oplus V_f\).
The identity holds for every idempotent operator other than zero and
the identity on an inner-product space~\cite{xu-zikatanov2003}; we
include the short argument. Let \(\gamma\) be the largest cosine
\(|(\nabla r,\nabla z)_\Omega|/(\|\nabla r\|_\Omega\|\nabla z\|_\Omega)\)
over nonzero \(r\in V_H\) and \(z\in V_f\). It is less than one, since
the spaces are finite dimensional and intersect trivially. For fixed
nonzero \(r\in V_H\), orthogonal projection onto \(V_f\) gives
\begin{equation}
 \min_{z\in V_f}\|\nabla(r+z)\|_\Omega^2=\big(1-\gamma_r^2\big)\|\nabla r\|_\Omega^2 \label{eq:projection-angle-identity}
\end{equation}
where \(\gamma_r\) is the largest cosine for this \(r\). Since
\(\Pi_H(r+z)=r\), taking the supremum over \(r\) gives
\(\|\Pi_H\|_\nabla^2=(1-\gamma^2)^{-1}\). Exchanging the roles of
\(V_H\) and \(V_f\) gives the same value for \(I-\Pi_H\), the
projection onto \(V_f\) along \(V_H\), because \(\gamma\) is symmetric
in the two spaces. If \(V_f=\{0\}\), then \(I-\Pi_H=0\); if
\(V_H=\{0\}\), then \(I-\Pi_H=I\). This proves the final assertion.
\end{proof}

We apply the weighted overlap inequality
(Lemma~\ref{lem:weighted-overlap}) to the evaluation error of
Section~\ref{subsec:inexact-solves}. Let \(a_i\) be the assembly norms
in \eqref{eq:matrix-assembly-certificate}, let \(\mathsf D_i\) be the
diagonal matrix of the nodal values of \(\phi_i\) on the interior
patch degrees of freedom, and define
\begin{equation}
 \mu_i^2=\lambda_{\max}(\mathsf D_i^T\mathsf K_i\mathsf D_i,\mathsf K_i) \label{eq:weight-multiplier-norm}
\end{equation}
with \(\mu_i=0\) for an empty patch space. Thus \(\mu_i\) is the
gradient norm of multiplication by \(\phi_i\) followed by \(I_h\) on
the patch space. Write \(\Theta_i^0\) for the overlap load of the
supports \(W_i=\omega_i\) and \(\Theta_i^+\) for that of
\(W_i=N(\omega_i)\).

\begin{prop}[Overlap-aware assembly certificates]
\label{prop:overlap-certificates}
For the constraint-preserving approximate corrections in
Section~\ref{subsec:inexact-solves} and all positive weights,
\begin{align}
 \|\nabla(\widehat u-M_h^{\mathcal D}(v_H))\|_\Omega&\le\frac1m\Big(\sum_i\frac{\Theta_i^+(\nu)}{\nu_i}a_i^2\sigma_i^2\Big)^{1/2} \label{eq:overlap-certificate-exact}\\
 \|\nabla(\widehat u-M_h^{\mathcal D}(v_H))\|_\Omega&\le\frac{\|I-\Pi_H\|_\nabla}m\Big(\sum_i\frac{\Theta_i^0(\nu)}{\nu_i}\mu_i^2\sigma_i^2\Big)^{1/2} \label{eq:overlap-certificate-product}
\end{align}
With \(\nu_i=a_i\sigma_i\), the right-hand side of
\eqref{eq:overlap-certificate-exact} is at most \(\zeta_{\rm mat}\).
With unit weights, \(\Theta_i^0=d+1\) in
\eqref{eq:overlap-certificate-product}.
\end{prop}

\begin{proof}
Let \(e_i=\widehat q_i-q_i(v_H)\). By \eqref{eq:inexact-patch-error},
\(\|\nabla e_i\|_{D_i}\le\sigma_i/m\). The function
\(I_h(\phi_ie_i)\) vanishes outside \(\omega_i\), because \(\phi_i\)
vanishes on every coarse element not containing the vertex \(i\).
Its coefficient vector is \(\mathsf D_i\) times that of \(e_i\), and
both functions are supported in \(D_i\); hence its gradient norm is at
most \(\mu_i\sigma_i/m\). By the support property in
Lemma~\ref{lem:local-projection} with \(r=1\),
\(S_ie_i=(I-\Pi_H)I_h(\phi_ie_i)\) vanishes outside \(N(\omega_i)\),
and \(\|\nabla S_ie_i\|_\Omega\le a_i\sigma_i/m\) by the definition of \(a_i\) in \eqref{eq:matrix-assembly-certificate}. Since
\(\widehat u-M_h^{\mathcal D}(v_H)=\sum_iS_ie_i\), Lemma~\ref{lem:weighted-overlap}
with \(W_i=N(\omega_i)\) proves \eqref{eq:overlap-certificate-exact}.
Writing the same difference as \((I-\Pi_H)\sum_iI_h(\phi_ie_i)\) and
applying Lemma~\ref{lem:weighted-overlap} with \(W_i=\omega_i\) proves
\eqref{eq:overlap-certificate-product}. The comparison with
\(\zeta_{\rm mat}\) is the final assertion of
Lemma~\ref{lem:weighted-overlap}. Finally, every coarse element lies in
the supports of exactly its \(d+1\) vertex hats, including boundary
hats.
\end{proof}

The bounds \eqref{eq:overlap-certificate-exact} and
\eqref{eq:overlap-certificate-product} use locality in different ways. The first keeps the
exact assembly norms and enlarges each support by the projection
radius. The second keeps the partition-of-unity supports and applies
the projection once to the whole sum; by
Lemma~\ref{lem:complementary-projection} its prefactor equals
\(\|\Pi_H\|_\nabla\) in nondegenerate cases. With unit weights it is
the computable counterpart of Lemma~\ref{lem:inexact-patch-assembly},
with \(N_{\rm pu}=d+1\). The minimum over both bounds and over any
finite set of computed weight choices is again a certificate, and the
self-weighted choice is never weaker than the summed bound. The norms
\(\mu_i\) require only diagonally scaled patch stiffness matrices. All
bounds retain the global monotonicity constant from
\eqref{eq:inexact-patch-error}.

\subsection{Residual Certificates for Inexact Patch Evaluation}
\label{subsec:numerical-certificates}

We evaluate the certificates of Section~\ref{subsec:computable-certificates}
for the baseline local weak projection with \(n_H=6\), \(n_h=24\)
and depths one and two. Let \(u_{\rm ref}\in V_h\) denote the
computed fine reference on this fine grid. At the fixed comparator coordinate
\(v_H=\Pi_Hu_{\rm ref}\), we vary the patch Newton tolerance and
recompute the state and tangent tests. These six evaluations are a
separate diagnostic, not additional outer optimization runs.
The saturated law of Example~\ref{example:saturated-law} has the valid
global lower Hessian bound \(m=1\)
in this experiment. All quantities in Table~\ref{tab:residual-certificates}
are absolute gradient norms, rather than the relative errors plotted in
Section~\ref{sec:numerical-experiments} of the main paper.

\begin{table}[htbp]
\centering
\small
\begin{tabular}{ccccccc}
\toprule
Depth & Tolerance & \(\zeta_{\rm mat}\) & \(\zeta_{\rm ov}\) & \begin{tabular}{@{}c@{}}Evaluation\\difference\end{tabular} & \(\rho/m\) & \begin{tabular}{@{}c@{}}Reference\\difference\end{tabular}\\
\midrule
1 & \(1\times10^{-3}\) & \(1.49\times10^{-2}\) & \(8.02\times10^{-3}\) & \(2.19\times10^{-4}\) & \(1.58\times10^{-1}\) & \(1.87\times10^{-2}\) \\
1 & \(1\times10^{-7}\) & \(2.70\times10^{-6}\) & \(1.62\times10^{-6}\) & \(6.08\times10^{-8}\) & \(1.58\times10^{-1}\) & \(1.87\times10^{-2}\) \\
1 & \(2\times10^{-11}\) & \(1.23\times10^{-10}\) & \(1.22\times10^{-10}\) & \(1.82\times10^{-12}\) & \(1.58\times10^{-1}\) & \(1.87\times10^{-2}\) \\
2 & \(1\times10^{-3}\) & \(1.87\times10^{-2}\) & \(9.44\times10^{-3}\) & \(1.41\times10^{-4}\) & \(2.26\times10^{-2}\) & \(2.61\times10^{-3}\) \\
2 & \(1\times10^{-7}\) & \(6.59\times10^{-6}\) & \(3.63\times10^{-6}\) & \(3.07\times10^{-8}\) & \(2.27\times10^{-2}\) & \(2.62\times10^{-3}\) \\
2 & \(2\times10^{-11}\) & \(4.98\times10^{-11}\) & \(4.80\times10^{-11}\) & \(1.17\times10^{-12}\) & \(2.27\times10^{-2}\) & \(2.62\times10^{-3}\) \\

\bottomrule
\end{tabular}
\caption{Computed residual certificates. The overlap-aware bound
\(\zeta_{\rm ov}\) is the smallest of \eqref{eq:overlap-certificate-exact}
and \eqref{eq:overlap-certificate-product} with unit weights and with
weights equal to the computed amplitudes. The evaluation difference compares
with a reconstruction at patch tolerance \(2\times10^{-13}\); its own
assembly-error bound is below \(2.0\times10^{-12}\). The reference
difference compares with the computed fine solution, whose residual
bounds its gradient error by \(9.1\times10^{-13}\). Thus each measured
difference is bounded by the displayed certificate plus the corresponding
comparison-state bound, up to floating-point error.}
\label{tab:residual-certificates}
\end{table}

Tightening the patch tolerance reduces the patch evaluation bound,
while the full-residual bound changes little. The latter still includes
the localization error. At depth two and the
production patch tolerance \(2\times10^{-11}\), the assembly-error
bound is below \(5.0\times10^{-11}\), compared with a measured
fine-reference difference of approximately \(2.62\times10^{-3}\).
This separates patch evaluation error from localization for the tested
coordinate. The full bound is conservative, as expected from using the
global monotonicity constant.

The overlap-aware bound \(\zeta_{\rm ov}\) of
Proposition~\ref{prop:overlap-certificates} is smaller than
\(\zeta_{\rm mat}\) by factors between \(1.66\) and \(1.99\) at the
two looser tolerances, where the patch residuals are spread over many
patches. At the production tolerance, a single patch carries about
\(50\%\) and \(96\%\) of \(\sum_ia_i\sigma_i\) at depths one and two.
The self-weighted bound then nearly coincides with \(\zeta_{\rm mat}\),
which it cannot exceed. The computed projection norms satisfy
\(\|I-\Pi_H\|_\nabla=\|\Pi_H\|_\nabla=1.7013\) to twelve digits, as
Lemma~\ref{lem:complementary-projection} predicts. As a diagnostic,
the patch corrections were also compared with those at tolerance
\(2\times10^{-13}\). In the Euclidean sum over patches, \(\sigma_i/m\)
exceeds these patch gradient errors by factors between \(5.8\) and
\(8.9\). Most of the remaining overestimation therefore comes from the
global monotonicity constant in \eqref{eq:inexact-patch-error}, not
from assembly.

Constrained Riesz norms were cross-checked with independent nullspace
representations. Kernel feasibility, test-map section properties and
linear-system residuals were also recorded. The diagnostic retains the
coarse dual residuals and induced test-map norms; it does not evaluate
the analytical localization constants or certify the oversampling
condition. These are floating-point checks, with the arithmetic
limitations stated in Section~\ref{subsec:computable-certificates}.

\paragraph{Certificates at Saved Outer-Solver States.}
\label{par:saved-state-certificates}
We next apply the certificates to all six saved local weak outer-solver
states on the fine grid with \(n_h=48\), using depths two and four and the
three initial guesses of Section~\ref{subsec:numerical-computed-states}.
We retain the original states, coordinates and optimizer flags and
evaluate the certificates without repeating the outer optimization. Patch corrections and test maps are
reevaluated at each saved coordinate. The full residual is evaluated
at the saved state itself.

For these larger patches the summed bound uses a product bound for the
assembly norms. With \(\mu_i\) from \eqref{eq:weight-multiplier-norm},
projection stability and the triangle inequality give
\begin{equation}
 a_i\le(1+\|\Pi_H\|_\nabla)\mu_i \label{eq:product-assembly-norm}
\end{equation}
where \(\|\Pi_H\|_\nabla\) is computed from the fine and coarse
stiffness matrices. This avoids forming dense fine-grid assembly
matrices. By Lemma~\ref{lem:complementary-projection}, the factor
\(1+\|\Pi_H\|_\nabla\) could be replaced by \(\|\Pi_H\|_\nabla\); the
summed column of Table~\ref{tab:saved-state-certificates} retains the
original factor. The overlap column uses the smallest of the bounds in
Proposition~\ref{prop:overlap-certificates}, with the exact norms
\(a_i\) evaluated by a matrix-free Lanczos iteration, since \(\Pi_H\)
has rank equal to the coarse dimension. Each resulting assembly bound
applies to the reevaluated reconstruction. Adding the measured gradient difference
between saved and reevaluated states bounds the evaluation error of
the saved state relative to the exact patch map at that coordinate.
This added difference accounts for the different patch-solver histories.

\begin{table}[htbp]
\centering
\small
\begin{tabular}{clcccc}
\toprule
Depth & Initial guess & \begin{tabular}{@{}c@{}}Summed\\bound\end{tabular} & \begin{tabular}{@{}c@{}}Overlap\\bound\end{tabular} & \(\rho/m\) & \begin{tabular}{@{}c@{}}Reference\\difference\end{tabular}\\
\midrule
2 & Zero & \(4.51\times10^{-10}\) & \(1.39\times10^{-10}\) & \(3.02\times10^{-2}\) & \(3.75\times10^{-3}\) \\
2 & Random & \(1.06\times10^{-13}\) & \(1.95\times10^{-14}\) & \(3.02\times10^{-2}\) & \(3.75\times10^{-3}\) \\
2 & Reference-informed & \(2.02\times10^{-12}\) & \(1.82\times10^{-12}\) & \(3.02\times10^{-2}\) & \(3.75\times10^{-3}\) \\
4 & Zero & \(8.76\times10^{-10}\) & \(1.98\times10^{-10}\) & \(1.02\times10^{-3}\) & \(1.52\times10^{-4}\) \\
4 & Random & \(1.31\times10^{-13}\) & \(2.33\times10^{-14}\) & \(1.02\times10^{-3}\) & \(1.52\times10^{-4}\) \\
4 & Reference-informed & \(6.69\times10^{-12}\) & \(6.42\times10^{-12}\) & \(1.02\times10^{-3}\) & \(1.52\times10^{-4}\) \\

\bottomrule
\end{tabular}
\caption{Absolute gradient certificates for the six saved computed states.
Both evaluation bounds include the saved-to-reevaluated state difference.
The overlap column was computed in a separate run; the summed bounds
recomputed there agree with the displayed values, except for the two
random-start entries at rounding level, which differ by about \(3\%\).
The fine reference has a gradient residual bound below
\(6.4\times10^{-12}\). All total-error checks pass after accounting
for that reference error.}
\label{tab:saved-state-certificates}
\end{table}

Table~\ref{tab:saved-state-certificates} summarizes the saved-state checks.
The total-error bound exceeds the measured reference difference by
factors between \(6.70\) and \(8.06\). Every summed evaluation-error
bound is below \(8.8\times10^{-10}\), and every overlap-aware bound is
below \(2.0\times10^{-10}\); both are well below the measured discrete
errors. The overlap-aware assembly bounds are smaller than the summed
product bounds by factors between \(3.4\) and \(6.1\). For the
reference-informed starts, the saved-to-reevaluated difference
dominates both totals. The induced test-map norms are approximately \(1.08\), and
the section property and the constrained Riesz identities hold to
floating-point accuracy. The coarse test residuals are evaluated at the saved states
using the reevaluated test maps. Their small values complement the
full-residual bounds; they do not certify the analytical oversampling
condition.


\section{Energies with \texorpdfstring{{\boldmath$p$}}{p}-Structure}
\label{sec:p-structure}

Sections~\ref{sec:energy-setting}--\ref{sec:stationary-residual} assume
the uniform Hessian bounds \eqref{eq:hessian-bounds}. They exclude laws
with unbounded growth, such as the quartic law. This section replaces
\eqref{eq:hessian-bounds} by a \(p\)-structure condition with \(p\ge2\)
and proves the corresponding energy and computed-state estimates. The secant identity of Lemma~\ref{lem:secant-harmonicity}, the cutoff
construction of Appendix~\ref{sec:buffered-cutoff} and the assembly
identity \eqref{eq:exact-assembly-error} are unchanged. The constant upper bound \(M\) is replaced by a weight that
depends on the ideal state and on the localization error. The
contraction factor then depends on scale-invariant local quantities of
that weight. For \(p=2\) and a constant coefficient, these quantities
are mesh constants, and the results reduce to those of
Sections~\ref{sec:localization}--\ref{sec:stationary-residual}.
Section~\ref{sec:p-constants} discusses the size of the weighted
quantities and their dependence on the load.

\paragraph{Summary of Results.}
\label{par:p-summary}
Theorem~\ref{thm:p-energy} extends \eqref{eq:intro-energy-rate} to
potentials with \(p\)-structure \eqref{eq:p-structure}, with a
piecewise constant coefficient \(a_\varepsilon\) and structural constants
\(\hat m,\hat M\). Its uniform-depth form \eqref{eq:p-rate} reads
\begin{equation}
 0\le\mathcal E(u_{\rm ms}^\ell)-\mathcal E(u_h)\le C\,\Gamma_\Omega(u_h)\Gamma_\Pi(u_h)(1+\ell)^de^{-b\ell}\frac{\hat M\hat L_F^2}{\hat m^2} \quad \forall\ell\ge1 \label{eq:intro-p-energy-rate}
\end{equation}
where \(\hat L_F\) is the load norm weighted by \(a_\varepsilon\). The
constant \(C\) depends on \(p\), \(\hat M/\hat m\), \(d\) and the mesh
family, and \(b\) is determined by \(p\), \(\hat M/\hat m\), a cutoff
constant and \(\Gamma_\Omega(u_h)\). The quantity \(\Gamma_\Omega\)
collects weighted Poincar\'e constants and Muckenhoupt-type ratios of
\(a_\varepsilon(1+|\nabla u_h|)^{p-2}\) at the coarse scale, and
\(\Gamma_\Pi\) is a weighted stability constant of the projection. Both
are invariant under scaling of the weight, so they measure its variation
rather than its size. For \(p=2\) and a constant coefficient they are
mesh constants, and \eqref{eq:intro-p-energy-rate} reduces to
\eqref{eq:intro-energy-rate}. For \(d=2\), and for \(d=3\) with
\(p\le6\), they are bounded independently of the fine mesh, with the crude
bound of Lemma~\ref{lem:p-crude}, which grows with the load and with a
negative power of \(H\). This
bound therefore does not by itself give logarithmic oversampling.
Section~\ref{sec:p-constants} shows that
the load dependence of the Muckenhoupt-type ratios reduces to the shape of
the gradient magnitude, and the experiments show that the quantities
saturate as the load increases. Theorem~\ref{thm:p-a2-uniform} gives
bounds independent of the load and of the fine mesh if
\(|\nabla u_h|^{p-2}\) satisfies a local Muckenhoupt \(A_2\) condition with a
load-independent constant. Theorem~\ref{thm:p-fixed-mesh} gives bounds
independent of the load on every fixed mesh for asymptotically
homogeneous laws, such as the quartic law, on meshes where the large-load
limit is nondegenerate. Theorem~\ref{thm:p-computed} extends the
computed-state estimate of Theorem~\ref{thm:computed-state-error} under compatibility and oversampling conditions at
the computed coordinate. By Corollary~\ref{cor:p-apriori}, for \(d=2\) or
\(p\le6\), both hold for all coordinates with bounded energy once the depth
exceeds a threshold. For fixed \(H\) this threshold grows like
\(\log(1/h)\), but it also grows with the load and like a power of
\(1/H\).

\paragraph{Weighted Localization and \texorpdfstring{{\boldmath$p$}}{p}-Structure.}
\label{par:literature-weighted}
For linear problems, the localization error of orthogonal decompositions
generally deteriorates with the coefficient contrast. Peterseim and
Scheichl~\cite{peterseim-scheichl2016} use coefficient-dependent
quasi-interpolation and obtain a corrector decay rate that depends on
constants of the interpolation operator, controlled by a weighted
Poincar\'e inequality, and not on the contrast. Hellman and
M\aa lqvist~\cite{hellman-malqvist2017} construct an interpolation
operator for two-valued coefficients and prove a decay rate independent
of the contrast under geometric assumptions on inclusions and channels.
The weighted quantities of this section play a
comparable role. Here, however, the projection is fixed and independent
of the coefficient, and the weight depends on the solution. For problems
with \(p\)-structure, the equivalence of the natural distance with a
weighted quadratic form, and with differences of a nonlinear function
\(V\) of the gradient, is a standard tool. It underlies the quasi-norm
error analysis of Barrett and Liu~\cite{barrett-liu1994} and the shifted
\(N\)-functions of Diening and Ettwein~\cite{diening-ettwein2008}.
Lemmas~\ref{lem:p-secant} and~\ref{lem:p-v-function} state these
equivalences with explicit constants.

\subsection{Assumption and Structural Lemmas}
\label{subsec:p-assumption}

Let \(p\ge2\). The potential satisfies the measurability, smoothness and
normalization conditions of Section~\ref{subsec:potential-bounds}, with
\eqref{eq:hessian-bounds} replaced as follows. There is a coefficient
\(a_\varepsilon\), piecewise constant on \(\mathcal T_h\), with
\(0<a_-\le a_\varepsilon\le a_+\), and there are structural constants
\(0<\hat m\le\hat M<\infty\) such that
\begin{equation}
 \hat m\,a_\varepsilon(x)(1+|\xi|)^{p-2}|\zeta|^2\le\zeta^TD_\xi a(x,\xi)\zeta\le\hat M\,a_\varepsilon(x)(1+|\xi|)^{p-2}|\zeta|^2 \quad \forall\xi,\zeta\in\mathbb R^d \label{eq:p-structure}
\end{equation}
for almost every \(x\). The bounds \(a_\pm\) are used only for
well-posedness and in the crude bounds of
Section~\ref{subsec:p-crude}. For \(p=2\),
\eqref{eq:p-structure} implies \eqref{eq:hessian-bounds} with
\(m=\hat ma_-\) and \(M=\hat Ma_+\). Conversely,
\eqref{eq:hessian-bounds} is \eqref{eq:p-structure} with \(p=2\),
\(a_\varepsilon=1\), \(\hat m=m\) and \(\hat M=M\).

\begin{example}[The quartic law]
\label{example:quartic-law}
The law \(a(x,\xi)=a_\varepsilon(x)(1+|\xi|^2)\xi\) has the potential
\(\psi=a_\varepsilon(|\xi|^2/2+|\xi|^4/4)\) and the Jacobian
\begin{equation}
 D_\xi a(x,\xi)=a_\varepsilon(x)\big[(1+|\xi|^2)I+2\xi\otimes\xi\big] \label{eq:quartic-jacobian}
\end{equation}
with eigenvalues \(a_\varepsilon(1+|\xi|^2)\) and \(a_\varepsilon(1+3|\xi|^2)\).
Since \((1+|\xi|)^2/2\le1+|\xi|^2\) and \(1+3|\xi|^2\le3(1+|\xi|)^2\),
condition \eqref{eq:p-structure} holds with \(p=4\), \(\hat m=1/2\) and
\(\hat M=3\), independently of the coefficient contrast. The saturated
law of Example~\ref{example:saturated-law} satisfies
\eqref{eq:p-structure} with \(p=2\), \(\hat m=1\) and
\(\hat M=1+9\alpha/(8\beta)\).
\end{example}

For vectors \(\xi,\eta\), states \(U\in V_h\) and loads \(F\in V_h'\), we write
\begin{align}
 \omega(\xi,\eta)&=(1+|\xi|+|\eta|)^{p-2},\qquad \widehat W_U=a_\varepsilon(1+|\nabla U|)^{p-2} \label{eq:p-weights}\\
 c_p&=\frac{2^{2-p}}{p-1},\qquad \hat\chi_p=\frac{\hat M}{c_p\hat m},\qquad \hat L_F=\sup_{0\ne w\in V_h}\frac{|F(w)|}{\big(\int_\Omega a_\varepsilon|\nabla w|^2\,dx\big)^{1/2}} \label{eq:p-constants}
\end{align}
with the convention \(|\zeta|^{p-2}=1\) for \(p=2\). The secant weight
\(\omega(\xi,\eta)\) is always written with its arguments, to distinguish
it from the patches \(\omega_i\) and \(\omega\). For piecewise
affine states, \(\widehat W_U\) and \(a_\varepsilon\omega(\nabla U,\nabla U_D)\) are
constant on every fine element. For \(p=2\), \(c_2=1\), so
\(\hat\chi_2=\hat M/\hat m\).

\begin{lem}[Secant equivalence]
\label{lem:p-secant}
For \(\xi,\eta\in\mathbb R^d\) and \(\xi_t=\xi+t(\eta-\xi)\),
\begin{equation}
 c_p\,\omega(\xi,\eta)\le\int_0^1(1+|\xi_t|)^{p-2}\,dt\le\omega(\xi,\eta) \label{eq:p-scalar-secant}
\end{equation}
Consequently, the averaged tensor
\(B=\int_0^1D_\xi a(x,\xi_t)\,dt\) satisfies
\begin{equation}
 c_p\hat m\,a_\varepsilon\,\omega(\xi,\eta)|\zeta|^2\le\zeta^TB\zeta\le\hat M\,a_\varepsilon\,\omega(\xi,\eta)|\zeta|^2 \quad \forall\zeta\in\mathbb R^d \label{eq:p-secant-tensor}
\end{equation}
and, for \(v,w\in V_h\) and every element region \(E\),
\begin{equation}
 c_p\hat m\int_Ea_\varepsilon\,\omega(\nabla v,\nabla w)|\nabla(v-w)|^2\,dx\le\mathscr D_E(v,w)\le\hat M\int_Ea_\varepsilon\,\omega(\nabla v,\nabla w)|\nabla(v-w)|^2\,dx \label{eq:p-natural-equivalence}
\end{equation}
\end{lem}

\begin{proof}
The upper bound in \eqref{eq:p-scalar-secant} follows from
\(|\xi_t|\le\max\{|\xi|,|\eta|\}\). For the lower bound, write
\(s=p-2\), \(\alpha=|\xi|\), \(\beta=|\eta|\) and assume
\(\alpha+\beta>0\); otherwise both sides equal one. The reverse triangle
inequality gives \(|\xi_t|\ge|(1-t)\alpha-t\beta|\), and the right side
is affine in \(t\) with slope \(-(\alpha+\beta)\). Substitution gives
\begin{equation}
 \int_0^1(1+|\xi_t|)^s\,dt\ge\frac{(1+\alpha)^{s+1}+(1+\beta)^{s+1}-2}{(s+1)(\alpha+\beta)}\ge\frac{(1+\gamma)^{s+1}-1}{(s+1)\gamma} \label{eq:p-secant-profile}
\end{equation}
with \(\gamma=(\alpha+\beta)/2\), where the second step uses convexity
of \(x\mapsto(1+x)^{s+1}\). Since
\((1+\gamma)^{s+1}-1\ge\gamma(1+\gamma)^s\) and
\((1+\gamma)^s\ge2^{-s}(1+2\gamma)^s\), the right side is at least
\(2^{-s}(1+\alpha+\beta)^s/(s+1)=c_p\,\omega(\xi,\eta)\). The tensor
bounds follow from \eqref{eq:p-structure}. The last statement follows
from \(\mathscr D_E(v,w)=\int_E\nabla(v-w)^TB\nabla(v-w)\,dx\), with
\(B\) evaluated at \(\xi=\nabla w\) and \(\eta=\nabla v\).
\end{proof}

Equivalences of this type are standard for problems with
\(p\)-structure~\cite{barrett-liu1994,diening-ettwein2008}; the proof is
included for the explicit constants.
The constant \(c_p\) is sharp; it is approached for \(\eta=-\xi\) with
\(|\xi|\to\infty\). For \(p=2\), \eqref{eq:p-natural-equivalence} is
\eqref{eq:natural-equivalence} with the coefficient \(a_\varepsilon\).

\begin{lem}[Weight splitting]
\label{lem:p-splitting}
For \(U,e\in V_h\) and \(U_D=U+e\),
\begin{equation}
 \max\big\{\widehat W_U,a_\varepsilon|\nabla e|^{p-2}\big\}\le a_\varepsilon\,\omega(\nabla U,\nabla U_D)\le3^{p-2}\big(\widehat W_U+a_\varepsilon|\nabla e|^{p-2}\big) \label{eq:p-splitting}
\end{equation}
pointwise. For \(p=2\), \(a_\varepsilon\omega(\nabla U,\nabla U_D)=\widehat W_U=a_\varepsilon\).
\end{lem}

\begin{proof}
The quantity \(1+|\nabla U|+|\nabla U_D|\) is at least \(1+|\nabla U|\)
and at least \(|\nabla e|\). It is at most
\(2(1+|\nabla U|)+|\nabla e|\le3\max\{1+|\nabla U|,|\nabla e|\}\).
\end{proof}

\begin{lem}[Coercivity and data bounds]
\label{lem:p-coercivity}
Let \(c_\psi=2^{-3-p}\), \(c_B=4^{-1-p}\) for \(p>2\), and \(c_B=1/2\)
for \(p=2\). Then
\begin{align}
 \psi(x,\xi)&\ge c_\psi\hat m\,a_\varepsilon\big(|\xi|^2+|\xi|^p\big) \quad \forall\xi\in\mathbb R^d \label{eq:p-potential-lower}\\
 \mathcal B(v,w)&\ge c_B\hat m\int_\Omega a_\varepsilon\,\omega(\nabla v,\nabla w)|\nabla(v-w)|^2\,dx \quad \forall v,w\in V_h \label{eq:p-bregman-lower}
\end{align}
The energy has a unique minimizer \(u_h\in V_h\), and
\begin{align}
 \int_\Omega a_\varepsilon|\nabla u_h|^2\,dx&\le\frac{\hat L_F^2}{\hat m^2} \label{eq:p-data-two}\\
 \int_\Omega a_\varepsilon\big(|\nabla u_h|^2+|\nabla u_h|^p\big)\,dx&\le\frac{\hat L_F^2}{c_\psi\hat m^2} \label{eq:p-data-p}\\
 \int_\Omega a_\varepsilon|\nabla u_h|^p\,dx&\le\Big(\frac{|\Omega|_a^{(p-2)/(2p)}\hat L_F}{c_\psi\hat m}\Big)^{p/(p-1)} \label{eq:p-data-large}
\end{align}
where \(|\Omega|_a=\int_\Omega a_\varepsilon\,dx\).
\end{lem}

\begin{proof}
Taylor's formula gives
\(\psi(x,\xi)=\int_0^1(1-t)\xi^TD_\xi a(x,t\xi)\xi\,dt\). Restricting the
integral to \(t\in[1/2,3/4]\) bounds it from below by
\(2^{-4}\hat ma_\varepsilon(1+|\xi|/2)^{p-2}|\xi|^2\ge2^{-2-p}\hat ma_\varepsilon(1+|\xi|)^{p-2}|\xi|^2\).
Since \((1+|\xi|)^{p-2}|\xi|^2\ge\max\{|\xi|^2,|\xi|^p\}\), this proves
\eqref{eq:p-potential-lower}. For \eqref{eq:p-bregman-lower}, the
Bregman integrand is
\(\int_0^1(1-t)\zeta^TD_\xi a(x,\xi_t)\zeta\,dt\) with \(\xi=\nabla w\),
\(\eta=\nabla v\) and \(\zeta=\eta-\xi\). For \(p=2\) it is at least
\(\hat ma_\varepsilon|\zeta|^2/2\). For \(p>2\), if \(|\xi|\ge|\eta|\),
then \(|\xi_t|\ge(|\xi|+|\eta|)/4\) and \(1-t\ge3/4\) for
\(t\in[0,1/4]\). If \(|\eta|>|\xi|\), the same lower bound for
\(|\xi_t|\) and \(1-t\ge1/8\) hold for \(t\in[3/4,7/8]\). In both cases
the integrand is at least \(2^{-6}4^{2-p}\hat ma_\varepsilon\omega(\xi,\eta)|\zeta|^2\).

The lower bound in \eqref{eq:p-structure} gives
\(a(x,\xi)\cdot\xi\ge\hat ma_\varepsilon|\xi|^2\) and strict convexity of
\(\psi(x,\cdot)\). With \eqref{eq:p-potential-lower},
\(\mathcal E(v)\ge c_\psi\hat m\int_\Omega a_\varepsilon|\nabla v|^2-\hat L_F(\int_\Omega a_\varepsilon|\nabla v|^2)^{1/2}\),
so \(\mathcal E\) is coercive and has a unique minimizer. Testing its
equation with \(u_h\) proves \eqref{eq:p-data-two}. Since
\(\mathcal E(u_h)\le\mathcal E(0)=0\),
\(c_\psi\hat m\int a_\varepsilon(|\nabla u_h|^2+|\nabla u_h|^p)\le F(u_h)\le\hat L_F(\int a_\varepsilon|\nabla u_h|^2)^{1/2}\),
which gives \eqref{eq:p-data-p} by \eqref{eq:p-data-two}. Bounding the
last factor instead by H\"older's inequality for the measure
\(a_\varepsilon\,dx\),
\((\int a_\varepsilon|\nabla u_h|^2)^{1/2}\le|\Omega|_a^{(p-2)/(2p)}(\int a_\varepsilon|\nabla u_h|^p)^{1/p}\),
proves \eqref{eq:p-data-large}.
\end{proof}

The bound \eqref{eq:p-data-large} grows like \(\hat L_F^{p/(p-1)}\)
rather than \(\hat L_F^2\), and it is the relevant one for large loads.
The global and patch corrections of Section~\ref{sec:reconstruction}
remain well defined, since the restricted energies are coercive and
strictly convex. The patch maps are continuously differentiable: in
fixed bases, continuity of \(D_\xi a\) and the local bound
\(|D_\xi a(x,\xi)|\le\hat Ma_+(1+|\xi|)^{p-2}\) justify differentiation
under the integral, and the lower bound in \eqref{eq:p-structure} makes
the patch Jacobians positive definite. The first part of
Theorem~\ref{thm:energy-comparison}, namely attainment,
\eqref{eq:bregman-best-approximation} and \eqref{eq:energy-comparator},
uses only coercivity, continuity and the section property
\eqref{eq:general-section}, and
therefore holds. Lemma~\ref{lem:secant-harmonicity} and the orthogonal
identity \eqref{eq:outer-pythagoras} use only symmetry and positivity of
the secant tensor, and also hold. The Lipschitz bound
\eqref{eq:patch-state-lipschitz} and the gradient statements of
Theorem~\ref{thm:energy-comparison} rely on \(M\) and are not used below.

\subsection{Weighted Local Quantities}
\label{subsec:p-quantities}

In Sections~\ref{sec:localization} and~\ref{sec:assembly-error}, the cutoff
and assembly estimates use the unweighted kernel Poincar\'e inequality
\eqref{eq:kernel-poincare}, and the bounds \(m\) and \(M\) convert gradient
norms into natural distances. Under \eqref{eq:p-structure}, the natural
distance is a weighted gradient norm, and the corresponding inequalities
carry the weight
\(\widehat W_U\) of the ideal state and, for the part of the weight that
comes from the localization error, the coefficient \(a_\varepsilon\) in
\(L^p\). Let \(\operatorname{star}z\) denote the support of the coarse
hat \(\lambda_z\).

\begin{definition}[Weighted local quantities]
\label{def:p-quantities}
For \(U\in V_h\), a coarse element \(T\) and an interior coarse vertex
\(z\), let \(\mathcal P_T(U)\) and \(\mathcal P_T^e\) be the smallest
constants such that
\begin{align}
 H^{-2}\int_T\widehat W_Uv_f^2\,dx&\le\mathcal P_T(U)\int_{N(T)}\widehat W_U|\nabla v_f|^2\,dx \quad \forall v_f\in V_f \label{eq:p-weighted-poincare}\\
 H^{-p}\int_Ta_\varepsilon|v_f|^p\,dx&\le\mathcal P_T^e\int_{N(T)}a_\varepsilon|\nabla v_f|^p\,dx \quad \forall v_f\in V_f \label{eq:p-error-poincare}
\end{align}
and define the ratios
\begin{align}
 \mathcal R_z(U)&=\Big(\frac1{|K_z|}\int_{K_z}\widehat W_U^{-1}\,dx\Big)\Big(\frac1{|\operatorname{star}z|}\int_{\operatorname{star}z}\widehat W_U\,dx\Big) \label{eq:p-ratio}\\
 \mathcal R_z^e&=\Big(\frac1{|K_z|}\int_{K_z}a_\varepsilon^{-1/(p-1)}\,dx\Big)^{p-1}\Big(\frac1{|\operatorname{star}z|}\int_{\operatorname{star}z}a_\varepsilon\,dx\Big) \label{eq:p-error-ratio}
\end{align}
For a set \(A\) of coarse elements, \(\mathcal P_A\), \(\mathcal P_A^e\),
\(\mathcal R_A\) and \(\mathcal R_A^e\) denote the maxima over \(T\in A\)
and over the vertices \(z\) with \(K_z\in A\), with the value zero if
there is no such vertex. Set
\begin{equation}
 \Gamma_A(U)=1+\mathcal P_A(U)\big(1+\mathcal R_A(U)\big)+\Gamma_A^e,\qquad \Gamma_A^e=\mathcal P_A^e\big(1+\mathcal R_A^e\big) \label{eq:p-gamma}
\end{equation}
for \(p>2\), and \(\Gamma_A^e=0\) for \(p=2\). Finally, let
\(\Gamma_\Pi(U)\ge1\) be the smallest constant such that
\begin{align}
 \int_\Omega\widehat W_U|\nabla\Pi_Hv|^2\,dx&\le\Gamma_\Pi(U)\int_\Omega\widehat W_U|\nabla v|^2\,dx \quad \forall v\in V_h \label{eq:p-projection-weighted}\\
 \int_\Omega a_\varepsilon|\nabla\Pi_Hv|^k\,dx&\le\Gamma_\Pi(U)\int_\Omega a_\varepsilon|\nabla v|^k\,dx \quad \forall v\in V_h,\quad k\in\{2,p\} \label{eq:p-projection-coefficient}
\end{align}
\end{definition}

All these quantities are finite. On each discrete configuration the
weights are bounded above and below by positive constants, and the
unweighted inequalities \eqref{eq:kernel-poincare},
\eqref{eq:projection-global-stability} and their \(L^p\) versions in
Lemma~\ref{lem:p-lq} hold. They are invariant under
\(\widehat W_U\mapsto c\widehat W_U\) and \(a_\varepsilon\mapsto ca_\varepsilon\) for
constants \(c>0\). They measure the variation of the weights at the
coarse scale, not their size. The quantities with superscript \(e\) and
the coefficient part of \(\Gamma_\Pi\) depend on \(a_\varepsilon\) only and
not on the load. The ratio satisfies
\(\mathcal R_z(U)\ge|K_z|/|\operatorname{star}z|\) by Jensen's inequality
on \(K_z\), and equals one for constant weights. For \(p=2\) and
\(a_\varepsilon=1\), \eqref{eq:kernel-poincare} gives
\(\mathcal P_T\le C_P^2\), and \(\mathcal R_z=1\) and
\(\Gamma_\Pi\le C_\Pi^2\), so \(\Gamma_A\le1+2C_P^2\) is a mesh constant.
The quantities are defined on the discrete kernel and may depend on the
fine mesh. Section~\ref{subsec:p-crude} gives bounds independent of
\(h\) for \(d=2\), and for \(d=3\) with \(p\le6\).

\subsection{Weighted Cutoff and Contraction}
\label{subsec:p-contraction}

We use the three-layer construction of Appendix~\ref{sec:buffered-cutoff}.
For a seed set \(S\) of coarse elements, \(A=S_3\setminus S\), the
coarse cut-off function \(\eta_H\) and its transition set \(\mathcal T_\eta\) are
defined in \eqref{eq:geometric-cutoff}--\eqref{eq:transition-buffer-inclusion},
and the corrected cutoff of \(e\in V_f\) is \(z_f=(I-\Pi_H)v\) with
\(v=I_h(\eta_He)\) as in \eqref{eq:corrected-cutoff-definition}. The
support and matching statements
\eqref{eq:buffer-inner-matching}--\eqref{eq:buffer-outer-support} are
algebraic and hold without change.

\begin{lem}[Weighted three-layer cutoff]
\label{lem:p-cutoff}
Let \(U\in V_h\), \(e\in V_f\), \(U_D=U+e\) and
\(w=a_\varepsilon\,\omega(\nabla U,\nabla U_D)\). Then
\begin{equation}
 \int_Aw|\nabla z_f|^2\,dx\le C_{{\rm cut},p}\Gamma_A(U)\int_Aw|\nabla e|^2\,dx \label{eq:p-cutoff}
\end{equation}
where \(C_{{\rm cut},p}\) depends only on \(p\), \(d\) and the mesh
family.
\end{lem}

\begin{proof}
All constants \(C\) depend only on \(p\), \(d\) and the mesh family. We
write \(\mathcal P=\mathcal P_A(U)\), \(\mathcal R=\mathcal R_A(U)\),
\(\mathcal P^e=\mathcal P_A^e\) and \(\mathcal R^e=\mathcal R_A^e\), and use
\(N(\mathcal T_\eta)\subset A\) and
\(\operatorname{star}z\subset N(K_z)\subset A\) for \(K_z\in\mathcal T_\eta\).
For \(p=2\) the error parts below are absent, since \(w=\widehat W_U\).

\textbf{Uncorrected product.} On the elements of
\(A\setminus\mathcal T_\eta\), \(v\) equals \(e\) or zero. On a fine
element \(t\) of a transition element \(T\), the weight \(w\) is
constant, and \eqref{eq:product-gradient} with
\(\|\nabla\eta_H\|_{L^\infty}\le CH^{-1}\) gives
\(\int_tw|\nabla v|^2\le C\int_tw(|\nabla e|^2+H^{-2}e^2)\). By
\eqref{eq:p-splitting}, it remains to bound the two parts of
\(H^{-2}\int_Twe^2\). The first part satisfies
\(H^{-2}\int_T\widehat W_Ue^2\le\mathcal P\int_{N(T)}\widehat W_U|\nabla e|^2\).
For the second part, H\"older's inequality for the measure
\(a_\varepsilon\,dx\) and \eqref{eq:p-error-poincare} give
\begin{align}
 H^{-2}\int_Ta_\varepsilon|\nabla e|^{p-2}e^2\,dx&\le\Big(\int_Ta_\varepsilon|\nabla e|^p\,dx\Big)^{(p-2)/p}\Big(H^{-p}\int_Ta_\varepsilon|e|^p\,dx\Big)^{2/p} \label{eq:p-error-product}\\
 &\le(\mathcal P^e)^{2/p}\int_{N(T)}a_\varepsilon|\nabla e|^p\,dx \label{eq:p-error-product-bound}
\end{align}
By the lower bound in \eqref{eq:p-splitting}, both right-hand sides are
at most multiples of \(\int_{N(T)}w|\nabla e|^2\). Summation with bounded
overlap gives
\(\int_Aw|\nabla v|^2\le C(1+\mathcal P+(\mathcal P^e)^{2/p})\int_Aw|\nabla e|^2\).

\textbf{Projection correction.} By
\eqref{eq:cutoff-moment-cancellation},
\(\Pi_Hv=\sum_zc_z\lambda_z\) with \(c_z=\ell_z(g_z)\),
\(g_z=I_h((\eta_H-\eta_H(z))e)\), and only vertices with
\(K=K_z\in\mathcal T_\eta\) contribute. Since \(|\eta_H-\eta_H(z)|\le1\),
nodal norm equivalence on fine elements gives
\(\int_K|g_z|\le C\int_K|e|\), and \(|\theta_z^K|\le CH^{-d}\). The
Cauchy--Schwarz and H\"older inequalities therefore give the two bounds
\begin{align}
 |c_z|^2&\le CH^{-2d}\Big(\int_K\widehat W_Ue^2\,dx\Big)\Big(\int_K\widehat W_U^{-1}\,dx\Big) \label{eq:p-coefficient-background}\\
 |c_z|^2&\le CH^{-2d}\Big(\int_Ka_\varepsilon|e|^p\,dx\Big)^{2/p}\Big(\int_Ka_\varepsilon^{-1/(p-1)}\,dx\Big)^{2(p-1)/p} \label{eq:p-coefficient-error}
\end{align}
By \eqref{eq:p-splitting} and \(|\nabla\lambda_z|\le CH^{-1}\),
\(\int w|\nabla\lambda_z|^2\le CH^{-2}\int_{\operatorname{star}z}(\widehat W_U+a_\varepsilon|\nabla e|^{p-2})\).
For the background part, \eqref{eq:p-coefficient-background},
\(|K||\operatorname{star}z|\le CH^{2d}\) and \eqref{eq:p-weighted-poincare}
give
\begin{equation}
 |c_z|^2H^{-2}\int_{\operatorname{star}z}\widehat W_U\,dx\le C\mathcal R_z(U)H^{-2}\int_K\widehat W_Ue^2\,dx\le C\mathcal R\mathcal P\int_{N(K)}\widehat W_U|\nabla e|^2\,dx \label{eq:p-background-correction}
\end{equation}
For the error part, H\"older's inequality on the star gives
\begin{equation}
 \int_{\operatorname{star}z}a_\varepsilon|\nabla e|^{p-2}\,dx\le\Big(\int_{\operatorname{star}z}a_\varepsilon\,dx\Big)^{2/p}\Big(\int_{\operatorname{star}z}a_\varepsilon|\nabla e|^p\,dx\Big)^{(p-2)/p} \label{eq:p-star-holder}
\end{equation}
Combining this with \eqref{eq:p-coefficient-error},
\(|K|^{p-1}|\operatorname{star}z|\le CH^{dp}\) and \eqref{eq:p-error-poincare},
\begin{align}
 |c_z|^2H^{-2}\int_{\operatorname{star}z}a_\varepsilon|\nabla e|^{p-2}\,dx&\le C(\mathcal R_z^e)^{2/p}\Big(H^{-p}\int_Ka_\varepsilon|e|^p\,dx\Big)^{2/p}\Big(\int_{\operatorname{star}z}a_\varepsilon|\nabla e|^p\,dx\Big)^{(p-2)/p} \label{eq:p-error-correction}\\
 &\le C(\mathcal R^e\mathcal P^e)^{2/p}\int_{N(K)}a_\varepsilon|\nabla e|^p\,dx \label{eq:p-error-correction-bound}
\end{align}
because the powers of \(H\) cancel. At most \(d+1\) hats overlap at a
point, and a simplex is selected by at most \(d+1\) vertices, so
\(\int_Aw|\nabla\Pi_Hv|^2\le C(\mathcal P\mathcal R+(\mathcal P^e\mathcal R^e)^{2/p})\int_Aw|\nabla e|^2\).
The triangle inequality and \(x^{2/p}\le1+x\) complete the proof.
\end{proof}

\begin{prop}[Weighted block contraction]
\label{prop:p-contraction}
Let \((U,U_D)\) be the pair \eqref{eq:localization-pair} for a patch
\(D\), let \(D^+=N^3(D^-)\subset D\) and \(A=D^+\setminus D^-\). With
\begin{equation}
 \kappa_A=\big(\hat\chi_pC_{{\rm cut},p}\Gamma_A(U)\big)^{1/2},\qquad \theta_A=\frac{\kappa_A}{1+\kappa_A} \label{eq:p-contraction-factor}
\end{equation}
we have \(\mathscr D_{D^-}(U_D,U)\le\theta_A\mathscr D_{D^+}(U_D,U)\).
\end{prop}

\begin{proof}
Let \(z_f\) be the corrected cutoff of \(e=U_D-U\in V_f\) with \(S=D^-\),
and let \(E_-=\int_{D^-}\nabla e^TB\nabla e\,dx\) and
\(E_A=\int_A\nabla e^TB\nabla e\,dx\), where \(B\) is the secant tensor
of the pair.
As in the proof of Proposition~\ref{prop:block-contraction},
\(E_-\le E_A^{1/2}(\int_A\nabla z_f^TB\nabla z_f)^{1/2}\). By
Lemma~\ref{lem:p-secant} and Lemma~\ref{lem:p-cutoff},
\begin{equation}
 \int_A\nabla z_f^TB\nabla z_f\,dx\le\hat MC_{{\rm cut},p}\Gamma_A(U)\int_Aw|\nabla e|^2\,dx\le\kappa_A^2E_A \label{eq:p-contraction-chain}
\end{equation}
so \(E_-\le\kappa_AE_A\), and rearrangement completes the proof.
\end{proof}

The factor depends on the ideal state \(U\) only through the
scale-invariant quantities in \(\Gamma_A(U)\), and not on the patch
state. For an inner element set \(S\), iterating
Proposition~\ref{prop:p-contraction} over the annuli
\(A_j=N^{3j}(S)\setminus N^{3(j-1)}(S)\) gives
\begin{equation}
 \mathscr D_S(U_D,U)\le\Big(\prod_{j=1}^n\theta_{A_j}\Big)\mathscr D_D(U_D,U)\le\theta^n\mathscr D_D(U_D,U) \quad \forall n\in\mathbb N_0,\quad N^{3n}(S)\subset D \label{eq:p-interior-decay}
\end{equation}
with \(\theta=\kappa/(1+\kappa)\) and
\(\kappa^2=\hat\chi_pC_{{\rm cut},p}\Gamma_\Omega(U)\). Decay is faster in
regions where the weights vary little. For \(p=2\) and
\(a_\varepsilon=1\), \(\Gamma\) is a mesh constant and
\eqref{eq:p-interior-decay} is Corollary~\ref{cor:interior-decay} with
\(C_{\rm cut}\) replaced by \(C_{{\rm cut},2}\Gamma\).

\subsection{Amplitude, Assembly and Minimizer Error}
\label{subsec:p-energy}

The amplitude estimate \eqref{eq:outer-local-amplitude} used the upper
bound \(M\). Under \eqref{eq:p-structure} the outer energy contains the
error weight \(a_\varepsilon|\nabla e|^{p-2}\), which is absorbed by
Young's inequality.

\begin{lem}[Local amplitude by absorption]
\label{lem:p-amplitude}
For the pair \eqref{eq:localization-pair} on a patch \(D\), with
\(q=q(v_H)\),
\begin{equation}
 \mathscr D_D(U_D,U)\le C_{\rm amp}\hat M\mathcal Q_D(q;U),\qquad \mathcal Q_D(q;U)=\int_D\widehat W_U|\nabla q|^2\,dx+\hat\chi_p^{(p-2)/2}\int_Da_\varepsilon|\nabla q|^p\,dx \label{eq:p-amplitude}
\end{equation}
for \(p>2\), with
\(C_{\rm amp}=2\cdot3^{p-2}\max\{1,(2/p)(2(p-2)3^{p-2}/p)^{(p-2)/2}\}\).
For \(p=2\), \eqref{eq:p-amplitude} holds with \(C_{\rm amp}=1\) and
\(\mathcal Q_D(q;U)=\int_Da_\varepsilon|\nabla q|^2\,dx\).
\end{lem}

\begin{proof}
By \eqref{eq:outer-pythagoras} and Lemma~\ref{lem:p-secant},
\(\mathscr D_D(U_D,U)\le\int_D\nabla q^TB\nabla q\le\hat M\int_Da_\varepsilon\omega(\nabla U,\nabla U_D)|\nabla q|^2\),
where \(B\) is the secant tensor of the pair. This proves the case
\(p=2\). For \(p>2\), \eqref{eq:p-splitting} reduces the claim to the
term \(\int_Da_\varepsilon|\nabla e|^{p-2}|\nabla q|^2\). Young's inequality
with exponents \(p/(p-2)\) and \(p/2\) gives, for every \(\epsilon>0\),
\begin{equation}
 |\nabla e|^{p-2}|\nabla q|^2\le\frac{p-2}p\epsilon^{p/(p-2)}|\nabla e|^p+\frac2p\epsilon^{-p/2}|\nabla q|^p \label{eq:p-young}
\end{equation}
pointwise. By \eqref{eq:p-splitting} and \eqref{eq:p-natural-equivalence},
\(\int_Da_\varepsilon|\nabla e|^p\le\mathscr D_D(U_D,U)/(c_p\hat m)\). The
choice \(\epsilon^{p/(p-2)}=p/(2(p-2)3^{p-2}\hat\chi_p)\) bounds the
contribution of the first term by \(\mathscr D_D(U_D,U)/2\). This finite
term is absorbed, and
\(\epsilon^{-p/2}=(2(p-2)3^{p-2}\hat\chi_p/p)^{(p-2)/2}\) gives
\eqref{eq:p-amplitude}.
\end{proof}

The amplitude depends only on the ideal state and its correction, so the
weighted overlap applies as in Lemma~\ref{lem:outer-amplitude}. The
factor \(\hat\chi_p^{(p-2)/2}\) records the cost of the absorption.

\begin{lem}[Weighted assembly]
\label{lem:p-assembly}
Let \(U\in V_h\) and \(w_i\in V_f\) for \(i\in\mathcal I\). Then
\begin{align}
 \int_\Omega\widehat W_U\Big|\nabla(I-\Pi_H)I_h\sum_i\phi_iw_i\Big|^2\,dx&\le C_{{\rm asm},p}\Gamma_\Omega(U)\sum_i\int_{N(\omega_i)}\widehat W_U|\nabla w_i|^2\,dx \label{eq:p-assembly-weighted}\\
 \int_\Omega a_\varepsilon\Big|\nabla(I-\Pi_H)I_h\sum_i\phi_iw_i\Big|^p\,dx&\le C_{{\rm asm},p}\Gamma_\Omega(U)\sum_i\int_{N(\omega_i)}a_\varepsilon|\nabla w_i|^p\,dx \label{eq:p-assembly-lp}
\end{align}
with \(C_{{\rm asm},p}\) depending only on \(p\),
\(d\) and the mesh family.
\end{lem}

\begin{proof}
Let \(S_iw_i=(I-\Pi_H)I_h(\phi_iw_i)\). The product \(I_h(\phi_iw_i)\)
is supported in \(\omega_i\). A projection coefficient
\(\ell_z(I_h(\phi_iw_i))\) can be nonzero only if \(K_z\subset\omega_i\),
and then \(\operatorname{supp}\lambda_z\subset N(K_z)\subset N(\omega_i)\). Thus
\(S_iw_i\) vanishes outside \(N(\omega_i)\), and a fine element meets boundedly
many \(N(\omega_i)\). It suffices to bound the weighted integrals of
\(\nabla S_iw_i\) over \(N(\omega_i)\). On a fine element \(t\subset T\subset\omega_i\),
\eqref{eq:product-gradient}, \(0\le\phi_i\le1\) and
\(|\nabla\phi_i|\le CH^{-1}\) give
\(\int_t\widehat W_U|\nabla I_h(\phi_iw_i)|^2\le C\int_t\widehat W_U(|\nabla w_i|^2+H^{-2}w_i^2)\),
and \eqref{eq:p-weighted-poincare} on \(T\), with \(N(T)\subset N(\omega_i)\),
bounds the mass term. For the projection, the coefficient bound
\eqref{eq:p-coefficient-background} with \(g=I_h(\phi_iw_i)\) and
\eqref{eq:p-background-correction} on \(K_z\subset\omega_i\) give the
factor \(\mathcal R_z\mathcal P_{K_z}\), since \(N(K_z)\subset N(\omega_i)\). This
proves \eqref{eq:p-assembly-weighted}. The \(L^p\) estimate follows in
the same way. The product estimate holds elementwise in \(L^p\) by
\eqref{eq:p-lq-product}. The
coefficient bound is \eqref{eq:p-coefficient-error}, and
\(\int a_\varepsilon|\nabla\lambda_z|^p\le CH^{-p}\int_{\operatorname{star}z}a_\varepsilon\).
The powers of \(H\) cancel as before and give the factor
\(1+\mathcal P^e(1+\mathcal R^e)\le\Gamma_\Omega(U)\).
\end{proof}

At the coarse coordinate of the fine solution, the ideal state is
\(u_h\) and its correction is \(q^*=(I-\Pi_H)u_h\). For fixed patches with
\(N(\omega_i)\subset D_i\), choose for each proper patch an
integer \(n_i\ge0\) with \(N^{3n_i+r_P}(\omega_i)\subset D_i\) and set
\begin{equation}
 t_i=\theta^{n_i},\qquad \theta=\frac\kappa{1+\kappa},\qquad \kappa^2=\hat\chi_pC_{{\rm cut},p}\Gamma_\Omega(u_h) \label{eq:p-uniform-factor}
\end{equation}
with \(t_i=0\) for full-domain patches. The weighted overlap
\(\mathcal N_t(\mathcal D)\) is \eqref{eq:weighted-overlap} with these
weights.

\begin{thm}[Energy and gradient error for $p$-structure]
\label{thm:p-energy}
Assume \eqref{eq:p-structure}, \ref{ass:mesh}--\ref{ass:patches} of
Assumption~\ref{ass:standing}, and the patch weights
\eqref{eq:p-uniform-factor}. A global minimizer \(u_{\rm ms}\) of
\(\mathcal E\) on \(M_h^{\mathcal D}(V_H)\) exists. Every such minimizer
satisfies, with \(\delta=u_{\rm ms}-u_h\),
\begin{align}
 0\le\mathcal E(u_{\rm ms})-\mathcal E(u_h)&\le C_E\hat\chi_p\Gamma_\Omega(u_h)\hat M\mathcal N_t(\mathcal D)\mathcal Q_\Omega(q^*;u_h) \label{eq:p-energy}\\
 \frac{c_B\hat m}2\int_\Omega a_\varepsilon\big(|\nabla\delta|^2+|\nabla\delta|^p\big)\,dx&\le\mathcal E(u_{\rm ms})-\mathcal E(u_h) \label{eq:p-gradient}
\end{align}
with \(C_E=2\cdot3^{p-2}C_{{\rm asm},p}C_{\rm amp}\). Moreover,
\begin{equation}
 \mathcal Q_\Omega(q^*;u_h)\le C_Q\Gamma_\Pi(u_h)\hat\chi_p^{(p-2)/2}\frac{\hat L_F^2}{\hat m^2} \label{eq:p-amplitude-data}
\end{equation}
where \(C_Q\) depends only on \(p\).
\end{thm}

\begin{proof}
Existence and the comparison
\(\mathcal E(u_{\rm ms})-\mathcal E(u_h)\le\mathcal B(w_*,u_h)\), with
\(w_*=M_h^{\mathcal D}(\Pi_Hu_h)\), follow from the first part of
Theorem~\ref{thm:energy-comparison}. Both Bregman terms in
\eqref{eq:bregman-symmetrization} are nonnegative, so
\(\mathcal B(w_*,u_h)\le\mathscr D_\Omega(w_*,u_h)\). By
\eqref{eq:exact-comparator-coordinate} and \eqref{eq:exact-assembly-error},
\(\delta_*=w_*-u_h=(I-\Pi_H)I_h\sum_i\phi_ie_i\) with
\(e_i=q_i(\Pi_Hu_h)-q^*\). Lemma~\ref{lem:p-secant} and
\eqref{eq:p-splitting} with \(U=u_h\) give
\begin{equation}
 \mathscr D_\Omega(w_*,u_h)\le3^{p-2}\hat M\Big(\int_\Omega\widehat W_{u_h}|\nabla\delta_*|^2\,dx+\int_\Omega a_\varepsilon|\nabla\delta_*|^p\,dx\Big) \label{eq:p-energy-split}
\end{equation}
for \(p>2\), and the same bound without the second term for \(p=2\).
Lemma~\ref{lem:p-assembly} bounds the right side by
\(3^{p-2}\hat MC_{{\rm asm},p}\Gamma_\Omega(u_h)\sum_i\int_{N(\omega_i)}(\widehat W_{u_h}|\nabla e_i|^2+a_\varepsilon|\nabla e_i|^p)\).
Let \(U_i=\Pi_Hu_h+q_i(\Pi_Hu_h)\), so that \(U_i-u_h=e_i\). By the
lower bounds in \eqref{eq:p-splitting} and \eqref{eq:p-natural-equivalence},
each summand is at most \(2\mathscr D_{N(\omega_i)}(U_i,u_h)/(c_p\hat m)\).
Interior decay \eqref{eq:p-interior-decay} with \(U=u_h\) gives
\(\mathscr D_{N(\omega_i)}(U_i,u_h)\le t_i\mathscr D_{D_i}(U_i,u_h)\), and the
same holds with \(t_i=0\) for full-domain patches.
Lemma~\ref{lem:p-amplitude} and the definition
\eqref{eq:weighted-overlap} of the weighted overlap
give \(\sum_it_i\mathscr D_{D_i}(U_i,u_h)\le C_{\rm amp}\hat M\mathcal N_t(\mathcal D)\mathcal Q_\Omega(q^*;u_h)\).
This proves \eqref{eq:p-energy}.

The identity \eqref{eq:energy-bregman}, the lower bound
\eqref{eq:p-bregman-lower} and
\(\omega(\nabla u_{\rm ms},\nabla u_h)\ge\max\{1,|\nabla\delta|^{p-2}\}\)
prove \eqref{eq:p-gradient}. For \eqref{eq:p-amplitude-data}, the
definition of \(\Gamma_\Pi\) and \(\Gamma_\Pi\ge1\) give
\begin{equation}
 \int_\Omega\widehat W_{u_h}|\nabla q^*|^2\,dx\le4\Gamma_\Pi(u_h)\int_\Omega\widehat W_{u_h}|\nabla u_h|^2\,dx,\qquad \int_\Omega a_\varepsilon|\nabla q^*|^p\,dx\le2^p\Gamma_\Pi(u_h)\int_\Omega a_\varepsilon|\nabla u_h|^p\,dx \label{eq:p-correction-stability}
\end{equation}
The elementary bound
\((1+s)^{p-2}s^2\le\max\{1,2^{p-3}\}(s^2+s^p)\), \(\hat\chi_p\ge1\) and
\eqref{eq:p-data-p} complete the proof.
\end{proof}

For uniform patches \(D_i=N^\ell(\omega_i)\), \(\ell\ge r_P\), the
packing bound \eqref{eq:uniform-weighted-overlap} with
\(n_\ell=\lfloor(\ell-1)/3\rfloor\) gives
\(\mathcal N_t\le C_{\rm ol}(1+\ell)^d\theta^{n_\ell}\). We may take
\(C_{{\rm cut},p}\ge1\). Then \(\kappa\ge1\), since \(c_p\le1\le\hat M/\hat m\)
and \(\Gamma\ge1\). Since \(n_\ell\ge\ell/3-1\), we obtain
\(\theta^{n_\ell}\le\theta^{-1}e^{-b\ell}\le2e^{-b\ell}\).
Hence
\begin{equation}
 \mathcal E(u_{\rm ms}^\ell)-\mathcal E(u_h)\le C\hat\chi_p^{p/2}\Gamma_\Omega(u_h)\Gamma_\Pi(u_h)(1+\ell)^de^{-b\ell}\frac{\hat M\hat L_F^2}{\hat m^2},\qquad b=\frac13\log\Big(1+\frac1\kappa\Big) \label{eq:p-rate}
\end{equation}
with \(C\) depending only on \(p\), \(d\) and the mesh family. By
\eqref{eq:p-gradient}, the \(a_\varepsilon\)-weighted gradient error
decays like \((1+\ell)^{d/2}e^{-b\ell/2}\) in \(L^2\) and like
\((1+\ell)^{d/p}e^{-b\ell/p}\) in \(L^p\). The energy rate \(b\) equals
the decay rate of the patch energies; there is no loss in the exponent.
The data enter the rate only through \(\Gamma_\Omega(u_h)\), and the
prefactor additionally through \(\Gamma_\Pi(u_h)\).

\begin{cor}[Bounded Hessians as the case $p=2$]
\label{cor:p-two}
Let \(p=2\). Then \(\widehat W_U=a_\varepsilon\), the quantities in
\(\Gamma_\Omega\) and \(\Gamma_\Pi\) depend on \(a_\varepsilon\) only, and
Theorem~\ref{thm:p-energy} gives
\begin{equation}
 \mathcal E(u_{\rm ms}^\ell)-\mathcal E(u_h)\le C\frac{\hat M}{\hat m}\Gamma_\Omega\Gamma_\Pi(1+\ell)^de^{-b\ell}\frac{\hat M\hat L_F^2}{\hat m^2},\qquad b=\frac13\log\Big(1+\big(\tfrac{\hat M}{\hat m}C_{{\rm cut},2}\Gamma_\Omega\big)^{-1/2}\Big) \label{eq:p-two-rate}
\end{equation}
for every \(\ell\ge r_P\), independently of the load. For
\(a_\varepsilon=1\), \(\hat M/\hat m=\chi\) and \(\Gamma_\Omega\), \(\Gamma_\Pi\)
are mesh constants, so \eqref{eq:p-two-rate} has the form of
\eqref{eq:exponential-energy-form}. For a general coefficient, the
contrast enters only through the scale-invariant quantities
\(\mathcal P_T\), \(\mathcal R_z\) and \(\Gamma_\Pi\) of \(a_\varepsilon\),
and not through \(a_+/a_-\).
\end{cor}

\begin{proof}
For \(p=2\), \(\hat\chi_2=\hat M/\hat m\) and \(\Gamma_\Omega(u_h)\),
\(\Gamma_\Pi(u_h)\) do not depend on \(u_h\). Insert these into
\eqref{eq:p-rate}.
\end{proof}

\begin{rem}[Rate loss without weighted assembly]
\label{rem:p-rate-loss}
Bounding \(\int\widehat W_{u_h}|\nabla\delta_*|^2\) in \eqref{eq:p-energy-split}
by H\"older's inequality and the \(L^p\) assembly estimate alone would
give a factor \(\mathcal N_t^{2/p}\) and reduce the exponential rate by
the factor \(2/p\). Lemma~\ref{lem:p-assembly} avoids this loss at the
price of the scale-invariant factor \(\Gamma_\Omega(u_h)\).
\end{rem}

\subsection{Computed Coarse States}
\label{subsec:p-computed}

Theorem~\ref{thm:computed-state-error} bounds the error at every coarse
coordinate by its reduced residual and a localization term. Its proof
uses the flux Lipschitz bound \eqref{eq:flux-lipschitz} and the tangent
stability of Lemma~\ref{lem:localized-tangent-stability}, both of which
rely on \(M\). Under \eqref{eq:p-structure}, the first is replaced by
the Cauchy--Schwarz inequality in the secant inner product. The second
becomes a weighted estimate that requires the patch states to be
compatible with the ideal state at the computed coordinate.

\begin{lem}[A metric equivalent to the natural distance]
\label{lem:p-v-function}
Let \(V(\xi)=(1+|\xi|)^{(p-2)/2}\xi\), \(c_V=2^{4-p}/p^2\) and
\(C_V=p^2/4\). Then
\begin{equation}
 c_V\,\omega(\xi,\eta)|\xi-\eta|^2\le|V(\xi)-V(\eta)|^2\le C_V\,\omega(\xi,\eta)|\xi-\eta|^2 \quad \forall\xi,\eta\in\mathbb R^d \label{eq:p-v-function}
\end{equation}
Consequently
\(d_V(v,w)=(\int_\Omega a_\varepsilon|V(\nabla v)-V(\nabla w)|^2\,dx)^{1/2}\)
is a metric on \(V_h\) with
\begin{equation}
 \frac{c_p\hat m}{C_V}d_V(v,w)^2\le\mathscr D_\Omega(v,w)\le\frac{\hat M}{c_V}d_V(v,w)^2 \quad \forall v,w\in V_h \label{eq:p-v-metric}
\end{equation}
\end{lem}

\begin{proof}
With \(\gamma=(p-2)/2\) and \(\hat\xi=\xi/|\xi|\),
\(DV(\xi)=(1+|\xi|)^\gamma[I+\gamma|\xi|(1+|\xi|)^{-1}\hat\xi\otimes\hat\xi]\)
is symmetric with eigenvalues between \((1+|\xi|)^\gamma\) and
\((1+\gamma)(1+|\xi|)^\gamma\). Hence \(V(\eta)-V(\xi)=\mathsf M(\eta-\xi)\),
where \(\mathsf M=\int_0^1DV(\xi_t)\,dt\) is symmetric with eigenvalues
between \(\mu=\int_0^1(1+|\xi_t|)^\gamma\,dt\) and \((1+\gamma)\mu\).
The proof of \eqref{eq:p-scalar-secant}, with exponent \(\gamma\) in place
of \(p-2\), gives
\(2^{-\gamma}(1+\gamma)^{-1}(1+|\xi|+|\eta|)^\gamma\le\mu\le(1+|\xi|+|\eta|)^\gamma\).
Squaring proves \eqref{eq:p-v-function}. The map \(V\) is injective,
so \(d_V\) is the metric induced by the weighted \(L^2\) norm.
Combining \eqref{eq:p-v-function} with \eqref{eq:p-natural-equivalence}
proves \eqref{eq:p-v-metric}.
\end{proof}

The function \(V\) and the equivalence \eqref{eq:p-v-function} are standard
for problems with \(p\)-structure, see~\cite{diening-ettwein2008}; the
lemma records explicit constants.

Fix a coarse coordinate \(v_H\) and write
\(U=M_h(v_H)\), \(U_i=v_H+q_i(v_H)\), \(\widetilde u=M_h^{\mathcal D}(v_H)\)
and \(T_v=DM_h^{\mathcal D}(v_H)\). We define the compatibility quantities
\begin{align}
 \lambda_C&=\max_i\max_{t\subset N(\omega_i)}\frac{|\nabla(U_i-U)|}{1+|\nabla U|},\qquad \lambda_{\rm rec}=\max_{t\subset\Omega}\frac{|\nabla(\widetilde u-U)|}{1+|\nabla U|} \label{eq:p-compatibility-pointwise}\\
 \varrho_D&=\max_i\max_{T\subset D_i}\frac{\int_T\widehat W_{U_i}\,dx}{\int_T\widehat W_U\,dx} \label{eq:p-compatibility-average}
\end{align}
where \(t\) ranges over fine elements and \(T\) over coarse elements.
The first quantity is pointwise but restricted to the active supports.
The third compares averages over coarse elements in the whole patches.
The second is pointwise on the domain; it enters below only through the
factor \((2+\lambda_{\rm rec})^{p-2}\) and need not be small. The interior decay \eqref{eq:p-interior-decay} controls \(U_i-U\) on \(N(\omega_i)\) in energy, not pointwise, so
\(\lambda_C<1\) is a hypothesis at the computed coordinate. For \(p=2\)
all three quantities enter only through powers \(p-2=0\) or through
\(\varrho_D=1\).

\begin{lem}[Weighted tangent stability]
\label{lem:p-tangent}
If \(\lambda_C<1\), then for every \(\tau_H\in V_H\),
\begin{equation}
 \int_\Omega\widehat W_U|\nabla T_v\tau_H|^2\,dx\le C_{\rm tan}\int_\Omega\widehat W_U|\nabla\tau_H|^2\,dx,\qquad C_{\rm tan}=2+2C_{{\rm asm},p}\Gamma_\Omega(U)\frac{\hat M}{\hat m}\frac{\varrho_D}{(1-\lambda_C)^{p-2}}\mathcal N_0(\mathcal D) \label{eq:p-tangent}
\end{equation}
with the unweighted overlap \(\mathcal N_0(\mathcal D)\) of
\eqref{eq:tangent-stability-constant}.
\end{lem}

\begin{proof}
Let \(s_i=Dq_i(v_H)[\tau_H]\) and \(A_i=D_\xi a(x,\nabla U_i)\). As in
\eqref{eq:patch-tangent-pythagoras}, the \(A_i\)-energy of \(s_i\) on
\(D_i\) is at most that of \(\tau_H\). Since \(\tau_H\) has a constant
gradient on each coarse element,
\begin{equation}
 \int_{D_i}\nabla\tau_H^TA_i\nabla\tau_H\,dx\le\hat M\sum_{T\subset D_i}|\nabla\tau_H|_T|^2\int_T\widehat W_{U_i}\,dx\le\hat M\varrho_D\int_{D_i}\widehat W_U|\nabla\tau_H|^2\,dx \label{eq:p-tangent-upper}
\end{equation}
On \(N(\omega_i)\), \(1+|\nabla U_i|\ge(1-\lambda_C)(1+|\nabla U|)\), so
\(A_i\ge\hat m(1-\lambda_C)^{p-2}\widehat W_U\) there. Hence
\begin{equation}
 \int_{N(\omega_i)}\widehat W_U|\nabla s_i|^2\,dx\le\frac{\hat M}{\hat m}\frac{\varrho_D}{(1-\lambda_C)^{p-2}}\int_{D_i}\widehat W_U|\nabla\tau_H|^2\,dx \label{eq:p-tangent-correction-bound}
\end{equation}
By \eqref{eq:assembled-tangent}, \(T_v\tau_H-\tau_H=(I-\Pi_H)I_h\sum_i\phi_is_i\)
with \(s_i\in V_f\). Lemma~\ref{lem:p-assembly} and
\(\sum_i\int_{D_i}\le\mathcal N_0(\mathcal D)\int_\Omega\) bound
\(\int\widehat W_U|\nabla(T_v\tau_H-\tau_H)|^2\), and the triangle
inequality completes the proof.
\end{proof}

Near the artificial patch boundary, a patch state may differ from the
ideal state by the size of the fine correction. Because tangent
directions are coarse, only the coarse averages \(\varrho_D\) enter
there. We measure the reduced residual in the \(a_\varepsilon\)-weighted
coarse dual norm,
\begin{equation}
 \rho_{H,a}(v_H)=\sup_{0\ne\tau_H\in V_H}\frac{|(\mathcal A(\widetilde u)-F)[T_v\tau_H]|}{\big(\int_\Omega a_\varepsilon|\nabla\tau_H|^2\,dx\big)^{1/2}} \label{eq:p-reduced-residual}
\end{equation}
and set \(\mathscr D_{\rm id}=\mathscr D_\Omega(U,u_h)\) and
\(\mathscr D_{\rm loc}=\mathscr D_\Omega(U,\widetilde u)\).

\begin{prop}[Residual transfer]
\label{prop:p-transfer}
If \(\lambda_C<1\), then
\begin{equation}
 \mathscr D_{\rm id}\le\frac{2\Gamma_\Pi(U)}{c_p\hat m}\rho_{H,a}(v_H)^2+K_v\mathscr D_{\rm loc},\qquad K_v=2\hat\chi_p(2+\lambda_{\rm rec})^{p-2}C_{\rm tan}\Gamma_\Pi(U) \label{eq:p-transfer}
\end{equation}
\end{prop}

\begin{proof}
Let \(e=U-u_h\), \(\tau_H=\Pi_He\) and \(r_U=\mathcal A(U)-F\). As in
\eqref{eq:ideal-residual-cancellation}, \(e-T_v\tau_H\in V_f\) and
\(r_U[e]=r_U[T_v\tau_H]\). The fine equation \eqref{eq:fine-equation} gives \(r_U[e]=\mathscr D_{\rm id}\), so
\begin{equation}
 \mathscr D_{\rm id}=(\mathcal A(\widetilde u)-F)[T_v\tau_H]+(\mathcal A(U)-\mathcal A(\widetilde u))[T_v\tau_H] \label{eq:p-transfer-identity}
\end{equation}
By \eqref{eq:p-natural-equivalence} and \eqref{eq:p-splitting}, both
\(\int a_\varepsilon|\nabla e|^2\) and \(\int\widehat W_U|\nabla e|^2\) are at
most \(\mathscr D_{\rm id}/(c_p\hat m)\). With \eqref{eq:p-projection-coefficient},
the first term is at most \(\rho_{H,a}(\Gamma_\Pi\mathscr D_{\rm id}/(c_p\hat m))^{1/2}\).
The second equals \(\int\nabla(T_v\tau_H)^T\widetilde B\nabla(U-\widetilde u)\),
where \(\widetilde B\) is the secant tensor of \(U\) and \(\widetilde u\).
By the Cauchy--Schwarz inequality in the \(\widetilde B\)-inner product, it
is at most \(\mathscr D_{\rm loc}^{1/2}\) times the \(\widetilde B\)-norm of \(T_v\tau_H\).
Since \(1+|\nabla U|+|\nabla\widetilde u|\le(2+\lambda_{\rm rec})(1+|\nabla U|)\),
Lemma~\ref{lem:p-secant} gives
\(\widetilde B\le\hat M(2+\lambda_{\rm rec})^{p-2}\widehat W_U\). Lemma~\ref{lem:p-tangent}
and \eqref{eq:p-projection-weighted} bound the second term by
\((\mathscr D_{\rm loc}K_v\mathscr D_{\rm id}/2)^{1/2}\). Dividing by \(\mathscr D_{\rm id}^{1/2}\) and squaring
complete the proof.
\end{proof}

Proposition~\ref{prop:p-transfer} uses neither global minimality nor a
Lipschitz bound on the flux.

\begin{thm}[Computed-state error for $p$-structure]
\label{thm:p-computed}
Assume \eqref{eq:p-structure} and \ref{ass:mesh}--\ref{ass:patches} of
Assumption~\ref{ass:standing}, with \(N(\omega_i)\subset D_i\). Let
\(v_H\in V_H\) satisfy \(\lambda_C<1\), and let the weights \(t_i\) be
defined by \eqref{eq:p-uniform-factor} with \(\Gamma_\Omega(U)\) in place
of \(\Gamma_\Omega(u_h)\). There is a constant \(C\), depending only on
\(p\), \(d\) and the mesh family, such that
\begin{equation}
 \beta_v=CC_EK_v\hat\chi_p^{(p+2)/2}\Gamma_\Omega(U)\Gamma_\Pi(U)\mathcal N_t(\mathcal D)\le\frac12 \label{eq:p-oversampling}
\end{equation}
implies
\begin{align}
 d_V(\widetilde u,u_h)^2&\le\frac{CC_V}{c_p\hat m}\Big(\frac{\Gamma_\Pi(U)}{c_p\hat m}\rho_{H,a}(v_H)^2+C_EK_v\hat\chi_p\Gamma_\Omega(U)\hat M\mathcal N_t(\mathcal D)\mathcal Q_\Omega(q^*;u_h)\Big) \label{eq:p-computed}\\
 0\le\mathcal E(\widetilde u)-\mathcal E(u_h)&\le\frac{\hat M}{c_V}d_V(\widetilde u,u_h)^2 \label{eq:p-computed-energy}
\end{align}
In particular, every tangent-stationary coordinate with \(\lambda_C<1\)
that satisfies \eqref{eq:p-oversampling} obeys \eqref{eq:p-computed}
with \(\rho_{H,a}(v_H)=0\).
\end{thm}

\begin{proof}
The proof of Theorem~\ref{thm:p-energy}, applied at \(v_H\) with the
ideal state \(U\) in place of \(u_h\), gives
\(\mathscr D_{\rm loc}\le C_E\hat\chi_p\Gamma_\Omega(U)\hat M\mathcal N_t\mathcal Q_\Omega(q;U)\)
with \(q=q(v_H)\). Write \(q=(I-\Pi_H)e+q^*\) with \(e=U-u_h\). By the
definition of \(\Gamma_\Pi\), the integrals of \(\widehat W_U|\nabla(I-\Pi_H)e|^2\)
and \(a_\varepsilon|\nabla(I-\Pi_H)e|^p\) are at most \(2^p\Gamma_\Pi(U)\mathscr D_{\rm id}/(c_p\hat m)\).
Moreover,
\(\widehat W_U\le\max\{1,2^{p-3}\}(\widehat W_{u_h}+a_\varepsilon|\nabla e|^{p-2})\), and
Young's inequality bounds \(\int a_\varepsilon|\nabla e|^{p-2}|\nabla q^*|^2\) by
\(\int a_\varepsilon(|\nabla e|^p+|\nabla q^*|^p)\). Therefore
\begin{equation}
 \mathcal Q_\Omega(q;U)\le C\Gamma_\Pi(U)\hat\chi_p^{(p-2)/2}\frac{\mathscr D_{\rm id}}{c_p\hat m}+C\mathcal Q_\Omega(q^*;u_h) \label{eq:p-computed-amplitude}
\end{equation}
Insert the resulting bound for \(\mathscr D_{\rm loc}\) into \eqref{eq:p-transfer}. The
coefficient of \(\mathscr D_{\rm id}\) on the right is at most \(\beta_v\le1/2\). With
\(\mathcal Q^*=\mathcal Q_\Omega(q^*;u_h)\), this gives
\begin{align}
 \mathscr D_{\rm id}&\le\frac{4\Gamma_\Pi(U)}{c_p\hat m}\rho_{H,a}^2+CC_EK_v\hat\chi_p\Gamma_\Omega(U)\hat M\mathcal N_t\mathcal Q^* \label{eq:p-computed-ideal}\\
 \mathscr D_{\rm loc}&\le\mathscr D_{\rm id}+CC_E\hat\chi_p\Gamma_\Omega(U)\hat M\mathcal N_t\mathcal Q^* \label{eq:p-computed-localization}
\end{align}
where the second bound uses \(\beta_v/K_v\le1\), since
\(K_v\ge2\hat\chi_pC_{\rm tan}\Gamma_\Pi\ge4\). By Lemma~\ref{lem:p-v-function},
\(d_V(\widetilde u,u_h)^2\le2d_V(U,u_h)^2+2d_V(U,\widetilde u)^2\), and
both terms are at most \(C_V/(c_p\hat m)\) times \(\mathscr D_{\rm id}\) and \(\mathscr D_{\rm loc}\),
which proves \eqref{eq:p-computed}. The energy bound follows from
\(\mathcal E(\widetilde u)-\mathcal E(u_h)=\mathcal B(\widetilde u,u_h)\le\mathscr D_\Omega(\widetilde u,u_h)\)
and \eqref{eq:p-v-metric}.
\end{proof}

The constants \(C\) and \(C_E\) depend only on \(p\), \(d\) and the mesh
family. The structural contrast enters through \(\hat\chi_p\), \(K_v\) and
the factor \(\hat\chi_p^{(p-2)/2}\) in \(\mathcal Q^*\). For \(p=2\), the
hypothesis \(\lambda_C<1\) can be dropped, since the proof uses it only
through \((1-\lambda_C)^{p-2}\). For uniform patches,
\(\mathcal N_0\le C(1+\ell)^d\) and
\(\mathcal N_t\le C(1+\ell)^d\theta^{n_\ell}\), so
\(K_v\mathcal N_t\lesssim(1+\ell)^{2d}\theta^{n_\ell}\). As in
Corollary~\ref{cor:stationary-root-rate}, the polynomial factor is one
power of \((1+\ell)^d\) larger than for minimizers. The oversampling
threshold depends, besides the geometry and the structural constants, on
the weighted quantities and on the compatibility quantities at the
computed coordinate. For \(p=2\) the compatibility quantities drop out,
and Theorem~\ref{thm:p-computed} has the structure of
Theorem~\ref{thm:computed-state-error}, with the contrast in weighted
rather than global constants. The bound
\(K_v\mathcal N_t\lesssim(1+\ell)^{2d}\theta^{n_\ell}\) above assumes that
\(\lambda_C\), \(\lambda_{\rm rec}\) and \(\varrho_D\) are bounded independently of
\(\ell\). For \(d=2\) or \(p\le6\), Proposition~\ref{prop:p-apriori} and
Corollary~\ref{cor:p-apriori} show that \(\lambda_C<1\) follows from
oversampling for coordinates with bounded energy.


\section{Size of the Weighted Constants}
\label{sec:p-constants}

Theorems~\ref{thm:p-energy} and~\ref{thm:p-computed} contain the
quantities \(\Gamma_\Omega\) and \(\Gamma_\Pi\) of an ideal state, and in
particular of the unknown fine solution. They are computable a
posteriori. This section bounds them independently of the fine mesh,
for \(d=2\) and for \(d=3\) with \(p\le6\), weakens the \(A_2\)-type
condition contained in \(\mathcal R_z\), and separates the size of the load
from the shape of the gradient in the ratios. It then gives conditions for
load-uniform constants, derives the compatibility hypothesis of
Theorem~\ref{thm:p-computed} from oversampling, and ends with tensor
versions of the contraction and the energy estimate. Throughout,
\(\kappa_a=a_+/a_-\) and \(c_2=\max\{1,2^{p-3}\}\), so that
\((1+t)^{p-2}\le c_2(1+t^{p-2})\) for \(t\ge0\).

\subsection{Crude Bounds}
\label{subsec:p-crude}

\begin{lem}[Unweighted \texorpdfstring{$L^k$}{Lk} estimates]
\label{lem:p-lq}
Let \(1\le k<\infty\). There are constants \(C_{P,k}\) and \(C_{\Pi,k}\),
depending only on \(k\), \(d\) and the mesh family, such that
\begin{align}
 H^{-1}\|v_f\|_{L^k(T)}&\le C_{P,k}\|\nabla v_f\|_{L^k(N(T))} \quad \forall v_f\in V_f \label{eq:p-lq-kernel}\\
 \|\nabla\Pi_Hv\|_{L^k(\Omega)}&\le C_{\Pi,k}\|\nabla v\|_{L^k(\Omega)} \quad \forall v\in V_h \label{eq:p-lq-stability}
\end{align}
for every coarse element \(T\). For a fine element \(t\), a coarse
piecewise-affine function \(\eta\) and \(w\in V_h\),
\begin{equation}
 \|\nabla I_h(\eta w)\|_{L^k(t)}\le C\big(\|\eta\|_{L^\infty(t)}\|\nabla w\|_{L^k(t)}+\|\nabla\eta\|_{L^\infty(t)}\|w\|_{L^k(t)}\big) \label{eq:p-lq-product}
\end{equation}
If \(2\le p<\infty\) for \(d=2\), or \(2\le p\le6\) for \(d=3\), there is
\(C_S\) such that
\begin{equation}
 H^{-1}\|v_f\|_{L^p(T)}\le C_SH^{d/p-d/2}\|\nabla v_f\|_{N(T)} \quad \forall v_f\in V_f \label{eq:p-sobolev-poincare}
\end{equation}
\end{lem}

\begin{proof}
Constants \(C\) depend only on \(k\), \(p\), \(d\) and the mesh family.
Write \(v_K=|K|^{-1}\int_Kv\,dx\).

\textbf{Patch inequalities in \(L^k\).} Let \(K\) and \(L\) be coarse
elements sharing a face \(F\), let \(v\in W^{1,k}(K\cup L)\), and let
\(v_F\) be the mean of its trace on \(F\). The trace inequality
\(\int_F|w|\le C(H^{-1}\|w\|_{L^1(K)}+\|\nabla w\|_{L^1(K)})\) on a
shape-regular simplex, applied to \(w=v-v_K\), the \(L^1\) Poincar\'e
inequality on \(K\) and H\"older's inequality give
\(|v_F-v_K|\le|F|^{-1}\int_F|v-v_K|\le CH^{1-d}\|\nabla v\|_{L^1(K)}\le CH^{1-d/k}\|\nabla v\|_{L^k(K)}\).
The same bound on \(L\) gives
\begin{equation}
 |v_K-v_L|\le CH^{1-d/k}\|\nabla v\|_{L^k(K\cup L)} \label{eq:p-lq-neighbour}
\end{equation}
By the proof of Lemma~\ref{lem:one-layer-patches}, every element of
\(N(T)\) is joined to \(T\) by a face chain of bounded length. Summing
\eqref{eq:p-lq-neighbour} along these chains and using the simplex
Poincar\'e inequality \(\|v-v_K\|_{L^k(K)}\le CH\|\nabla v\|_{L^k(K)}\)
gives
\begin{equation}
 \|v-v_T\|_{L^k(N(T))}\le CH\|\nabla v\|_{L^k(N(T))} \quad \forall v\in W^{1,k}(N(T)) \label{eq:p-lq-patch}
\end{equation}
If \(T\) has a boundary vertex, \(N(T)\) contains an element \(K_0\) with a
physical boundary facet of measure comparable to \(H^{d-1}\). For \(v\)
vanishing on that facet, the trace mean is zero, so
\(|v_{K_0}|\le CH^{1-d/k}\|\nabla v\|_{L^k(K_0)}\), and chaining from
\(K_0\) gives \(\|v\|_{L^k(N(T))}\le CH\|\nabla v\|_{L^k(N(T))}\).

\textbf{Projection.} By \eqref{eq:local-dual-formula},
\(|\theta_z^K|\le C|K|^{-1}\), so
\(|\ell_z(v)|\le C|K_z|^{-1/k}\|v\|_{L^k(K_z)}\), and
\(|\nabla\lambda_z|\le CH^{-1}\). Since \(K_z\subset\operatorname{star}z\subset N(T)\)
for every vertex \(z\) of \(T\),
\(\|\Pi_Hv\|_{L^k(T)}+H\|\nabla\Pi_Hv\|_{L^k(T)}\le C\|v\|_{L^k(N(T))}\). If
every vertex of \(T\) is interior, \(\Pi_H\) reproduces constants on \(T\)
by \eqref{eq:local-constant-reproduction}. Applying the last bound to
\(v-v_T\) and using \eqref{eq:p-lq-patch} gives
\(\|\nabla\Pi_Hv\|_{L^k(T)}+H^{-1}\|v-\Pi_Hv\|_{L^k(T)}\le C\|\nabla v\|_{L^k(N(T))}\).
If \(T\) has a boundary vertex, the same bound follows from the zero-trace
inequality applied to \(v\) itself. For \(v=v_f\in V_f\), \(\Pi_Hv_f=0\), which
proves \eqref{eq:p-lq-kernel}. Raising the gradient bound to the power
\(k\) and summing with bounded overlap of the patches \(N(T)\) proves
\eqref{eq:p-lq-stability}.

\textbf{Products.} On \(t\), the edge identity
\eqref{eq:product-edge-identity} and fine-element scaling give
\begin{equation}
 |\nabla I_h(\eta w)|\le C\big(\|\eta\|_{L^\infty(t)}|\nabla w|+\|\nabla\eta\|_{L^\infty(t)}\max_t|w|\big) \label{eq:p-lq-product-pointwise}
\end{equation}
where \(\nabla w\) is constant on \(t\). Since
\(\max_t|w|\le C|t|^{-1/k}\|w\|_{L^k(t)}\) for affine \(w\), this proves
\eqref{eq:p-lq-product}.

\textbf{Sobolev--Poincar\'e inequality.} On a shape-regular simplex \(K\),
scaling the embedding \(H^1\subset L^p\) on the reference simplex, valid
in the stated range of \(p\), and the Poincar\'e inequality give
\(\|v-v_K\|_{L^p(K)}\le CH^{1+d/p-d/2}\|\nabla v\|_K\). The neighbour
estimate \eqref{eq:neighbor-mean-control} gives
\(|K|^{1/p}|v_K-v_T|\le CH^{1+d/p-d/2}\|\nabla v\|_{N(T)}\) for
\(K\subset N(T)\). Hence
\(\|v-v_T\|_{L^p(N(T))}\le CH^{1+d/p-d/2}\|\nabla v\|_{N(T)}\), and, by the
zero-trace argument above, \(\|v\|_{L^p(N(T))}\) obeys the same bound if
\(T\) has a boundary vertex and \(v\) vanishes on \(\partial\Omega\). Let
\(v_f\in V_f\). If \(T\) has a boundary vertex, this proves
\eqref{eq:p-sobolev-poincare}. Otherwise, let \(z\) be a vertex of \(T\).
Since \(\ell_z(v_f)=0\) and \(\ell_z(1)=1\), we have
\(v_T=-\ell_z(v_f-v_T)\), so
\(|v_T|\le C|K_z|^{-1/p}\|v_f-v_T\|_{L^p(K_z)}\). Therefore
\(\|v_f\|_{L^p(T)}\le C\|v_f-v_T\|_{L^p(N(T))}\), which proves
\eqref{eq:p-sobolev-poincare}.
\end{proof}

For \(U\in V_h\) and a coarse element \(T\), let
\(\Lambda_T(U)=(|T|^{-1}\int_T|\nabla U|^p\,dx)^{(p-2)/p}\) be the local
\(p\)-energy density, and \(\Lambda_A(U)\) its maximum over \(T\in A\).
It is not scale invariant.

\begin{lem}[Crude bounds]
\label{lem:p-crude}
Let \(p>2\). For every \(U\in V_h\),
\begin{align}
 \mathcal R_z(U)&\le C\kappa_a\Big(1+\max_{T\subset\operatorname{star}z}\Lambda_T(U)\Big),\qquad \mathcal P_T^e\le\kappa_aC_{P,p}^p,\qquad \mathcal R_z^e\le\kappa_a \label{eq:p-crude-ratio}\\
 \Gamma_\Pi(U)&\le C\kappa_a\big(1+\Lambda_\Omega(U)\big) \label{eq:p-crude-projection}
\end{align}
and, under the restriction on \(p\) in Lemma~\ref{lem:p-lq},
\begin{equation}
 \mathcal P_T(U)\le C\kappa_a\big(1+\Lambda_T(U)\big) \label{eq:p-crude-poincare}
\end{equation}
with constants depending only on \(p\), \(d\) and the mesh family.
Consequently \(\Gamma_\Omega(U)\le C\kappa_a^2(1+\Lambda_\Omega(U))^2\).
\end{lem}

\begin{proof}
Write \(W_U=(1+|\nabla U|)^{p-2}\ge1\). Then
\(\mathcal R_z\le\kappa_a\,{\rm avg}_{\operatorname{star}z}W_U\), and
\(W_U\le c_2(1+|\nabla U|^{p-2})\) and Jensen's inequality for the
concave power \((p-2)/p\) give the first bound. The bounds for
\(\mathcal P_T^e\) and \(\mathcal R_z^e\) follow from
\eqref{eq:p-lq-kernel} with \(k=p\) and from \(a_-\le a_\varepsilon\le a_+\).
Since \(\nabla\Pi_Hv\) is constant on \(T\),
\(\int_TW_U|\nabla\Pi_Hv|^2=({\rm avg}_TW_U)\|\nabla\Pi_Hv\|_T^2\le C(1+\Lambda_T)\|\nabla v\|_{N(T)}^2\)
by \eqref{eq:local-approximation}, and \(W_U\ge1\) gives the weighted
part of \eqref{eq:p-crude-projection}. The coefficient part
\eqref{eq:p-projection-coefficient} follows from \eqref{eq:p-lq-stability}
with \(k\in\{2,p\}\) and \(a_-\le a_\varepsilon\le a_+\). For
\eqref{eq:p-crude-poincare}, split
\(W_U\) in the same way. The constant part is bounded by
\eqref{eq:kernel-poincare}. For the other part, H\"older's inequality and
\eqref{eq:p-sobolev-poincare} give
\begin{equation}
 H^{-2}\int_T|\nabla U|^{p-2}v_f^2\,dx\le\Big(\int_T|\nabla U|^p\,dx\Big)^{(p-2)/p}\big(H^{-1}\|v_f\|_{L^p(T)}\big)^2\le C\Lambda_T(U)\|\nabla v_f\|_{N(T)}^2 \label{eq:p-crude-step}
\end{equation}
since \(|T|^{(p-2)/p}H^{2d/p-d}\) is bounded. Finally
\(\|\nabla v_f\|_{N(T)}^2\le a_-^{-1}\int_{N(T)}\widehat W_U|\nabla v_f|^2\).
\end{proof}

At the fine solution, \(\Lambda_\Omega(u_h)\le(CH^{-d}a_-^{-1}\int_\Omega a_\varepsilon|\nabla u_h|^p)^{(p-2)/p}\),
which \eqref{eq:p-data-large} bounds by the data. Under the restriction on
\(p\) in Lemma~\ref{lem:p-lq}, all constants of Theorem~\ref{thm:p-energy}
are therefore bounded independently of the fine mesh. The resulting rate
satisfies
\(b\gtrsim(\hat\chi_p^{1/2}\kappa_a(1+\Lambda_\Omega(u_h)))^{-1}\). It
depends on the contrast and on the load. Since \(\Lambda\) can grow like a
negative power of \(H\) when the energy concentrates in one coarse
element, this crude bound alone does not guarantee logarithmic
oversampling. The scale-invariant quantities themselves need not grow
with the load. The remainder of this section, and the computations in
Section~\ref{subsec:p-numerics}, concern their size.

\subsection{Refined Projection Constants}
\label{subsec:p-refined}

The ratio \(\mathcal R_z\) arises from the Cauchy--Schwarz inequality
\eqref{eq:p-coefficient-background}, which compares an unweighted moment
with a weighted energy. Two observations weaken the resulting
\(A_2\)-type condition. First, the projection coefficients do not change
when a constant is subtracted from the kernel function. Second, a
weighted Sobolev--Poincar\'e inequality with an exponent \(\sigma>2\) allows
H\"older's inequality with a weaker integrability requirement on the
inverse weight.

\begin{lem}[Moment shift]
\label{lem:p-moment-shift}
Let \(z\) be an interior coarse vertex, let \(\vartheta\) be a coarse
function that is affine on \(K_z\), and let \(v\in V_h\) satisfy
\(\ell_z(v)=0\). Then
\begin{equation}
 \ell_z\big(I_h(\vartheta v)\big)=\ell_z\big(I_h((\vartheta-\vartheta(z))(v-c))\big) \quad \forall c\in\mathbb R \label{eq:p-moment-shift}
\end{equation}
\end{lem}

\begin{proof}
The term \(\vartheta(z)v\) has moment \(\vartheta(z)\ell_z(v)=0\). For a
constant \(c\), \(I_h((\vartheta-\vartheta(z))c)=c(\vartheta-\vartheta(z))\)
is affine on \(K_z\), and the local duality \eqref{eq:local-duality} gives
its moment \(c(\vartheta(z)-\vartheta(z))=0\).
\end{proof}

The moment-shift lemma (Lemma~\ref{lem:p-moment-shift}) applies with
\(\vartheta=\eta_H\) in the weighted cutoff (Lemma~\ref{lem:p-cutoff}) and
with \(\vartheta=\phi_i\) in the weighted assembly
(Lemma~\ref{lem:p-assembly}). In both cases
\(|\vartheta-\vartheta(z)|\le1\) on \(K_z\).

\begin{definition}[Sobolev constant and \texorpdfstring{$A_\sigma$}{A-sigma}-type ratio]
\label{def:p-sobolev-ratio}
Let \(\sigma\ge2\), let \(T\) be a coarse element, \(\omega=N(T)\), and let
\(W\) be positive and constant on fine elements, with
\(W(\omega)=\int_\omega W\,dx\). Let \(C_{S,T}(W)\) be the smallest
constant such that every \(v\in V_h\) admits a constant \(c\) with
\begin{equation}
 \Big(\frac1{W(\omega)}\int_\omega W|v-c|^\sigma\,dx\Big)^{1/\sigma}\le C_{S,T}(W)H\Big(\frac1{W(\omega)}\int_\omega W|\nabla v|^2\,dx\Big)^{1/2} \label{eq:p-sobolev-gain}
\end{equation}
where \(c=0\) is required if \(T\) has a boundary vertex. For a coarse
element \(K\subset\omega\), define
\begin{equation}
 R_\sigma(K,\omega;W)=\Big(\frac1{|\omega|}\int_\omega W\,dx\Big)\Big(\frac1{|K|}\int_KW^{-1/(\sigma-1)}\,dx\Big)^{\sigma-1} \label{eq:p-aq-ratio}
\end{equation}
For a region \(A\), \(C_{S,A}\) and \(R_{\sigma,A}\) denote the maxima over
\(T\in A\) and over the pairs \((K_z,\omega)\) with \(K_z\in A\) and
\(\omega\) equal to \(N(K_z)\) or to \(N(T)\) for a coarse element
\(T\in A\) with vertex \(z\). A maximum over an empty set is zero.
\end{definition}

For \(\sigma=2\), \eqref{eq:p-sobolev-gain} is a weighted Poincar\'e
inequality, and \(R_2(K_z,\operatorname{star}z;W)=\mathcal R_z\). Both
quantities are invariant under \(W\mapsto cW\). By the monotonicity of
power means, \(R_\sigma\) is nonincreasing in \(\sigma\), so a larger exponent in
\eqref{eq:p-sobolev-gain} weakens the requirement on \(W^{-1}\). The
requirement \(c=0\) on elements with a boundary vertex is a weighted
Friedrichs inequality, which uses the zero trace of \(v\in V_h\). It is
needed for the projection stability below and for elements without an
interior vertex, such as two corner elements of the mesh in
Section~\ref{sec:numerical-experiments}.

\begin{lem}[Moment bound]
\label{lem:p-moment-bound}
If \(v\) and \(c\) satisfy \eqref{eq:p-sobolev-gain}, then for every
coarse element \(K\subset\omega\),
\begin{equation}
 H^{-2}W(\omega)\Big(\frac1{|K|}\int_K|v-c|\,dx\Big)^2\le\Big(\frac{|\omega|}{|K|}\Big)^{2/\sigma}C_{S,T}^2R_\sigma(K,\omega;W)^{2/\sigma}\int_\omega W|\nabla v|^2\,dx \label{eq:p-moment-bound}
\end{equation}
\end{lem}

\begin{proof}
H\"older's inequality with exponents \(\sigma\) and \(\sigma/(\sigma-1)\) gives
\(\int_K|v-c|\le(\int_KW|v-c|^\sigma)^{1/\sigma}(\int_KW^{-1/(\sigma-1)})^{(\sigma-1)/\sigma}\).
Bound the first factor by the integral over \(\omega\) and apply
\eqref{eq:p-sobolev-gain}. Squaring and collecting the powers of
\(W(\omega)=|\omega|\,{\rm avg}_\omega W\) and of \(|K|\) gives
\eqref{eq:p-moment-bound}.
\end{proof}

\begin{prop}[Refined cutoff and assembly constants]
\label{prop:p-refined}
Let \(\sigma\ge2\). Lemmas~\ref{lem:p-cutoff} and~\ref{lem:p-assembly},
Proposition~\ref{prop:p-contraction} and Theorems~\ref{thm:p-energy}
and~\ref{thm:p-computed} hold with \(\Gamma_A(U)\) replaced by
\begin{equation}
 \Gamma_A^{(\sigma)}(U)=1+C_{S,A}(\widehat W_U)^2\big(1+R_{\sigma,A}(\widehat W_U)^{2/\sigma}\big)+\Gamma_A^e \label{eq:p-refined-gamma}
\end{equation}
with constants depending additionally on \(\sigma\). Moreover,
\begin{equation}
 \Gamma_\Pi(U)\le\max\big\{CC_{S,\Omega}(\widehat W_U)^2R_{\sigma,\Omega}(\widehat W_U)^{2/\sigma},\Gamma_\Pi^a\big\} \label{eq:p-refined-projection}
\end{equation}
where \(\Gamma_\Pi^a\ge1\) is the smallest constant in
\eqref{eq:p-projection-coefficient} and depends only on \(a_\varepsilon\).
\end{prop}

\begin{proof}
It suffices to bound the two places where \(\mathcal P\) and
\(\mathcal R\) enter: the mass term \(H^{-2}\int_T\widehat W_Ue^2\) and the
background part of the projection correction.

\textbf{Mass term.} Let \(\omega=N(T)\) and let \(c\) be the constant in
\eqref{eq:p-sobolev-gain} for \(v=e\in V_f\). Then
\(H^{-2}\int_T\widehat W_Ue^2\le2H^{-2}\int_\omega\widehat W_U|e-c|^2+2H^{-2}\widehat W_U(T)c^2\).
H\"older's inequality for the probability measure
\(\widehat W_U/\widehat W_U(\omega)\) and \eqref{eq:p-sobolev-gain} bound the
first term by \(2C_{S,T}^2\int_\omega\widehat W_U|\nabla e|^2\). If \(T\)
has a boundary vertex, \(c=0\). Otherwise, let \(z\) be a vertex of
\(T\). Then \(K_z\subset\omega\) and \(\ell_z(e)=0\), so
\(c=-\ell_z(e-c)\) and \(|c|\le C|K_z|^{-1}\int_{K_z}|e-c|\).
Lemma~\ref{lem:p-moment-bound} bounds the second term by
\(CC_{S,T}^2R_\sigma^{2/\sigma}\int_\omega\widehat W_U|\nabla e|^2\).

\textbf{Projection correction.} For \(K=K_z\in\mathcal T_\eta\), apply
Lemma~\ref{lem:p-moment-shift} with \(\vartheta=\eta_H\) and the constant
\(c\) of \eqref{eq:p-sobolev-gain} on \(\omega=N(K)\). Then
\(|c_z|\le C|K|^{-1}\int_K|e-c|\), and
\(|c_z|^2H^{-2}\widehat W_U(\operatorname{star}z)\le C|c_z|^2H^{-2}\widehat W_U(\omega)\)
is bounded by \(CC_S^2R_\sigma^{2/\sigma}\int_\omega\widehat W_U|\nabla e|^2\) by
Lemma~\ref{lem:p-moment-bound}. The error part of the weight is treated
as in Lemma~\ref{lem:p-cutoff}. The weighted assembly lemma (Lemma~\ref{lem:p-assembly}) is modified in
the same way with \(\vartheta=\phi_i\).

\textbf{Projection stability.} Let \(v\in V_h\), \(\omega=N(T)\), and let
\(c\) be the constant in \eqref{eq:p-sobolev-gain} for \(v\). If every vertex
of \(T\) is interior, \eqref{eq:local-constant-reproduction} gives
\(\nabla\Pi_Hv=\sum_z\ell_z(v-c)\nabla\lambda_z\) on \(T\), with the sum over
the vertices of \(T\). If \(T\) has a boundary vertex, then \(c=0\) and the
same formula holds with the sum over the interior vertices of \(T\). In
both cases \(|\ell_z(v-c)|\le C|K_z|^{-1}\int_{K_z}|v-c|\) with
\(K_z\subset\omega\). Since \(\nabla\Pi_Hv\) is constant on \(T\),
\begin{equation}
 \int_T\widehat W_U|\nabla\Pi_Hv|^2\,dx\le CH^{-2}\widehat W_U(\omega)\sum_z\Big(\frac1{|K_z|}\int_{K_z}|v-c|\,dx\Big)^2\le CC_{S,T}^2R_{\sigma,\Omega}^{2/\sigma}\int_\omega\widehat W_U|\nabla v|^2\,dx \label{eq:p-projection-moment}
\end{equation}
by Lemma~\ref{lem:p-moment-bound}. Summation with bounded overlap proves
\eqref{eq:p-refined-projection}.
\end{proof}

Compared with \(\Gamma_A\), the ratio still multiplies a Poincar\'e-type
constant. The gain is the power \(2/\sigma\le1\) and, above all, the
\(A_\sigma\)-type condition with \(\sigma>2\) in place of the \(A_2\)-type
condition.

\subsection{Load and Shape}
\label{subsec:p-load-shape}

\begin{lem}[Scaling and subadditivity]
\label{lem:p-ratio-structure}
For positive weights \(W,W_1,W_2\) constant on fine elements, \(c>0\), and
a pair \(K\subset\omega\),
\begin{align}
 R_\sigma(K,\omega;cW)&=R_\sigma(K,\omega;W) \label{eq:p-ratio-scaling}\\
 R_\sigma(K,\omega;W_1+W_2)&\le R_\sigma(K,\omega;W_1)+R_\sigma(K,\omega;W_2) \label{eq:p-ratio-subadditive}
\end{align}
and \(R_\sigma(K,\omega;W)\le(c''/c')R_\sigma(K,\omega;W')\) whenever
\(c'W'\le W\le c''W'\).
\end{lem}

\begin{proof}
Scaling is immediate. With \(\beta=1/(\sigma-1)\),
\((W_1+W_2)^{-\beta}\le W_i^{-\beta}\) for \(i=1,2\). Splitting the
average of \(W_1+W_2\) over \(\omega\) and bounding the second factor
for each part by the one of \(W_i\) proves \eqref{eq:p-ratio-subadditive}.
The comparison follows from monotonicity of both factors.
\end{proof}

\begin{cor}[Load and shape]
\label{cor:p-load-shape}
For \(U\in V_h\) let \(G_U=a_\varepsilon|\nabla U|^{p-2}\). Then
\begin{equation}
 R_\sigma(K,\omega;\widehat W_U)\le c_2\big(R_\sigma(K,\omega;a_\varepsilon)+R_\sigma(K,\omega;G_U)\big) \label{eq:p-load-shape}
\end{equation}
The first term depends only on the coefficient. The second is invariant
under \(U\mapsto sU\) for \(s>0\), so it depends on the shape of the
gradient magnitude and not on its size.
\end{cor}

\begin{proof}
We have \(\widehat W_U\le c_2(a_\varepsilon+G_U)\) and
\(\widehat W_U^{-\beta}\le\min\{a_\varepsilon^{-\beta},G_U^{-\beta}\}\), since
\((1+t)^{p-2}\ge\max\{1,t^{p-2}\}\). Splitting the first factor of
\eqref{eq:p-aq-ratio} and bounding the second factor by the matching
term of the minimum gives \eqref{eq:p-load-shape}.
\end{proof}

At small loads, \(\widehat W_U\) is close to \(a_\varepsilon\) and
Lemma~\ref{lem:p-crude} applies. At large loads, only the shape of
\(a_\varepsilon|\nabla u_h|^{p-2}\) matters. The Sobolev constant
\(C_{S,T}(\widehat W_U)\) is also invariant under scaling of the weight, but
it is not subadditive, so it remains a quantity of \(\widehat W_U\) itself.

\subsection{Flux Balance and Model Weights}
\label{subsec:p-model-weights}

The ratios are large where \(|\nabla u_h|\) is small compared with its
coarse-scale average. Near a nondegenerate critical point the weight
behaves like \((1+s|x-x_0|^2)^{(p-2)/2}\) with \(s\) proportional to the
squared Hessian, which would make the \(A_2\)-type ratio grow with the
load when \(p-2\ge d\), logarithmically when \(p-2=d\). A balance of
fluxes shows that this
is not the behaviour where the load is positive.

\begin{prop}[Flux balance]
\label{prop:p-flux-balance}
Let \(B_\rho=B(x_0,\rho)\subset\Omega\), let \(f\ge f_0>0\) almost
everywhere in \(B_\rho\), and let \(u\in W^{1,\infty}(B_\rho)\) satisfy
\(\int a(x,\nabla u)\cdot\nabla\varphi\,dx=\int f\varphi\,dx\) for all
\(\varphi\in W_0^{1,\infty}(B_\rho)\). Then, for every \(r\in(0,\rho)\),
\begin{equation}
 \operatorname*{ess\,sup}_{B_r}(1+|\nabla u|)^{p-2}|\nabla u|\ge\frac{f_0r}{d\hat Ma_+} \label{eq:p-flux-balance}
\end{equation}
In particular, \(\operatorname*{ess\,sup}_{B_r}|\nabla u|\ge(f_0r/(2^{p-2}d\hat Ma_+))^{1/(p-1)}\)
whenever the left side of this inequality is at least one.
\end{prop}

\begin{proof}
Test with \(\varphi_k(x)=\min\{1,k(r-|x-x_0|)_+\}\). Its gradient is
\(-k(x-x_0)/|x-x_0|\) on the annulus \(r-1/k<|x-x_0|<r\) and zero
elsewhere. Since \(|a(x,\xi)|\le\hat Ma_+(1+|\xi|)^{p-2}|\xi|\), the left
side is at most \(k|B_r\setminus B_{r-1/k}|\) times the essential
supremum in \eqref{eq:p-flux-balance}, and the right side is at least
\(f_0|B_{r-1/k}|\). Letting \(k\to\infty\) and using
\(|\partial B_r|=d|B_r|/r\) proves \eqref{eq:p-flux-balance}. If the
supremum \(X\) of \(|\nabla u|\) is at least one, then
\((1+X)^{p-2}X\le2^{p-2}X^{p-1}\).
\end{proof}

Consider loads \(sf\) with \(s\to\infty\) and states
\(u_s=s^{1/(p-1)}\widetilde u\), the scaling of the large-load limit of
\eqref{eq:p-structure}. If \(\widetilde u\) had a critical point with
\(|\nabla\widetilde u|\le\gamma|x-x_0|\) near \(x_0\) where \(f\ge f_0\),
then \eqref{eq:p-flux-balance} as \(r\to0\) would require
\(s^{1/(p-1)}\gamma\ge sf_0/(d\hat Ma_+)\), which fails for large \(s\).
Instead, \eqref{eq:p-flux-balance} allows the profile
\(|\nabla u_s|\sim(sr)^{1/(p-1)}\) outside a core where \(|\nabla u_s|\le1\).
The weight then behaves like \(1+(s|x-x_0|)^{(p-2)/(p-1)}\). This is a
heuristic description, not a regularity result. The following lemma
treats both profiles.

\begin{lem}[Model weights]
\label{lem:p-model-weights}
Let \(\alpha\ge0\), \(s>0\) and \(x_0\in\mathbb R^d\), and let \(W\) be
constant on fine elements with
\begin{equation}
 c'g_-\big(1+(sd_t^2)^{\alpha/2}\big)\le W|_t\le c''g_+\big(1+(sd_t^2)^{\alpha/2}\big),\qquad d_t=\operatorname{dist}(x_0,t)+h \label{eq:p-model-weight}
\end{equation}
on every fine element \(t\), where \(h\) is the largest fine element
diameter. If \(\alpha<d(\sigma-1)\), then
\(R_\sigma(K,\omega;W)\le C\) for every pair \(K\subset\omega\), with \(C\)
depending only on \(\sigma\), \(\alpha\), \(d\), \(c''g_+/(c'g_-)\) and the
mesh family, and not on \(s\), \(h\) or \(x_0\). If \(\sigma=2\),
\(\alpha=d\), \(x_0\in\overline K\) and \(sH^2\ge1\), then
\(c\log\min\{sH^2,H^2/h^2\}-C\le R_2(K,\omega;W)\le C(1+\log(sH^2))\).
\end{lem}

\begin{proof}
Let \(D\) be the diameter of \(\omega\), \(\rho\) the distance from
\(x_0\) to \(\omega\), and \(\beta=1/(\sigma-1)\). For \(x\in t\),
\(|x-x_0|\le d_t\le|x-x_0|+h\). If \(\rho\ge D\), then
\(d_t\in[\rho,4\rho]\) on \(\omega\), so \(W\) varies by a bounded factor
there and the ratio is bounded. If \(\rho<D\), then
\(\omega\subset B=B(x_0,2D)\), \({\rm avg}_\omega W\le c''g_+(1+(9sD^2)^{\alpha/2})\),
and \({\rm avg}_KW^{-\beta}\le(c'g_-)^{-\beta}(|B|/|K|){\rm avg}_B(1+(s|x-x_0|^2)^{\alpha/2})^{-\beta}\).
With \(R=2D\) and \(sR^2>1\), split the radial integral at
\(r=s^{-1/2}\) and use \(\alpha\beta<d\):
\begin{equation}
 \frac d{R^d}\int_0^Rr^{d-1}\big(1+(sr^2)^{\alpha/2}\big)^{-\beta}\,dr\le(sR^2)^{-d/2}+\frac d{d-\alpha\beta}(sR^2)^{-\alpha\beta/2}\le C\big(1+(sR^2)^{\alpha/2}\big)^{-\beta} \label{eq:p-model-radial}
\end{equation}
The same bound holds trivially for \(sR^2\le1\). Raising to the power
\(\sigma-1\) and multiplying by the bound for \({\rm avg}_\omega W\) bounds
\(R_\sigma\). For \(\sigma=2\) and \(\alpha=d\), the same radial integral, now
logarithmic, gives the upper bound. For the lower bound, shape
regularity gives \(|\{x\in K:|x-x_0|<r\}|\ge cr^d\) for \(r\le cH\). For
\(|x-x_0|\ge r_0=\max\{s^{-1/2},h\}\), \(d_t\le2|x-x_0|\). The integral over
\(r_0<|x-x_0|<cH\) gives
\({\rm avg}_KW^{-1}\ge c(sH^2)^{-d/2}\log(cH/r_0)\). Since
\(d_t\ge|x-x_0|\) and a fixed fraction of \(\omega\) lies at distance at
least \(cH\) from \(x_0\), \({\rm avg}_\omega W\ge c(sH^2)^{d/2}\).
\end{proof}

For the quartic law in two dimensions, a nondegenerate critical point
gives \(\alpha=2=d\). The \(A_2\)-type ratio is then borderline and grows
logarithmically, whereas \(R_\sigma\) with \(\sigma>2\) stays bounded. The
flux-balanced profile gives \(\alpha=(p-2)/(p-1)<1\), so every \(R_\sigma\),
including \(R_2\), is bounded uniformly in the load. The hypothesis
\eqref{eq:p-model-weight} excludes gradients that vanish on whole fine
elements; if \(\nabla U=0\) on some element, \(R_\sigma(W_U)\) can grow without
bound. For \(\sigma=2\) and \(\alpha<d\), Corollary~\ref{cor:p-model-a2} below
proves the inequality \eqref{eq:p-sobolev-gain} for these model weights
uniformly in \(s\) and \(h\). For \(\sigma>2\), which is needed in the borderline
case \(\alpha=d\), we have not proved it.

\subsection{Load-Uniform Constants}
\label{subsec:p-uniform}

The quantities in \(\Gamma^{(\sigma)}\) are invariant under scaling of the
weight. We give two sufficient conditions under which they are bounded
uniformly in the load. The first is a local Muckenhoupt condition on the
shape of the gradient magnitude. It gives bounds that are also independent
of the fine mesh, without the restriction on \(p\) in
Lemma~\ref{lem:p-lq}. The second requires only that the large-load limit
on a fixed mesh is nondegenerate.

For a weight \(w\), that is, a positive locally integrable function on
\(\mathbb R^d\), let
\begin{equation}
 [w]_{A_2}=\sup_B\Big(\frac1{|B|}\int_Bw\,dx\Big)\Big(\frac1{|B|}\int_Bw^{-1}\,dx\Big) \label{eq:p-a2-constant}
\end{equation}
where the supremum is over all balls \(B\subset\mathbb R^d\). The argument of
Lemma~\ref{lem:p-ratio-structure} shows that
\([w_1+w_2]_{A_2}\le[w_1]_{A_2}+[w_2]_{A_2}\), \([cw]_{A_2}=[w]_{A_2}\), and
\([w]_{A_2}\le(c''/c')[w']_{A_2}\) if \(c'w'\le w\le c''w'\). For
\(E\subset B\), the Cauchy--Schwarz inequality gives
\begin{equation}
 \Big(\frac{|E|}{|B|}\Big)^2\le[w]_{A_2}\frac{w(E)}{w(B)} \label{eq:p-a2-doubling}
\end{equation}
with \(w(E)=\int_Ew\,dx\). In particular \(w\) is doubling:
\(w(B_{2r})\le4^d[w]_{A_2}w(B_r)\) for concentric balls.

Let \(\mathcal T_h^\partial\) be the set of fine elements whose vertices all
lie on \(\partial\Omega\). Every \(v\in V_h\) vanishes on these elements, and
so does \(\nabla U\) for \(U\in V_h\). A function \(W\ge0\) on \(\Omega\),
constant on fine elements and positive outside \(\mathcal T_h^\partial\),
satisfies the \emph{local \(A_2\) condition with constant \(A\)} if, for every
coarse element \(T\), there is a weight \(w_T\) on \(\mathbb R^d\) with
\([w_T]_{A_2}\le A\), \(w_T=W\) on the fine elements of \(N(T)\) not in
\(\mathcal T_h^\partial\), and \(w_T\ge W\) on those in \(\mathcal T_h^\partial\).

Near the boundary we use the following property (B), with constants
\(C_G\) and \(H_0\). If \(H\le H_0\), then for every coarse element \(T\)
with a boundary vertex there is an open set
\(E_T\subset\mathbb R^d\setminus\overline\Omega\), contained in a ball of radius
\(C_GH\) around \(T\), with \(|E_T|\ge C_G^{-1}H^d\), such that the interior
\(\omega_T^+\) of \(\overline{N(T)\cup E_T}\) is a John domain with
constant \(C_G\). For a polyhedral Lipschitz domain, the open half ball
outside \(\Omega\) centred at the barycentre of a boundary facet of
\(N(T)\), with radius \(c_\Omega H\), has this property. Here \(c_\Omega\)
and \(H_0\) depend on the dihedral angles of \(\Omega\) and on the
separation of its nonadjacent faces.

\begin{lem}[Weighted patch inequalities for Muckenhoupt weights]
\label{lem:p-a2-poincare}
Assume property (B), and that no element of \(\mathcal T_h^\partial\) lies in
a selected simplex \(K_z\). Let \(W\) satisfy the local \(A_2\) condition
with constant \(A\). Then, for \(\sigma=2\), \(C_{S,T}(W)\le C_A\) for every
coarse element \(T\). Moreover, \(R_2(K,\omega;W)\le C_A\) for every coarse
element \(K\subset\omega=N(T)\) that contains no element of
\(\mathcal T_h^\partial\), in particular for \(K=K_z\). The constant \(C_A\)
depends only on \(d\), \(A\), \(C_G\) and the mesh family.
\end{lem}

\begin{proof}
Let \(\omega=N(T)\) and \(w=w_T\). Since \(v\) and \(\nabla v\) vanish on
\(\mathcal T_h^\partial\) for \(v\in V_h\), and \(W\le w\) on \(\omega\) with
equality off \(\mathcal T_h^\partial\), we have
\(\int_\omega W|\nabla v|^2=\int_\omega w|\nabla v|^2\) and
\(\int_\omega W|v-c|^2\le\int_\omega w|v-c|^2\) for every \(c\). For \(\sigma=2\)
the normalization in \eqref{eq:p-sobolev-gain} cancels, so it suffices to
prove that inequality with \(w\) in place of \(W\). Let \(B\) be a ball of
radius \(CH\) containing \(\omega\) and, if \(T\) has a boundary vertex,
\(\omega_T^+\). If \(K\subset\omega\) contains no element of
\(\mathcal T_h^\partial\), then \(W=w\) on \(K\), and \(|B|\le C|K|\) gives
\({\rm avg}_\omega W\le C\,{\rm avg}_Bw\) and
\({\rm avg}_KW^{-1}\le C\,{\rm avg}_Bw^{-1}\). Hence \(R_2(K,\omega;W)\le CA\).

By shape regularity, each simplex contains a ball of radius \(cH\) around
its incentre, and two face neighbours contain a ball of radius \(cH\)
around the barycentre of their common face. Join a point \(x\in\omega\)
to the incentre of its element by a segment and continue along the
polygonal path through face barycentres and incentres of a face chain to
\(T\), which has bounded length by the proof of
Lemma~\ref{lem:one-layer-patches}. On the first segment, the distance to
\(\partial\omega\) is at least a fixed multiple of the distance travelled,
by convexity of the element. On the polygonal part it is at least a fixed
multiple of \(H\), while the total length is at most \(CH\). Hence
\(\omega\) is a John domain with a constant depending only on \(d\) and
the mesh family.

Muckenhoupt \(A_2\) weights are doubling and satisfy the weighted
Poincar\'e inequality on balls for Lipschitz
functions~\cite{fabes-kenig-serapioni1982}. On John domains, a doubling
measure that supports a Poincar\'e inequality on balls supports it on the
whole domain, by the chaining argument of Haj{\l}asz and
Koskela~\cite[Theorem~9.7]{hajlasz-koskela2000}; the lower mass bound used
there follows from doubling. See also Chua~\cite{chua1993} for weighted
inequalities on domains satisfying the Boman chain condition. Since the
\(A_2\) constant is invariant under dilations, scaling to \(H=1\) gives
\begin{equation}
 \inf_{c\in\mathbb R}\int_Dw|v-c|^2\,dx\le CH^2\int_Dw|\nabla v|^2\,dx \quad \forall v\in W^{1,\infty}(D),\quad D\in\{\omega,\omega_T^+\} \label{eq:p-a2-poincare}
\end{equation}
with \(C\) depending only on \(d\), \(A\) and the John constant. For \(D=\omega\)
this is \eqref{eq:p-sobolev-gain} with \(\sigma=2\) when every vertex of \(T\)
is interior. If \(T\) has a boundary vertex, extend \(v\in V_h\) by zero to
\(\omega_T^+\). The extension is Lipschitz, because \(v\) vanishes on
\(\partial\Omega\), and its gradient vanishes on \(E_T\). Let \(c\) attain
the infimum in \eqref{eq:p-a2-poincare} with \(D=\omega_T^+\). Since \(v=0\)
on \(E_T\), we have \(c^2w(E_T)\le CH^2\int_\omega w|\nabla v|^2\). By
\eqref{eq:p-a2-doubling}, \(w(\omega)\le w(B)\le CAw(E_T)\). Hence
\begin{equation}
 \int_\omega wv^2\,dx\le2\int_{\omega_T^+}w|v-c|^2\,dx+2c^2w(\omega)\le CH^2\int_\omega w|\nabla v|^2\,dx \label{eq:p-a2-friedrichs}
\end{equation}
which is \eqref{eq:p-sobolev-gain} with \(\sigma=2\) and \(c=0\).
\end{proof}

\begin{thm}[Load-uniform constants under a Muckenhoupt condition]
\label{thm:p-a2-uniform}
Assume property (B), and that no element of \(\mathcal T_h^\partial\) lies in a
selected simplex. Let \(U\in V_h\), and assume that for every coarse
element \(T\) there is a weight \(g_T\) on \(\mathbb R^d\) with
\([g_T]_{A_2}\le A_G\) and \(g_T=|\nabla U|^{p-2}\) on the fine elements of
\(N(T)\) not in \(\mathcal T_h^\partial\). Then \(\widehat W_U\) satisfies the
local \(A_2\) condition with constant \(A\le2c_2\kappa_a(1+A_G)\), and
\begin{equation}
 \Gamma_\Omega^{(2)}(U)\le C_A'+\Gamma_\Omega^e,\qquad \Gamma_\Pi(U)\le\max\{C_A',\Gamma_\Pi^a\},\qquad \Gamma_\Omega^e\le C\kappa_a^2,\qquad \Gamma_\Pi^a\le C\kappa_a \label{eq:p-a2-uniform}
\end{equation}
where \(C_A'\) depends only on \(d\), \(A\), \(C_G\) and the mesh family,
and \(C\) only on \(p\), \(d\) and the mesh family. In particular, if this
holds for \(U=u_h\) with \(A_G\) independent of the load, then
Theorem~\ref{thm:p-energy}, with \(\Gamma_\Omega\) replaced by
\(\Gamma_\Omega^{(2)}\) as in Proposition~\ref{prop:p-refined}, has a
contraction factor and a prefactor bounded independently of the load and
of the fine mesh.
\end{thm}

\begin{proof}
On the fine elements of \(N(T)\) not in \(\mathcal T_h^\partial\),
\((1+g_T)/2\le\max\{1,g_T\}\le(1+|\nabla U|)^{p-2}\le c_2(1+g_T)\). Set
\(w_T=\widehat W_U\) on these elements and \(w_T=a_+(1+g_T)\) elsewhere, in
particular on \(\mathcal T_h^\partial\), where \(\widehat W_U=a_\varepsilon\le w_T\).
Then \(a_-(1+g_T)/2\le w_T\le c_2a_+(1+g_T)\) on \(\mathbb R^d\), so
\([w_T]_{A_2}\le2c_2\kappa_a[1+g_T]_{A_2}\le2c_2\kappa_a(1+A_G)\).
Lemma~\ref{lem:p-a2-poincare} and Proposition~\ref{prop:p-refined} with
\(\sigma=2\) give the first two bounds. Lemma~\ref{lem:p-lq} gives
\(\mathcal P_T^e\le\kappa_aC_{P,p}^p\), \(\mathcal R_z^e\le\kappa_a\) and
\(\Gamma_\Pi^a\le\kappa_a\max\{C_{\Pi,2}^2,C_{\Pi,p}^p\}\). None of these
bounds uses Lemma~\ref{lem:p-crude} or the Sobolev--Poincar\'e inequality
\eqref{eq:p-sobolev-poincare}.
\end{proof}

The condition on \(|\nabla u_h|^{p-2}\) is invariant under scaling of
\(u_h\). It concerns only the shape of the gradient magnitude, in
agreement with Corollary~\ref{cor:p-load-shape}. For the model profiles of
Section~\ref{subsec:p-model-weights} it can be verified.

\begin{cor}[Model weights below the critical exponent]
\label{cor:p-model-a2}
If \(W\) satisfies \eqref{eq:p-model-weight} with \(0\le\alpha<d\), then
\(W\) satisfies the local \(A_2\) condition with a constant depending only
on \(\alpha\), \(d\) and \(c''g_+/(c'g_-)\), and not on \(s\), \(h\) or
\(x_0\). Consequently, under property~(B), \(C_{S,T}(W)\) and
\(R_2(K,\omega;W)\) are bounded uniformly in \(s\) and \(h\), by a constant
that depends in addition on \(C_G\) and the mesh family.
\end{cor}

\begin{proof}
Let \(w(x)=1+s^{\alpha/2}(|x-x_0|+h)^\alpha\) on \(\mathbb R^d\). For \(x\) in
a fine element \(t\), \(|x-x_0|\le d_t\le|x-x_0|+h\le2d_t\), so
\(2^{-\alpha}c'g_-w\le W\le c''g_+w\) on \(\Omega\). Extend \(W\) by
\(c''g_+w\) outside \(\Omega\); the extension lies between the same
multiples of \(w\). Since \(\max\{r,h\}\le r+h\le2\max\{r,h\}\), the weight
\(w\) lies between fixed multiples, depending only on \(\alpha\), of
\(c_0+c_1|x-x_0|^\alpha\) with \(c_0=1+s^{\alpha/2}h^\alpha\) and
\(c_1=s^{\alpha/2}\). By subadditivity and invariance under scaling and
translation, \([c_0+c_1|x-x_0|^\alpha]_{A_2}\le1+[|x|^\alpha]_{A_2}\), which is
finite for \(0\le\alpha<d\), see~\cite{grafakos2014}. Lemma~\ref{lem:p-a2-poincare}
completes the proof. Since \(w_T=W\) on every fine element, its bound for
\(R_2\) holds for every coarse element \(K\subset\omega\), and its condition
on the selected simplices is not needed.
\end{proof}

The flux-balanced profile of Section~\ref{subsec:p-model-weights} has
\(\alpha=(p-2)/(p-1)<1\le d\), so it is covered for every \(p\) and \(d\).
A nondegenerate critical point has \(\alpha=p-2\). This exponent is
covered when \(p-2<d\), but not for the quartic law in two dimensions,
where \(R_2\) grows logarithmically by Lemma~\ref{lem:p-model-weights}.

On a fixed pair of meshes, a weaker hypothesis suffices. We say that the
law is asymptotically \(p\)-homogeneous if there is a potential
\(\psi_\infty(x,\xi)\), strictly convex in \(\xi\) and with
\(\psi_\infty(x,t\xi)=t^p\psi_\infty(x,\xi)\) for \(t>0\), such that
\(t^{-p}\psi(x,t\xi)\to\psi_\infty(x,\xi)\) as \(t\to\infty\), uniformly in
\(x\) and in \(\xi\) on compact sets. The quartic law has
\(\psi_\infty=a_\varepsilon|\xi|^4/4\).

\begin{thm}[Load uniformity on a fixed mesh]
\label{thm:p-fixed-mesh}
Let the law be asymptotically \(p\)-homogeneous, fix \(F\in V_h'\setminus\{0\}\),
and let \(u_s\) be the fine solution for the load \(sF\), \(s>0\). Let
\(\widetilde u_h\) minimize \(\int_\Omega\psi_\infty(x,\nabla v)\,dx-F(v)\) over
\(V_h\). Assume that \(\nabla\widetilde u_h\ne0\) on every fine element not in
\(\mathcal T_h^\partial\), and that no element of \(\mathcal T_h^\partial\) lies in a
selected simplex \(K_z\). Then
\begin{equation}
 \sup_{s>0}\big(\Gamma_\Omega(u_s)+\Gamma_\Pi(u_s)\big)<\infty \label{eq:p-fixed-mesh}
\end{equation}
so the contraction factor and the prefactor in Theorem~\ref{thm:p-energy}
are bounded uniformly in the load on these meshes.
\end{thm}

\begin{proof}
Let \(\tau_s=s^{1/(p-1)}\) and \(w_s=u_s/\tau_s\). Since
\(s\tau_s=\tau_s^p\), the function \(w_s\) minimizes
\(J_s(w)=\int_\Omega\tau_s^{-p}\psi(x,\tau_s\nabla w)\,dx-F(w)\) over \(V_h\).
By \eqref{eq:p-potential-lower},
\(\tau^{-p}\psi(x,\tau\xi)\ge c_\psi\hat ma_-|\xi|^p\), so the minimizers
\(w_s\), \(s\ge1\), are bounded in the finite-dimensional space \(V_h\). On
bounded sets of \(V_h\), \(J_s\) converges uniformly to
\(J_\infty(w)=\int_\Omega\psi_\infty(x,\nabla w)\,dx-F(w)\). Every limit point
\(w\) of \(w_s\) as \(s\to\infty\) therefore satisfies
\(J_\infty(w)=\lim J_s(w_s)\le\lim J_s(v)=J_\infty(v)\) for all \(v\in V_h\),
and uniqueness gives \(w_s\to\widetilde u_h\).

The quantities \(\mathcal P_T\), \(\mathcal R_z\) and \(\Gamma_\Pi\) are
invariant under scaling of the weight. They are continuous functions of
the nonnegative elementwise weight values on the set where these values
are positive outside \(\mathcal T_h^\partial\). Indeed, \(\mathcal P_T\) is the
largest eigenvalue of a pencil on the restrictions of \(V_f\) to \(N(T)\),
whose forms depend linearly on the weight values. Its second form
\(\int_{N(T)}\widehat W|\nabla v|^2\) is positive definite there, because
\(\nabla v\) vanishes on \(\mathcal T_h^\partial\): if the form vanishes, \(v\) is
constant on \(N(T)\), and \eqref{eq:kernel-poincare} forces \(v=0\). The
same holds for \(\Gamma_\Pi\) on \(V_h\). In \(\mathcal R_z\), the average of
\(\widehat W^{-1}\) over \(K_z\) involves only elements outside
\(\mathcal T_h^\partial\) by assumption, and the average of \(\widehat W\) is
linear in the weight values.

The map \(s\mapsto u_s\) is continuous. The minimizers are unique, they
are bounded uniformly for \(s\) in compact subsets of \((0,\infty)\) by
\eqref{eq:p-data-p}, and the energies depend continuously on \(s\), so the
compactness argument above applies at every \(s>0\).
Hence \(f(s)=\Gamma_\Omega(u_s)+\Gamma_\Pi(u_s)\) is continuous on
\((0,\infty)\). As \(s\to0\), \(u_s\to0\) by \eqref{eq:p-data-two}, so
\(\widehat W_{u_s}\to a_\varepsilon\). As \(s\to\infty\), scale invariance gives
\(f(s)=\Gamma_\Omega+\Gamma_\Pi\) evaluated at the weight
\(\tau_s^{2-p}\widehat W_{u_s}=a_\varepsilon(\tau_s^{-1}+|\nabla w_s|)^{p-2}\),
which converges to \(a_\varepsilon|\nabla\widetilde u_h|^{p-2}\). This limit is
positive outside \(\mathcal T_h^\partial\). Thus \(f\) has finite limits at
both ends and is bounded.
\end{proof}

The bound depends on the meshes, and the hypothesis can be checked a
posteriori on \(\widetilde u_h\). For the pixel coefficient of
Section~\ref{subsec:p-numerics} with \(n_H=12\) and \(n_h=48\), the rescaled
solution at load \(2\times10^6\) has \(\min|\nabla w_s|/\max|\nabla w_s|=0.080\)
outside \(\mathcal T_h^\partial\), and this ratio changes by less than
\(2\times10^{-4}\) up to load \(2\times10^9\). For \(n_h=96\) and \(192\),
the ratios at load \(2\times10^9\) are \(0.106\) and \(0.041\). The two
corner elements of \(\mathcal T_h^\partial\) lie in no selected simplex.
These computations indicate that the hypothesis holds on the three meshes,
with a margin that depends on the fine mesh. This is consistent with
the saturation in Table~\ref{tab:p-load}.
Theorems~\ref{thm:p-a2-uniform} and~\ref{thm:p-fixed-mesh} bracket the
remaining question: whether \(|\nabla u_h|^{p-2}\) satisfies the local
\(A_2\) condition outside \(\mathcal T_h^\partial\) uniformly in the load and the
mesh.

\subsection{A Priori Compatibility}
\label{subsec:p-apriori}

Theorem~\ref{thm:p-computed} assumes \(\lambda_C<1\) at the computed
coordinate. For coordinates with bounded energy, this follows from
oversampling. For fixed \(H\), the required depth grows logarithmically
with the fine resolution. Through the crude constants, it also grows with
the load and like a power of \(1/H\). Let \(h_{\min}\) be the smallest fine element diameter, so
that \(|t|\ge c_hh_{\min}^d\) for every fine element \(t\), with \(c_h\)
depending only on shape regularity.

\begin{prop}[A priori compatibility]
\label{prop:p-apriori}
Let \(p>2\) satisfy the restriction of Lemma~\ref{lem:p-lq}, let
\(E_0\in\mathbb R\), and let \(v_H\in V_H\) satisfy
\(\mathcal E(M_h^{\mathcal D}(v_H))\le E_0\). Define
\begin{align}
 K_0&=\Big(\frac{\hat L_F+(c_\psi\hat mE_0^+)^{1/2}}{c_\psi\hat m}\Big)^2,\qquad \Lambda_0=\Big(\frac{C_\Lambda K_0}{a_-H^d}\Big)^{(p-2)/p} \label{eq:p-apriori-energy}\\
 \Gamma_0&=C\kappa_a^2(1+\Lambda_0)^2,\qquad \Gamma_{\Pi,0}=C\kappa_a(1+\Lambda_0),\qquad \theta_0=\frac{\kappa_0}{1+\kappa_0},\qquad \kappa_0^2=\hat\chi_pC_{{\rm cut},p}\Gamma_0 \label{eq:p-apriori-constants}\\
 Q_0&=(4c_2+2^p)\hat\chi_p^{(p-2)/2}\Gamma_{\Pi,0}K_0,\qquad Y_0=\frac{C_{\rm amp}\hat\chi_pQ_0}{c_ha_-} \label{eq:p-apriori-amplitude}
\end{align}
with the constants of Lemma~\ref{lem:p-crude} and \(C_\Lambda\) depending
only on the mesh family. Let \(n=\min_in_i\) over the proper patches,
with \(N^{3n_i+r_P}(\omega_i)\subset D_i\), and let \(\mathcal N_t\) be computed with
\(t_i=\theta_0^{n_i}\). Then
\begin{equation}
 \lambda_C^2\le\frac{\theta_0^nY_0}{h_{\min}^d},\qquad \lambda_{\rm rec}^2\le\frac{C_{{\rm asm},p}\Gamma_0\mathcal N_t(\mathcal D)Y_0}{h_{\min}^d},\qquad \varrho_D\le c_2\Big(1+\Big(\frac{CY_0}{H^d}\Big)^{(p-2)/p}\Big) \label{eq:p-apriori}
\end{equation}
In particular, \(\lambda_C\le1/2\) whenever
\(n\ge(d\log(1/h_{\min})+\log(4Y_0))/|\log\theta_0|\).
\end{prop}

\begin{proof}
Let \(U=M_h(v_H)\) and \(\widetilde u=M_h^{\mathcal D}(v_H)\). Since
\(\widetilde u\in v_H+V_f\) and \(U\) minimizes the energy on this affine
space, \(\mathcal E(U)\le E_0\). With
\(X=\int_\Omega a_\varepsilon(|\nabla U|^2+|\nabla U|^p)\),
\eqref{eq:p-potential-lower} and \(F(U)\le\hat L_FX^{1/2}\) give
\(c_\psi\hat mX-\hat L_FX^{1/2}\le E_0\), hence \(X\le K_0\). Since
\(|T|\ge cH^d\), \(\Lambda_\Omega(U)\le\Lambda_0\), and
Lemma~\ref{lem:p-crude} gives \(\Gamma_\Omega(U)\le\Gamma_0\) and
\(\Gamma_\Pi(U)\le\Gamma_{\Pi,0}\). The contraction factor of
\eqref{eq:p-interior-decay} is therefore at most \(\theta_0\), and interior
decay holds with \(\theta_0\).

As in the proof of \eqref{eq:p-amplitude-data}, the correction
\(q=(I-\Pi_H)U\) satisfies \(\mathcal Q_\Omega(q;U)\le Q_0\).
Lemma~\ref{lem:p-amplitude} and the interior decay
\eqref{eq:p-interior-decay} give, with \(e_i=U_i-U\),
\begin{equation}
 \mathscr D_{D_i}(U_i,U)\le C_{\rm amp}\hat M\mathcal Q_{D_i}(q;U),\qquad \mathscr D_{N(\omega_i)}(U_i,U)\le\theta_0^{n_i}C_{\rm amp}\hat MQ_0 \label{eq:p-apriori-decay}
\end{equation}
On a fine element \(t\subset N(\omega_i)\), \(\nabla e_i\) is constant, so
\begin{equation}
 a_-|t||\nabla e_i|^2\le\int_{N(\omega_i)}a_\varepsilon|\nabla e_i|^2\,dx\le\frac{\mathscr D_{N(\omega_i)}(U_i,U)}{c_p\hat m} \label{eq:p-apriori-pointwise}
\end{equation}
by \eqref{eq:p-natural-equivalence} and \(\omega(\nabla U,\nabla U_i)\ge1\).
Since \(1+|\nabla U|\ge1\), this proves the bound for \(\lambda_C\). For
\(\lambda_{\rm rec}\), \(\widetilde u-U=(I-\Pi_H)I_h\sum_i\phi_ie_i\).
Lemma~\ref{lem:p-assembly}, \(\widehat W_U\ge a_\varepsilon\),
\eqref{eq:p-splitting} and \eqref{eq:p-apriori-decay} give
\begin{equation}
 a_-|t||\nabla(\widetilde u-U)|^2\le\int_\Omega\widehat W_U|\nabla(\widetilde u-U)|^2\,dx\le C_{{\rm asm},p}\Gamma_0\sum_i\theta_0^{n_i}C_{\rm amp}\hat\chi_p\mathcal Q_{D_i}(q;U) \label{eq:p-apriori-rec}
\end{equation}
and the definition of the weighted overlap bounds the sum by
\(\mathcal N_tQ_0\). For \(\varrho_D\), let \(T\subset D_i\). Then
\(\widehat W_{U_i}\le c_2(\widehat W_U+a_\varepsilon|\nabla e_i|^{p-2})\), and
H\"older's inequality gives
\begin{equation}
 \int_Ta_\varepsilon|\nabla e_i|^{p-2}\,dx\le\Big(\int_Ta_\varepsilon\,dx\Big)^{2/p}\Big(\int_Ta_\varepsilon|\nabla e_i|^p\,dx\Big)^{(p-2)/p} \label{eq:p-apriori-holder}
\end{equation}
The last factor is at most \((\mathscr D_{D_i}(U_i,U)/(c_p\hat m))^{(p-2)/p}\).
Dividing by \(\int_T\widehat W_U\ge\int_Ta_\varepsilon\ge ca_-H^d\) proves the
bound for \(\varrho_D\).
\end{proof}

\begin{cor}[A priori computed-state error]
\label{cor:p-apriori}
Assume the hypotheses of Theorem~\ref{thm:p-computed} except
\(\lambda_C<1\), let \(p\) and \(E_0\) be as in Proposition~\ref{prop:p-apriori},
and let \(D_i=N^\ell(\omega_i)\) with \(t_i=\theta_0^{n_\ell}\). There is a
depth \(\ell_*\) such that, for every \(\ell\ge\ell_*\), the bound
\eqref{eq:p-computed} holds at every coordinate \(v_H\) with
\(\mathcal E(M_h^\ell(v_H))\le E_0\). In particular, it holds at every
tangent-stationary coordinate in this energy sublevel set, and at every
iterate of an outer descent method started from such a coordinate. The
depth \(\ell_*\) depends on the data, \(E_0\), \(H\), the contrast and the
structural constants. For fixed \(H\) it grows like
\(\log(1/h_{\min})\); through \(\theta_0\) it grows with the load and like a
power of \(1/H\).
\end{cor}

\begin{proof}
For \(n_\ell\) as in Proposition~\ref{prop:p-apriori}, \(\lambda_C\le1/2\)
once \(n_\ell\) exceeds the stated threshold. Then
\eqref{eq:p-apriori}, \(\mathcal N_0\le C_{\rm ol}(1+\ell)^d\) and
\(\mathcal N_t\le C_{\rm ol}(1+\ell)^d\theta_0^{n_\ell}\) bound
\(\lambda_{\rm rec}\), \(\varrho_D\), \(C_{\rm tan}\) and \(K_v\) by quantities that
grow at most polynomially in \(\ell\); in particular
\((2+\lambda_{\rm rec})^{p-2}\le C(1+\ell)^{d(p-2)/2}\). Hence
\(\beta_v\le C(1+\ell)^{d(p+2)/2}\theta_0^{n_\ell}\) with \(C\) independent of
\(\ell\) and \(v_H\), and \(\beta_v\le1/2\) for large \(\ell\). The weights
computed with \(\theta_0\ge\theta\) are admissible in the proof of
Theorem~\ref{thm:p-computed}, since the interior decay \eqref{eq:p-interior-decay} holds with every factor
at least the true one.
\end{proof}

The logarithmic growth comes from converting the energy decay into the
pointwise quantity \(\lambda_C\) by an inverse estimate. If the shape
weight \(|\nabla U|^{p-2}\) of the ideal states satisfies the local \(A_2\)
condition with a load-independent constant, Theorem~\ref{thm:p-a2-uniform}
replaces \(\Gamma_0\) and \(\Gamma_{\Pi,0}\) by load-independent constants,
but the factor \(h_{\min}^{-d}\) remains. Whether a depth independent of \(h\)
suffices is open.

\subsection{Anisotropy and the Linearized Tensor}
\label{subsec:p-anisotropy}

The scalar condition \eqref{eq:p-structure} bounds the secant tensor
by multiples of a scalar weight, so the anisotropy of the Jacobian enters
only through \(\hat M/\hat m\). We now compare tensors directionally. In
this subsection we assume in addition that \(x\mapsto a(x,\xi)\) is
constant on every fine element. This holds for laws of the form
\(a(x,\xi)=a_\varepsilon(x)\bar a(\xi)\), and in particular for
Examples~\ref{example:saturated-law} and~\ref{example:quartic-law}. For
\(U\in V_h\) let \(A_U=D_\xi a(x,\nabla U)\), with largest and smallest
eigenvalues \(\lambda_U^+\) and \(\lambda_U^-\). They are then constant on
every fine element. We say that the law has
tensor comparability constants \(0<c_{\rm tc}\le C_{\rm tc}\) and error
constant \(c_e>0\) if
\begin{align}
 c_{\rm tc}D_\xi a(x,\xi)&\le B\le C_{\rm tc}\big(D_\xi a(x,\xi)+a_\varepsilon|\eta-\xi|^{p-2}I\big) \label{eq:p-tensor-comparability}\\
 B&\ge c_ea_\varepsilon|\eta-\xi|^{p-2}I \label{eq:p-tensor-error}
\end{align}
for all \(\xi,\eta\), where \(B=\int_0^1D_\xi a(x,\xi+t(\eta-\xi))\,dt\)
and the inequalities hold for symmetric matrices.

\begin{lem}[Tensor comparability]
\label{lem:p-tensor-comparability}
Every law with \eqref{eq:p-structure} satisfies
\eqref{eq:p-tensor-comparability}--\eqref{eq:p-tensor-error} with
\(c_{\rm tc}=1/\hat\chi_p\), \(C_{\rm tc}=3^{p-2}\max\{\hat M/\hat m,\hat M\}\)
and \(c_e=c_p\hat m\). The quartic law satisfies them with
\(c_{\rm tc}=1/4\), \(C_{\rm tc}=3\) and \(c_e=1/12\).
\end{lem}

\begin{proof}
Write \(\delta=\eta-\xi\) and \(\xi_t=\xi+t\delta\). For the general
statement, \(\omega(\xi,\eta)\ge\max\{(1+|\xi|)^{p-2},|\delta|^{p-2}\}\),
Lemma~\ref{lem:p-secant} and the upper bound in \eqref{eq:p-structure}
give \(B\ge c_p\hat ma_\varepsilon\omega(\xi,\eta)I\). Hence
\(B\ge D_\xi a(x,\xi)/\hat\chi_p\) and
\(B\ge c_p\hat ma_\varepsilon|\delta|^{p-2}I\). Lemma~\ref{lem:p-splitting} with
\(e\) replaced by \(\delta\) gives
\(B\le3^{p-2}\hat Ma_\varepsilon((1+|\xi|)^{p-2}+|\delta|^{p-2})I\), and the
lower bound in \eqref{eq:p-structure},
\(a_\varepsilon(1+|\xi|)^{p-2}I\le D_\xi a(x,\xi)/\hat m\), proves the upper bound.

For the quartic law,
\(B=a_\varepsilon[(1+\overline{|\xi_t|^2})I+2\,\overline{\xi_t\otimes\xi_t}]\),
where the bar denotes the average over \(t\in[0,1]\). For scalars
\(\alpha\) and \(\beta\), the average of \((\alpha+t\beta)^2\) equals
\(\alpha^2+\alpha\beta+\beta^2/3\), which is at least \(\alpha^2/4\) and
at least \(\beta^2/12\). Hence
\(\overline{(\xi_t\cdot\zeta)^2}\ge(\xi\cdot\zeta)^2/4\) for every
\(\zeta\), \(\overline{|\xi_t|^2}\ge\max\{|\xi|^2/4,|\delta|^2/12\}\), and
\(B\ge D_\xi a(x,\xi)/4\) and \(B\ge a_\varepsilon|\delta|^2I/12\). For the
upper bound, \(\overline{|\xi_t|^2}\le2|\xi|^2+|\delta|^2\), and
\(\xi\otimes\delta+\delta\otimes\xi\le\xi\otimes\xi+\delta\otimes\delta\)
gives \(\overline{\xi_t\otimes\xi_t}\le\frac32\xi\otimes\xi+|\delta|^2I\).
Therefore
\(B\le a_\varepsilon[(1+2|\xi|^2)I+3\xi\otimes\xi]+3a_\varepsilon|\delta|^2I\le2D_\xi a(x,\xi)+3a_\varepsilon|\delta|^2I\).
\end{proof}

\begin{definition}[Tensor quantities]
\label{def:p-tensor-quantities}
For a coarse element \(T\), let \(\mathcal P^{\rm ten}_T(U)\) be the smallest
constant with
\begin{equation}
 H^{-2}\int_T\lambda_U^+v_f^2\,dx\le\mathcal P^{\rm ten}_T(U)\int_{N(T)}\nabla v_f^TA_U\nabla v_f\,dx \quad \forall v_f\in V_f \label{eq:p-tensor-poincare}
\end{equation}
For a one-layer patch \(\omega\), let
\(\alpha_\omega(U)=\int_\omega\lambda_U^+\,dx/\int_\omega\lambda_U^-\,dx\). With
maxima over a region \(A\) as in Definition~\ref{def:p-sobolev-ratio}, set
\begin{equation}
 \Gamma^{\rm T}_A(U)=1+\mathcal P^{\rm ten}_A(U)+\alpha_A(U)C_{S,A}(\lambda_U^-)^2R_{\sigma,A}(\lambda_U^-)^{2/\sigma}+\Gamma^e_A \label{eq:p-tensor-gamma}
\end{equation}
\end{definition}

For an isotropic tensor \(A_U=\widehat W_UI\), \(\mathcal P^{\rm ten}_T\) is the
weighted Poincar\'e constant \(\mathcal P_T\) and \(\alpha_\omega=1\). For
an anisotropic tensor, \(\mathcal P^{\rm ten}_T\) compares the largest eigenvalue in
the mass term with the full tensor in the energy.

\begin{prop}[Contraction with the linearized tensor]
\label{prop:p-tensor-contraction}
Assume \eqref{eq:p-structure} and
\eqref{eq:p-tensor-comparability}--\eqref{eq:p-tensor-error}. For the
pair \eqref{eq:localization-pair} with \(D^+=N^3(D^-)\subset D\) and
\(A=D^+\setminus D^-\),
\(\mathscr D_{D^-}(U_D,U)\le\theta_A\mathscr D_{D^+}(U_D,U)\) with
\(\theta_A=\kappa_A/(1+\kappa_A)\) and
\begin{equation}
 \kappa_A^2=C\frac{C_{\rm tc}}{c_*}\Gamma^{\rm T}_A(U),\qquad c_*=\min\{c_{\rm tc},c_e\} \label{eq:p-tensor-kappa}
\end{equation}
where \(C\) depends only on \(p\), \(\sigma\), \(d\) and the mesh family.
\end{prop}

\begin{proof}
By \eqref{eq:p-tensor-comparability}, the \(B\)-energy of \(z_f\) on
\(A\) is at most \(C_{\rm tc}\) times the sum of its \(A_U\)-energy and of
\(\int_Aa_\varepsilon|\nabla e|^{p-2}|\nabla z_f|^2\). The error parts of
the proof of Lemma~\ref{lem:p-cutoff} bound the second term by
\(C(1+\Gamma^e_A)\int_Aa_\varepsilon|\nabla e|^p\). For the first, let \(t\)
be a fine element of a transition element and \(\eta_t\) the value of
\(\eta_H\) at a vertex of \(t\). Then
\(\nabla I_h(\eta_He)=\eta_t\nabla e+\nabla I_h((\eta_H-\eta_t)e)\) on
\(t\), and \(|\nabla I_h((\eta_H-\eta_t)e)|\le CH^{-1}\max_t|e|\). Since
\(|\eta_t|\le1\) and \(A_U\) is constant on \(t\),
\begin{equation}
 \int_t\nabla I_h(\eta_He)^TA_U\nabla I_h(\eta_He)\,dx\le2\int_t\nabla e^TA_U\nabla e\,dx+CH^{-2}\int_t\lambda_U^+e^2\,dx \label{eq:p-tensor-product}
\end{equation}
and \eqref{eq:p-tensor-poincare} bounds the last term. For the projection
correction, with \(\omega=N(K_z)\),
\begin{equation}
 \int_\Omega\nabla\lambda_z^TA_U\nabla\lambda_z\,dx\le CH^{-2}\int_{\operatorname{star}z}\lambda_U^+\,dx\le CH^{-2}\alpha_\omega\int_\omega\lambda_U^-\,dx \label{eq:p-tensor-hat}
\end{equation}
Lemma~\ref{lem:p-moment-shift} and
Lemma~\ref{lem:p-moment-bound} with the weight \(\lambda_U^-\) give
\begin{equation}
 |c_z|^2H^{-2}\int_{\operatorname{star}z}\lambda_U^+\,dx\le C\alpha_\omega C_{S}(\lambda_U^-)^2R_\sigma(K_z,\omega;\lambda_U^-)^{2/\sigma}\int_\omega\lambda_U^-|\nabla e|^2\,dx \label{eq:p-tensor-correction}
\end{equation}
Finally,
\(\int\lambda_U^-|\nabla e|^2\le\int\nabla e^TA_U\nabla e\le c_{\rm tc}^{-1}\int\nabla e^TB\nabla e\)
and \(\int a_\varepsilon|\nabla e|^p\le c_e^{-1}\int\nabla e^TB\nabla e\), by
\eqref{eq:p-tensor-comparability}--\eqref{eq:p-tensor-error}. The argument
of Proposition~\ref{prop:p-contraction} closes with \eqref{eq:p-tensor-kappa}.
\end{proof}

In general the structural prefactor is not smaller than in the scalar
version: with the constants of the first part of
Lemma~\ref{lem:p-tensor-comparability}, \(C_{\rm tc}/c_*\ge\hat\chi_p\).
For the quartic law, \(C_{\rm tc}/c_*=36\) and \(\hat\chi_4=72\); since both
factors are multiplied by unspecified constants, this comparison has limited
meaning. The gain lies in the local quantities: the factor contains the
tensor Poincar\'e constant of the actual Jacobian and an averaged
anisotropy, not the global ratio \(\hat M/\hat m\).
Section~\ref{subsec:p-numerics} shows that \(\mathcal P^{\rm ten}\) separates the
full Jacobian from isotropic coefficients at every contrast. The
amplitude, assembly and energy estimates also have tensor versions. For a
state \(U\) and a correction \(q\), let
\begin{equation}
 \mathcal Q^{\rm T}_D(q;U)=\int_D\nabla q^TA_U\nabla q\,dx+\Big(\frac{C_{\rm tc}}{c_e}\Big)^{(p-2)/2}\int_Da_\varepsilon|\nabla q|^p\,dx \label{eq:p-tensor-amplitude-quantity}
\end{equation}

\begin{lem}[Tensor amplitude and assembly]
\label{lem:p-tensor-amplitude}
Assume \eqref{eq:p-structure} and
\eqref{eq:p-tensor-comparability}--\eqref{eq:p-tensor-error}. For the pair
\eqref{eq:localization-pair} on a patch \(D\), with \(q=q(v_H)\),
\begin{equation}
 \mathscr D_D(U_D,U)\le C_{\rm amp}^{\rm T}C_{\rm tc}\mathcal Q^{\rm T}_D(q;U) \label{eq:p-tensor-amplitude}
\end{equation}
with \(C_{\rm amp}^{\rm T}=2\max\{1,(2/p)(2(p-2)/p)^{(p-2)/2}\}\) for \(p>2\)
and \(C_{\rm amp}^{\rm T}=1\) for \(p=2\). For \(U\in V_h\) and
\(w_i\in V_f\),
\begin{equation}
 \int_\Omega\nabla\delta^TA_U\nabla\delta\,dx\le C\,\Gamma^{\rm T}_\Omega(U)\sum_i\int_{N(\omega_i)}\nabla w_i^TA_U\nabla w_i\,dx,\qquad \delta=(I-\Pi_H)I_h\sum_i\phi_iw_i \label{eq:p-tensor-assembly}
\end{equation}
with \(C\) depending only on \(p\), \(d\), the Sobolev exponent in
\(\Gamma^{\rm T}_\Omega\) and the mesh family.
\end{lem}

\begin{proof}
By \eqref{eq:outer-pythagoras} and \eqref{eq:p-tensor-comparability},
\(\mathscr D_D(U_D,U)\le\int_D\nabla q^TB\nabla q\le C_{\rm tc}\int_D(\nabla q^TA_U\nabla q+a_\varepsilon|\nabla e|^{p-2}|\nabla q|^2)\).
By \eqref{eq:p-tensor-error},
\(\int_Da_\varepsilon|\nabla e|^p\le c_e^{-1}\int_D\nabla e^TB\nabla e=c_e^{-1}\mathscr D_D(U_D,U)\).
Young's inequality \eqref{eq:p-young} with
\(\epsilon^{p/(p-2)}=pc_e/(2(p-2)C_{\rm tc})\) absorbs half of the natural
distance for \(p>2\), as in Lemma~\ref{lem:p-amplitude}, and leaves the factor
\(\epsilon^{-p/2}=(2(p-2)C_{\rm tc}/(pc_e))^{(p-2)/2}\) on
\(\int_Da_\varepsilon|\nabla q|^p\). This proves \eqref{eq:p-tensor-amplitude}.

For \eqref{eq:p-tensor-assembly}, argue as in Lemma~\ref{lem:p-assembly}
with the tensor energy. On a fine element \(t\subset\omega_i\), let
\(\phi_t\) be the value of \(\phi_i\) at a vertex of \(t\). As in
\eqref{eq:p-tensor-product},
\(\int_t\nabla I_h(\phi_iw)^TA_U\nabla I_h(\phi_iw)\le2\int_t\nabla w^TA_U\nabla w+CH^{-2}\int_t\lambda_U^+w^2\),
and \eqref{eq:p-tensor-poincare} on \(T\subset\omega_i\), with \(N(T)\subset N(\omega_i)\),
bounds the mass term. For the projection, Lemma~\ref{lem:p-moment-shift}
with \(\vartheta=\phi_i\) and \eqref{eq:p-tensor-correction} on
\(\omega=N(K_z)\subset N(\omega_i)\) bound each coefficient term by
\(C\alpha_\omega C_S(\lambda_U^-)^2R_\sigma(K_z,\omega;\lambda_U^-)^{2/\sigma}\int_\omega\nabla w^TA_U\nabla w\).
Bounded overlap completes the proof.
\end{proof}

\begin{thm}[Energy error with the linearized tensor]
\label{thm:p-tensor-energy}
Assume \eqref{eq:p-structure},
\eqref{eq:p-tensor-comparability}--\eqref{eq:p-tensor-error} and the
hypotheses of Theorem~\ref{thm:p-energy}, with the weights \(t_i\) computed
from \(\theta=\kappa/(1+\kappa)\) and
\(\kappa^2=C(C_{\rm tc}/c_*)\Gamma^{\rm T}_\Omega(u_h)\) as in
\eqref{eq:p-tensor-kappa}. Every global minimizer \(u_{\rm ms}\) of
\(\mathcal E\) on \(M_h^{\mathcal D}(V_H)\) satisfies
\begin{equation}
 0\le\mathcal E(u_{\rm ms})-\mathcal E(u_h)\le C\frac{C_{\rm tc}^2}{c_*}\Gamma^{\rm T}_\Omega(u_h)\mathcal N_t(\mathcal D)\mathcal Q^{\rm T}_\Omega(q^*;u_h) \label{eq:p-tensor-energy}
\end{equation}
with \(C\) depending only on \(p\), \(d\), the Sobolev exponent and the mesh family, and
\(\mathcal Q^{\rm T}_\Omega(q^*;u_h)\le C\Gamma_\Pi(u_h)(\hat M+(C_{\rm tc}/c_e)^{(p-2)/2})\hat L_F^2/\hat m^2\).
\end{thm}

\begin{proof}
As in the proof of Theorem~\ref{thm:p-energy},
\(\mathcal E(u_{\rm ms})-\mathcal E(u_h)\le\mathscr D_\Omega(w_*,u_h)=\int_\Omega\nabla\delta_*^TB\nabla\delta_*\)
with \(\delta_*=(I-\Pi_H)I_h\sum_i\phi_ie_i\). By
\eqref{eq:p-tensor-comparability} with \(\xi=\nabla u_h\), this is at most
\(C_{\rm tc}(\int\nabla\delta_*^TA_{u_h}\nabla\delta_*+\int a_\varepsilon|\nabla\delta_*|^p)\).
Lemma~\ref{lem:p-tensor-amplitude} and the proof of
\eqref{eq:p-assembly-lp}, whose factor is at most \(C(1+\Gamma^e_\Omega)\),
bound the two terms by \(C\Gamma^{\rm T}_\Omega(u_h)\) times
\(\sum_i\int_{N(\omega_i)}(\nabla e_i^TA_{u_h}\nabla e_i+a_\varepsilon|\nabla e_i|^p)\le2c_*^{-1}\sum_i\mathscr D_{N(\omega_i)}(U_i,u_h)\),
by \eqref{eq:p-tensor-comparability}--\eqref{eq:p-tensor-error}.
Proposition~\ref{prop:p-tensor-contraction} gives
\(\mathscr D_{N(\omega_i)}(U_i,u_h)\le t_i\mathscr D_{D_i}(U_i,u_h)\), and
\eqref{eq:p-tensor-amplitude} and the weighted overlap complete the proof of
\eqref{eq:p-tensor-energy}. For the amplitude,
\(A_{u_h}\le\hat M\widehat W_{u_h}I\) and the proof of
\eqref{eq:p-amplitude-data} give the stated bound.
\end{proof}

All factors in \eqref{eq:p-tensor-energy} are local quantities of the
Jacobian at the fine solution, apart from the structural constants
\(C_{\rm tc}\), \(c_e\) and \(c_*\).


\section{Numerical Experiments for the Quartic Law}
\label{subsec:p-numerics}

We now use the quartic law of Example~\ref{example:quartic-law},
\(a=a_\varepsilon(1+|\xi|^2)\xi\), with \(p=4\), \(\hat m=1/2\) and
\(\hat M=3\), and the saturated law of Example~\ref{example:saturated-law} as a
bounded-Hessian control. Unless
stated otherwise, the coefficient is the pixel field of
Section~\ref{subsec:numerical-setup} with contrast 10, \(n_H=12\) and
\(n_h=48\), and depths zero to four contain no full-domain patches. The
channel coefficient equals the contrast on the six horizontal strips
\(j/6+1/24\le x_2<j/6+1/12\), \(j=0,\ldots,5\), which span the domain, and
one elsewhere. For each load we compute comparator reconstructions
\(M_h^\ell(\Pi_Hu_h)\) and report the relative gradient error \(e_4\) at
depth four and the geometric mean \(\bar r=(e_4/e_0)^{1/4}\) of the error
ratios between successive depths. The weighted quantities are evaluated
at the quartic fine solution for
\(\widehat W=a_\varepsilon(1+|\nabla u_h|)^2\) and for its scale-free part
\(G=a_\varepsilon|\nabla u_h|^2\), and for the saturated law for
\(\widehat W=a_\varepsilon\). The Poincar\'e constants are computed from local
generalized eigenvalue problems on \(N(T)\) with the moment constraints
of the interior vertices of \(T\). This local space contains the restrictions of
\(V_f\), so the values are upper bounds for \(\mathcal P_T\). The ratios
are maxima of \(R_\sigma(K_z,N(K_z);\cdot)\) over interior vertices.

\paragraph{Load.}
\label{par:p-load}
Table~\ref{tab:p-load} increases the load by six orders of magnitude.
The quartic decay slows mildly up to load \(2\times10^3\) and then
saturates; the depth-four errors at loads \(2\times10^4\), \(2\times10^5\)
and \(2\times10^6\) agree to within one percent. The saturated law shows no
trend. Over the same range, the local \(p\)-energy density \(\max\Lambda_T\)
grows by a factor of more than \(6\times10^4\), whereas \(\max\mathcal P_T\)
and the ratios of \(\widehat W\) stay between \(0.21\) and \(0.24\) and
between \(5.8\) and \(8.2\). The shape ratios of \(G\) decrease towards
the ratios of \(\widehat W\), as Corollary~\ref{cor:p-load-shape} suggests.
They are large only at small loads, where \(\widehat W\approx a_\varepsilon\).
The crude bound of Lemma~\ref{lem:p-crude} is therefore far from the
observed behaviour, while the scale-invariant quantities describe it.
Refining the fine mesh at load \(2\times10^4\) to \(n_h=96\) and \(192\)
does not increase the ratios of \(W=(1+|\nabla u_h|)^2\): \(R_2(W)\)
equals \(4.43\), \(4.28\) and \(4.22\) for \(n_h=48\), \(96\) and \(192\).
This agrees with the flux-balanced profile of
Section~\ref{subsec:p-model-weights} rather than with nondegenerate
critical points.

\begin{table}[htbp]
\centering
\small
\begin{tabular}{rcccccccccc}
\toprule
 & \multicolumn{2}{c}{Quartic} & \multicolumn{2}{c}{Saturated} & & \multicolumn{3}{c}{\(\widehat W\)} & \multicolumn{2}{c}{\(G\)} \\
\cmidrule(lr){2-3}\cmidrule(lr){4-5}\cmidrule(lr){7-9}\cmidrule(lr){10-11}
Load & \(e_4\) & \(\bar r\) & \(e_4\) & \(\bar r\) & \(\max\Lambda_T\) & \(\max\mathcal P_T\) & \(R_2\) & \(R_4\) & \(R_2\) & \(R_4\) \\
\midrule
\(2\) & \(2.68\times10^{-4}\) & \(0.206\) & \(2.57\times10^{-4}\) & \(0.206\) & \(0.104\) & \(0.232\) & \(8.11\) & \(8.11\) & \(53.0\) & \(22.2\) \\
\(2\times10^{1}\) & \(2.33\times10^{-4}\) & \(0.209\) & \(2.60\times10^{-4}\) & \(0.214\) & \(2.10\) & \(0.220\) & \(6.65\) & \(6.62\) & \(14.5\) & \(8.45\) \\
\(2\times10^{2}\) & \(2.94\times10^{-4}\) & \(0.226\) & \(3.18\times10^{-4}\) & \(0.218\) & \(13.3\) & \(0.219\) & \(6.09\) & \(6.07\) & \(7.30\) & \(5.87\) \\
\(2\times10^{3}\) & \(3.29\times10^{-4}\) & \(0.234\) & \(2.78\times10^{-4}\) & \(0.207\) & \(65.3\) & \(0.218\) & \(5.98\) & \(5.96\) & \(6.67\) & \(5.85\) \\
\(2\times10^{4}\) & \(3.37\times10^{-4}\) & \(0.236\) & \(2.74\times10^{-4}\) & \(0.206\) & \(307\) & \(0.217\) & \(5.92\) & \(5.90\) & \(6.61\) & \(5.84\) \\
\(2\times10^{5}\) & \(3.39\times10^{-4}\) & \(0.236\) & \(2.74\times10^{-4}\) & \(0.206\) & \(1430\) & \(0.217\) & \(5.97\) & \(5.87\) & \(6.59\) & \(5.84\) \\
\(2\times10^{6}\) & \(3.40\times10^{-4}\) & \(0.236\) & \(2.74\times10^{-4}\) & \(0.206\) & \(6630\) & \(0.217\) & \(6.27\) & \(5.86\) & \(6.59\) & \(5.84\) \\
\bottomrule
\end{tabular}
\caption{Comparator localization and weighted quantities versus load,
 with the pixel coefficient of contrast 10, \(n_H=12\) and \(n_h=48\).
 The weighted quantities are those of the quartic fine solution, with
 \(\widehat W=a_\varepsilon(1+|\nabla u_h|)^2\) and \(G=a_\varepsilon|\nabla u_h|^2\).
 For the coefficient alone, \(\max\mathcal P_T=0.239\) and \(R_2=R_4=8.62\).}
\label{tab:p-load}
\end{table}

\paragraph{Contrast and Channels.}
\label{par:p-contrast}
Table~\ref{tab:p-contrast} varies the contrast and the coefficient
geometry at a moderate and a large load. For the pixel field, the
quantities of \(\widehat W\) are nearly load-independent and close to those of
\(a_\varepsilon\), and the decay slows with the contrast as these quantities
grow. At contrast 1, \(\bar r\) equals \(0.189\) and \(0.200\) for the
quartic law at loads \(2\) and \(2\times10^4\). The channels connect the
two Dirichlet edges through the high-conductivity phase, and they decay
markedly more slowly than the pixel field at the same contrast. At
contrast 100 their Poincar\'e constant is also larger, but at contrast 10 it
is smaller than for the pixel field. The local quantities, which enter the
contraction through their maxima, therefore bound the decay but do not
rank coefficient geometries. With \(n_H=24\), the coarse mesh resolves the
channel width, \(\max\mathcal P_T(a_\varepsilon)\) drops to \(0.31\), and the
decay ratios at contrast 100 and load \(2\times10^4\) improve to \(0.363\)
and \(0.280\) for the quartic and saturated laws.

\begin{table}[htbp]
\centering
\small
\begin{tabular}{llrcccccc}
\toprule
 & & & \multicolumn{2}{c}{\(\bar r\)} & \multicolumn{2}{c}{\(\widehat W\), quartic} & \multicolumn{2}{c}{\(a_\varepsilon\)} \\
\cmidrule(lr){4-5}\cmidrule(lr){6-7}\cmidrule(lr){8-9}
Coefficient & Contrast & Load & Quartic & Saturated & \(\max\mathcal P_T\) & \(R_2\) & \(\max\mathcal P_T\) & \(R_2\) \\
\midrule
Pixels & 10 & \(2\times10^1\) & \(0.209\) & \(0.214\) & \(0.220\) & \(6.6\) & \(0.239\) & \(8.6\) \\
Pixels & 10 & \(2\times10^4\) & \(0.236\) & \(0.206\) & \(0.217\) & \(5.9\) & \(0.239\) & \(8.6\) \\
Pixels & 100 & \(2\times10^1\) & \(0.297\) & \(0.299\) & \(0.388\) & \(71.8\) & \(0.404\) & \(84.8\) \\
Pixels & 100 & \(2\times10^4\) & \(0.348\) & \(0.292\) & \(0.401\) & \(87.6\) & \(0.404\) & \(84.8\) \\
Channels & 10 & \(2\times10^1\) & \(0.324\) & \(0.315\) & \(0.172\) & \(4.1\) & \(0.167\) & \(3.8\) \\
Channels & 10 & \(2\times10^4\) & \(0.354\) & \(0.277\) & \(0.175\) & \(6.8\) & \(0.167\) & \(3.8\) \\
Channels & 100 & \(2\times10^1\) & \(0.422\) & \(0.395\) & \(0.751\) & \(22.2\) & \(1.035\) & \(31.5\) \\
Channels & 100 & \(2\times10^4\) & \(0.407\) & \(0.417\) & \(0.622\) & \(36.5\) & \(1.035\) & \(31.5\) \\
\bottomrule
\end{tabular}
\caption{Geometric-mean decay ratio \(\bar r\) of comparator errors over
 depths zero to four and weighted quantities, with \(n_H=12\) and
 \(n_h=48\). For the saturated law, \(\widehat W=a_\varepsilon\).}
\label{tab:p-contrast}
\end{table}

\paragraph{Coarse and Fine Refinement.}
\label{par:p-refinement}
Table~\ref{tab:p-refinement} refines the coarse and the fine mesh with
the pixel field fixed. On the \(6/24\) mesh, depth-four patches cover
almost the whole domain, and on the \(12/48\) mesh they still meet the
boundary on most sides. The Dirichlet condition then reduces the errors.
On the \(24/96\) mesh, the mean ratio over depths zero to four agrees with
the ratio \((e_6/e_4)^{1/2}\) over depths four to six to within \(0.02\):
\(0.234\) and \(0.228\), \(0.284\) and \(0.301\), \(0.242\) and \(0.240\),
and \(0.220\) and \(0.219\) for the four columns. In three columns the
interior rate is therefore reached, and the increase of \(\bar r\) from the
\(6/24\) to the \(24/96\) mesh is a boundary effect. For the quartic law at
load \(2\times10^4\) the ratio is still increasing, and its load
dependence, \(0.234\) to \(0.284\), is about twice that on the \(12/48\)
mesh. The load study above therefore understates the load dependence of
the interior rate. From \(n_H=12\) on, the quartic decay
at load \(2\times10^4\) is the slowest in every row, in agreement with the
anisotropy study below. Refining the fine mesh at \(n_H=12\) changes
\(\bar r\) by at most \(0.019\) from \(n_h=48\) to \(96\) and by at most
\(0.005\) from \(96\) to \(192\). The depth-four errors of the quartic law
grow by \(64\%\) and \(56\%\) for the two loads from \(n_h=48\) to \(192\),
by \(48\%\) and \(42\%\) in the first refinement and by about \(10\%\) in
the second. The decay ratios and the slowing growth of the errors are
consistent with constants independent of the fine-to-coarse mesh ratio.

\begin{table}[htbp]
\centering
\small
\begin{tabular}{rrcccccccc}
\toprule
 & & \multicolumn{4}{c}{Quartic} & \multicolumn{4}{c}{Saturated} \\
\cmidrule(lr){3-6}\cmidrule(lr){7-10}
 & & \multicolumn{2}{c}{Load \(2\times10^1\)} & \multicolumn{2}{c}{Load \(2\times10^4\)} & \multicolumn{2}{c}{Load \(2\times10^1\)} & \multicolumn{2}{c}{Load \(2\times10^4\)} \\
\cmidrule(lr){3-4}\cmidrule(lr){5-6}\cmidrule(lr){7-8}\cmidrule(lr){9-10}
\(n_H\) & \(n_h\) & \(e_4\) & \(\bar r\) & \(e_4\) & \(\bar r\) & \(e_4\) & \(\bar r\) & \(e_4\) & \(\bar r\) \\
\midrule
6 & 24 & \(1.28\times10^{-4}\) & \(0.167\) & \(1.21\times10^{-4}\) & \(0.171\) & \(1.31\times10^{-4}\) & \(0.168\) & \(2.51\times10^{-4}\) & \(0.188\) \\
12 & 48 & \(2.33\times10^{-4}\) & \(0.209\) & \(3.37\times10^{-4}\) & \(0.236\) & \(2.60\times10^{-4}\) & \(0.214\) & \(2.74\times10^{-4}\) & \(0.206\) \\
24 & 96 & \(2.09\times10^{-4}\) & \(0.234\) & \(4.08\times10^{-4}\) & \(0.284\) & \(2.40\times10^{-4}\) & \(0.242\) & \(2.21\times10^{-4}\) & \(0.220\) \\
12 & 96 & \(3.46\times10^{-4}\) & \(0.228\) & \(4.79\times10^{-4}\) & \(0.255\) & \(3.80\times10^{-4}\) & \(0.233\) & \(3.94\times10^{-4}\) & \(0.222\) \\
12 & 192 & \(3.82\times10^{-4}\) & \(0.232\) & \(5.26\times10^{-4}\) & \(0.259\) & & & & \\
\bottomrule
\end{tabular}
\caption{Comparator localization under coarse and fine refinement with the
 pixel coefficient of contrast 10. The relative errors are measured
 against the fine solution on the same fine mesh.}
\label{tab:p-refinement}
\end{table}

\paragraph{Frozen Linearization and Anisotropy.}
\label{par:p-frozen}
To isolate the source of the load dependence, we freeze the Jacobian of
the quartic law at the fine solution for load \(2\times10^4\). We then run
the same comparator test for linear problems with the full tensor
\(A=a_\varepsilon[(1+|\nabla u_h|^2)I+2\nabla u_h\otimes\nabla u_h]\), with the
isotropic tensors \(a_\varepsilon(1+k|\nabla u_h|^2)I\), \(k=1,2,3\), and with
\(a_\varepsilon I\). Table~\ref{tab:p-frozen} also reports the tensor Poincar\'e
constant \eqref{eq:p-tensor-poincare}. The frozen full Jacobian reproduces the
nonlinear decay to within \(0.014\) in \(\bar r\), so localization of the
nonlinear reconstruction is governed by the linearization at the fine
solution. The isotropic coefficients decay slightly faster than
\(a_\varepsilon I\) for every \(k\), so the size and variation of the scalar
weight do not cause the slowdown. The full tensor decays more slowly at every
contrast, and its tensor Poincar\'e constant is about twice the isotropic
one. The slowdown is thus an effect of the anisotropy of the Jacobian, and
\(\mathcal P^{\rm ten}_T\) detects it. It does not resolve the small difference
between the isotropic coefficients and \(a_\varepsilon I\) at contrast 1.

\begin{table}[htbp]
\centering
\small
\begin{tabular}{rccccccc}
\toprule
 & \multicolumn{4}{c}{\(\bar r\)} & \multicolumn{3}{c}{\(\max\mathcal P^{\rm ten}_T\)} \\
\cmidrule(lr){2-5}\cmidrule(lr){6-8}
Contrast & Nonlinear & Full \(A\) & Isotropic & \(a_\varepsilon I\) & Full \(A\) & Isotropic & \(a_\varepsilon I\) \\
\midrule
1 & \(0.200\) & \(0.201\) & \(0.182\) & \(0.189\) & \(0.162\) & \(0.083\) & \(0.077\) \\
10 & \(0.236\) & \(0.235\) & \(0.192\) & \(0.206\) & \(0.423\) & \(0.217\) & \(0.239\) \\
100 & \(0.348\) & \(0.362\) & \(0.283\) & \(0.292\) & \(0.807\) & \(0.404\) & \(0.404\) \\
\bottomrule
\end{tabular}
\caption{Decay ratio \(\bar r\) for the quartic law at load \(2\times10^4\)
 and for linear problems with coefficients frozen at its fine solution,
 with \(n_H=12\) and \(n_h=48\). The isotropic decay ratios agree to three
 digits for \(k=1,2,3\); the tensor constants are given for \(k=2\). The
 largest averaged anisotropy of the full tensor is \(2.99\).}
\label{tab:p-frozen}
\end{table}

\paragraph{Compatibility and Computed States.}
\label{par:p-compatibility}
Table~\ref{tab:p-compatibility} reports the compatibility quantities
\eqref{eq:p-compatibility-pointwise}--\eqref{eq:p-compatibility-average}
at the comparator coordinate, where the ideal state is \(u_h\). The
pointwise quantity on the active supports, \(\lambda_C\), decreases with
the depth like the localization error and is below one from depth three
on at load \(2\times10^4\). The pointwise ratio \(\lambda_D\) over the whole
patch does not decrease, because near the artificial boundary the patch
state differs from the ideal state by the fine correction. The averaged
ratio \(\varrho_D\) stays below \(1.9\). This is why
Lemma~\ref{lem:p-tangent} uses \(\lambda_C\) and \(\varrho_D\) rather than
pointwise compatibility on the patch.

\begin{table}[htbp]
\centering
\small
\begin{tabular}{rcccccccc}
\toprule
 & \multicolumn{4}{c}{Load \(2\times10^1\)} & \multicolumn{4}{c}{Load \(2\times10^4\)} \\
\cmidrule(lr){2-5}\cmidrule(lr){6-9}
Depth & \(\lambda_C\) & \(\lambda_D\) & \(\lambda_{\rm rec}\) & \(\varrho_D\) & \(\lambda_C\) & \(\lambda_D\) & \(\lambda_{\rm rec}\) & \(\varrho_D\) \\
\midrule
1 & \(0.602\) & \(0.602\) & \(0.095\) & \(1.37\) & \(2.55\) & \(2.55\) & \(0.715\) & \(1.83\) \\
2 & \(0.284\) & \(0.602\) & \(0.021\) & \(1.34\) & \(1.04\) & \(2.54\) & \(0.190\) & \(1.73\) \\
3 & \(0.073\) & \(0.602\) & \(0.004\) & \(1.34\) & \(0.384\) & \(2.54\) & \(0.069\) & \(1.73\) \\
4 & \(0.013\) & \(0.602\) & \(0.001\) & \(1.34\) & \(0.076\) & \(2.54\) & \(0.005\) & \(1.73\) \\
\bottomrule
\end{tabular}
\caption{Compatibility quantities at the comparator coordinate for the
 quartic law with the pixel coefficient of contrast 10, \(n_H=12\) and
 \(n_h=48\). The quantity \(\lambda_D\) is the pointwise ratio over the
 whole patch \(D_i\), which the theory does not need.}
\label{tab:p-compatibility}
\end{table}

We also solved the reduced problem for the quartic law at load
\(2\times10^4\) by BFGS, starting from zero and from \(\Pi_Hu_h\). At depth
two, both runs reach the same stationary state, with relative reduced
residuals \(2.9\times10^{-11}\) and \(5.2\times10^{-8}\). Both end with the
optimizer's precision-loss flag. The computed-state error
\(7.576\times10^{-3}\) agrees with the comparator error
\(7.571\times10^{-3}\). At the computed coordinate, \(\lambda_C=1.03\),
\(\varrho_D=1.73\) and \(\lambda_{\rm rec}=0.19\), so the hypothesis \(\lambda_C<1\) of Theorem~\ref{thm:p-computed} is narrowly violated although the computed state is
accurate. At depth three, both starts again reach the same stationary
state, with relative reduced residuals \(4.4\times10^{-11}\) and
\(7.3\times10^{-8}\). The computed-state error \(1.840\times10^{-3}\) agrees
with the comparator error \(1.839\times10^{-3}\), and \(\lambda_C=0.38\),
\(\varrho_D=1.73\) and \(\lambda_{\rm rec}=0.069\) at the computed
coordinate. These checks do not evaluate the constants in
\eqref{eq:p-oversampling} and do not establish uniqueness of the
stationary state.

\paragraph{Limitations.}
\label{par:p-limitations}
These computations are two-dimensional and use one pixel field and one
channel geometry. Apart from the outer solves above, they use comparator
reconstructions. They do not evaluate \(C_{{\rm cut},p}\), \(C_E\) or the
Sobolev constants of Definition~\ref{def:p-sobolev-ratio}, and the
theoretical factor \(\theta\) of \eqref{eq:p-uniform-factor} is much closer
to one than the observed ratios. The \(24/96\) and \(12/192\) runs use a
patch-by-patch evaluation of the same reconstruction. It reproduces the
errors of the original implementation exactly on the \(6/24\) mesh and for
the quartic law at load \(2\times10^4\) on the \(12/48\) mesh.


\section{Discussion of the \texorpdfstring{{\boldmath$p$}}{p}-Structure Results}
\label{sec:p-discussion}

For the quartic law, the largest eigenvalue of the Jacobian
\eqref{eq:quartic-jacobian} grows with \(|\xi|\), so the state-independent
upper constant in \eqref{eq:hessian-bounds} is unavailable.
Section~\ref{sec:p-structure} replaces the global ratio \(M/m\) by
scale-invariant local quantities of the weight
\(a_\varepsilon(1+|\nabla u_h|)^{p-2}\), which control the cutoff and the
outer amplitudes together. The price is that the constants are
quantities of the fine solution. For \(d=2\), and for \(d=3\) with
\(p\le6\), they are bounded independently of the fine mesh, but the crude
bound of Lemma~\ref{lem:p-crude} grows with the load and with a negative
power of \(H\), so it does not by itself give logarithmic oversampling.
Theorem~\ref{thm:p-a2-uniform} reduces load- and mesh-uniform bounds,
including the projection constant \(\Gamma_\Pi\), to a local Muckenhoupt
\(A_2\) condition on the shape weight \(|\nabla u_h|^{p-2}\), and
Theorem~\ref{thm:p-fixed-mesh} proves load-uniform bounds on every fixed
mesh under a nondegeneracy condition on the large-load limit. Whether the
shape weight of the fine solution satisfies the \(A_2\) condition uniformly
in the load and the mesh remains open. This is a question of gradient
regularity for problems with rough coefficients. Corollary~\ref{cor:p-model-a2}
covers the model profiles suggested by flux balance at extrema where the
load is positive. Nondegenerate critical points with \(p-2\ge d\) violate
the condition, and saddle points, regions without load, coefficient
interfaces and boundary points where the gradient vanishes, such as convex
corners of the domain, are not analysed. The computations show saturation of these
quantities and of the decay.

The compatibility condition \(\lambda_C<1\) in Theorem~\ref{thm:p-computed}
follows from oversampling for coordinates with bounded energy
(Corollary~\ref{cor:p-apriori}). For fixed \(H\), the depth grows like
\(\log(1/h)\). Through the crude constants, it also grows with the load and
like a power of \(1/H\). For the quartic law in two dimensions, the patch
diameter it requires therefore need not shrink as \(H\to0\). Whether a
depth independent of the fine mesh and the load suffices is open. The tensor estimates of Section~\ref{subsec:p-anisotropy} improve
the local quantities but not the structural prefactor; a contraction
argument that keeps the direction of the cutoff might relate the tensor
Poincar\'e constant quantitatively to the observed decay. Degenerate laws such as the \(p\)-Laplacian, for
which \((1+|\xi|)^{p-2}\) is replaced by \(|\xi|^{p-2}\), lose the lower
bound at critical points and require coercivity in the \(p\)-energy. Laws
with \(p<2\) are not covered.

\paragraph{Use of Generative AI.}
The use of generative AI in this work is described in the main paper.

\end{document}